\documentclass[draft]{amsart}

\usepackage{amssymb}
\usepackage{url}
\usepackage{amscd}

\theoremstyle{plain}

\theoremstyle{remark}

\begin{document}

\date{} 

\title[Calabi-Yau structure and Special Lagrangian submanifold]
{Calabi-Yau structure and \\
special Lagrangian submanifold \\
of complexified symmeric space}

\author{Naoyuki Koike}

\address{Department of Mathematics, Faculty of Science\\
Tokyo University of Science, 1-3 Kagurazaka\\
Shinjuku-ku, Tokyo 162-8601 Japan}

\email{koike@rs.kagu.tus.ac.jp}

\thanks{}

\subjclass{53D12, 53C35}

\begin{abstract}
It is known that there exist Calabi-Yau structures on the complexifications of symmetric spaces 
of compact type.  In this paper, we describe the Calabi-Yau structures of the complexified symmetric spaces 
in terms of the Schwarz's theorem in detail.  We consider the case where the Calabi-Yau structure arises from 
the Riemannian metric corresponding to the Stenzel metric.  In the complexified symmetric spaces equipped with 
such a Calabi-Yau structure, we give constructions of special Lagrangian submanifolds of any given phase 
which are invariant under the actions of symmetric subgroups of the isometry group of the original symmetric space 
of compact type.  
\end{abstract}

\maketitle


\newcommand\sfrac[2]{{#1/#2}}

\newcommand\cont{\operatorname{cont}}
\newcommand\diff{\operatorname{diff}}


\section{Introduction}
An $2n$-dimensional Riemannian manifold is called a {\it Calabi-Yau manifold} if the holonomy group is 
a subgroup of $SU(n)$.  A Kaehler manifold is Calabi-Yau if and only if it is Ricci-flat.  
Let $(M,J,\omega)$ be a complex $n$-dimensional Kaehler manifold, where $J$ is the complex structure and 
$\omega$ is the Kaehler form.  Also, let $g$ be the Kaehler metric associated to $(J,\omega)$.  
If there exists a non-vanishing holomorphic $(n,0)$-form $\Omega$ on $M$ (i.e., 
the holomorphic complex line bundle $\bigwedge_h^{(n,0)}(M)$ is trivial), then $(M,J,\omega)$ is called a 
{\it almost Calabi-Yau manifold}.  In particular, if $(\omega,\Omega)$ satisfies 
$$\omega^n=(-1)^{n(n-1)/2}(\sqrt{-1})^nc(\Omega\wedge\overline{\Omega})$$
for some positive real constant $c$, then $(M,J,\omega)$ is Ricci-flat and hence it is Calabi-Yau.  
By replacing $\Omega$ to a suitable positive real constant-multiple of $\Omega$ if necessary, 
we may assume that $c=\frac{n!}{2^n}$.  
In the sequel, the Calabi-Yau manifold (resp. the Calabi-Yau structure) means a quadruple $(M,J,\omega,\Omega)$ 
(resp. a triple $(J,\omega,\Omega)$) such that $(J,\omega)$ is a Kaehler structure and that $(\omega,\Omega)$ 
satisfies 
$$\omega^n=(-1)^{n(n-1)/2}n!\left(\frac{\sqrt{-1}}{2}\right)^n\Omega\wedge\overline{\Omega}.\leqno{(1.1)}$$
Let $(J,\omega,\Omega)$ be a Calabi-Yau structure on $M$ and $g$ the Kaehler metric associated to 
$(J,\omega)$.  Then, for any real constant $\theta$, a $n$-form ${\rm Re}(e^{\sqrt{-1}\theta}\Omega)$ is 
a calibration on $(M,g)$.  A submanifold calibrated by ${\rm Re}(e^{\sqrt{-1}\theta}\Omega)$ is called a 
{\it special Lagrangian submanifold of phase} $\theta$.  
According to Strominger-Yau-Zaslov's conjecture (\cite{SYZ}) for the Mirror symmetry in the string theory, 
it is important to construct special Lagrangian submanifolds in a Calabi-Yau manifold.  

Let $M$ be $C^{\omega}$-Riemannian manifold and $M^{\Bbb C}$ its complexification.  
In 1991, V. Gillemin and M. Stenzel (\cite{GS}) gave a construction of Ricci-flat metrics on a sufficiently small 
tubular neighborhood of $M$ in $M^{\Bbb C}$.  
Let $G/K$ be a (Reimannian) symmetric space of compact type.  The complexification $(G/K)^{\mathbb C}$ of $G/K$ is 
defined as the complexified symmetric space $G^{\Bbb C}/K^{\Bbb C}$ equipped with the $G^{\mathbb C}$-invariant 
anti-Kaehler metric $\beta_A$.  The anti-Kaheler manifold $(G^{\mathbb C}/K^{\mathbb C},\beta_A)$ is called 
an {\it anti-Kaehler symmetric space}.  This space $(G^{\mathbb C}/K^{\mathbb C},\beta_A)$ is identified with 
the tangent bundle $T(G/K)$ of $G/K$ under the one-to-one correspondence 
$\displaystyle{\Psi:T(G/K)\mathop{\longrightarrow}_{\cong}G^{\mathbb C}/K^{\mathbb C}}$ defined by 
$$\Psi(v):={\rm Exp}_p((J_0)_p(v))\quad\,\,(p\in G/K,\,\,v\in T_p(G/K))$$
(see Figure 1), where ${\rm Exp}_p$ denotes the exponential map of $(G^{\mathbb C}/K^{\mathbb C},\beta_A)$ at $p$, 
$J_0$ denotes the natural complex structure of $G^{\mathbb C}/K^{\mathbb C}$ and $v$ is regarded as a tangent vector 
of the submanifold $G\cdot o(\approx G/K)$ ($o=eK^{\mathbb C}$) in $G^{\mathbb C}/K^{\mathbb C}$.  
For each $p\in G/K(\approx G\cdot o)$ set $\Psi_p:=\Psi|_{T_p(G/K)}(={\rm Exp}_p\circ(J_0)_p)$ and 
$(G/K)^d_p:=\Psi(T_p(G/K))$.  Note that $(G/K)^d_p$'s equipped with the (Riemannian) metric induced from $\beta_A$ 
are isometric to the symmetric space $G^d/K$ of non-compact type given as the dual of $G/K$ and they 
are totally geodesic submanifolds in $(G^{\mathbb C}/K^{\mathbb C},\beta_A)$.  

\vspace{0.5truecm}

\centerline{
\unitlength 0.1in
\begin{picture}( 38.4200, 19.1600)( -7.4600,-30.3100)
%
\special{pn 8}%
\special{ar 2092 1892 558 132  6.2831853 6.2831853}%
\special{ar 2092 1892 558 132  0.0000000 3.1415927}%
%
\special{pn 8}%
\special{ar 2092 1892 558 132  3.1415927 3.1764257}%
\special{ar 2092 1892 558 132  3.2809250 3.3157581}%
\special{ar 2092 1892 558 132  3.4202574 3.4550905}%
\special{ar 2092 1892 558 132  3.5595898 3.5944228}%
\special{ar 2092 1892 558 132  3.6989221 3.7337552}%
\special{ar 2092 1892 558 132  3.8382545 3.8730876}%
\special{ar 2092 1892 558 132  3.9775868 4.0124199}%
\special{ar 2092 1892 558 132  4.1169192 4.1517523}%
\special{ar 2092 1892 558 132  4.2562516 4.2910847}%
\special{ar 2092 1892 558 132  4.3955839 4.4304170}%
\special{ar 2092 1892 558 132  4.5349163 4.5697494}%
\special{ar 2092 1892 558 132  4.6742487 4.7090818}%
\special{ar 2092 1892 558 132  4.8135810 4.8484141}%
\special{ar 2092 1892 558 132  4.9529134 4.9877465}%
\special{ar 2092 1892 558 132  5.0922458 5.1270789}%
\special{ar 2092 1892 558 132  5.2315781 5.2664112}%
\special{ar 2092 1892 558 132  5.3709105 5.4057436}%
\special{ar 2092 1892 558 132  5.5102429 5.5450760}%
\special{ar 2092 1892 558 132  5.6495752 5.6844083}%
\special{ar 2092 1892 558 132  5.7889076 5.8237407}%
\special{ar 2092 1892 558 132  5.9282400 5.9630731}%
\special{ar 2092 1892 558 132  6.0675723 6.1024054}%
\special{ar 2092 1892 558 132  6.2069047 6.2417378}%
%
\special{pn 8}%
\special{ar 48 1892 1488 1140  5.5349878 6.2831853}%
\special{ar 48 1892 1488 1140  0.0000000 0.7853982}%
%
\special{pn 8}%
\special{ar 4136 1892 1486 1140  2.3561945 3.8894741}%
%
\special{pn 8}%
\special{ar 3626 2014 1154 1334  2.5895814 3.8660538}%
%
\special{pn 8}%
\special{ar 606 2014 1152 1334  5.5595863 6.2831853}%
\special{ar 606 2014 1152 1334  0.0000000 0.5509426}%
%
\special{pn 20}%
\special{sh 1}%
\special{ar 2474 1988 10 10 0  6.28318530717959E+0000}%
\special{sh 1}%
\special{ar 2474 1988 10 10 0  6.28318530717959E+0000}%
%
\special{pn 20}%
\special{sh 1}%
\special{ar 1748 1996 10 10 0  6.28318530717959E+0000}%
\special{sh 1}%
\special{ar 1748 1996 10 10 0  6.28318530717959E+0000}%
%
\special{pn 8}%
\special{pa 1200 1524}%
\special{pa 1610 1410}%
\special{dt 0.045}%
\special{sh 1}%
\special{pa 1610 1410}%
\special{pa 1540 1408}%
\special{pa 1558 1424}%
\special{pa 1550 1446}%
\special{pa 1610 1410}%
\special{fp}%
%
\special{pn 8}%
\special{pa 3050 1348}%
\special{pa 2724 1216}%
\special{dt 0.045}%
\special{sh 1}%
\special{pa 2724 1216}%
\special{pa 2778 1260}%
\special{pa 2774 1236}%
\special{pa 2794 1222}%
\special{pa 2724 1216}%
\special{fp}%
\put(17.6700,-19.7000){\makebox(0,0)[lb]{$o$}}%
\put(25.1000,-20.3100){\makebox(0,0)[lt]{$p$}}%
%
\special{pn 8}%
\special{pa 1172 2110}%
\special{pa 1590 1952}%
\special{dt 0.045}%
\special{sh 1}%
\special{pa 1590 1952}%
\special{pa 1522 1958}%
\special{pa 1540 1972}%
\special{pa 1536 1994}%
\special{pa 1590 1952}%
\special{fp}%
\put(11.4400,-15.2300){\makebox(0,0)[rt]{$\Psi(T_o(G\cdot o))$}}%
\put(11.1600,-21.0200){\makebox(0,0)[rt]{$G\cdot o=G/K$}}%
\put(30.9600,-13.3900){\makebox(0,0)[lt]{$\Psi(T_p(G\cdot o))$}}%
\put(24.1700,-16.8900){\makebox(0,0)[rb]{$(J_0)_p(v)$}}%
\put(28.2600,-18.5600){\makebox(0,0)[lt]{$v$}}%
\put(18.7800,-28.4700){\makebox(0,0)[lt]{$G^{\mathbb C}/K^{\mathbb C}$}}%
\put(30.1200,-30.3100){\makebox(0,0)[lt]{$\,$}}%
%
\special{pn 13}%
\special{pa 2482 1996}%
\special{pa 2798 1910}%
\special{fp}%
\special{sh 1}%
\special{pa 2798 1910}%
\special{pa 2728 1908}%
\special{pa 2748 1924}%
\special{pa 2740 1946}%
\special{pa 2798 1910}%
\special{fp}%
%
\special{pn 13}%
\special{pa 2464 1996}%
\special{pa 2474 1620}%
\special{fp}%
\special{sh 1}%
\special{pa 2474 1620}%
\special{pa 2452 1686}%
\special{pa 2472 1672}%
\special{pa 2492 1686}%
\special{pa 2474 1620}%
\special{fp}%
%
\special{pn 20}%
\special{sh 1}%
\special{ar 2538 1576 10 10 0  6.28318530717959E+0000}%
\special{sh 1}%
\special{ar 2538 1576 10 10 0  6.28318530717959E+0000}%
%
\special{pn 8}%
\special{pa 2910 1672}%
\special{pa 2548 1576}%
\special{dt 0.045}%
\special{sh 1}%
\special{pa 2548 1576}%
\special{pa 2606 1612}%
\special{pa 2600 1590}%
\special{pa 2618 1574}%
\special{pa 2548 1576}%
\special{fp}%
\put(29.4700,-16.4600){\makebox(0,0)[lt]{$\Psi(v)={\rm Exp}_p((J_0)_p(v))$}}%
\end{picture}%
\hspace{5.25truecm}}

\vspace{0.1truecm}

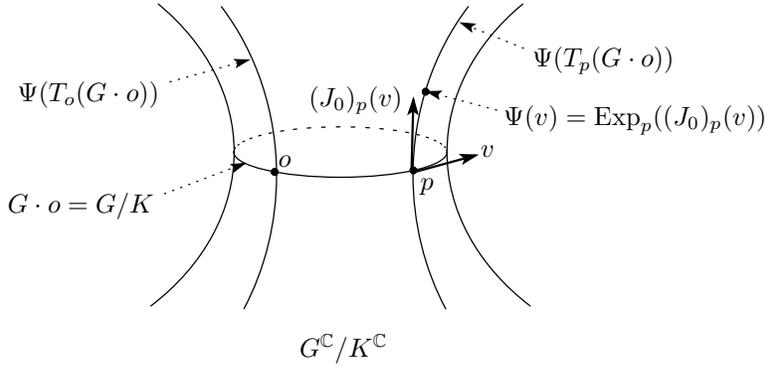
\centerline{{\bf Figure 1.}}

\vspace{0.5truecm}

\noindent
We consider the case where $G/K$ is the sphere $SO(n+1)/SO(n)(=S^n(1))$.  Then the complexification 
$SO(n+1,\mathbb C)/SO(n,\mathbb C)$ of $SO(n+1)/SO(n)$ is embedded into $\mathbb C^{n+1}$ as the complex sphere 
$S^n_{\mathbb C}(1):=\{(z_1,\cdots,z_{n+1})\,|\,\sum\limits_{i=1}^{n+1}z_i^2=1\}$ of complex radius $1$.  
The natural embedding $\iota$ of $SO(n+1,\mathbb C)/SO(n,\mathbb C)$ into $\mathbb C^{n+1}$ is given by 
$$\begin{array}{r}
\displaystyle{\iota(q):=\cosh(\|\Psi^{-1}(q)\|)\cdot\overrightarrow{Op}+\sqrt{-1}\cdot\frac{\sinh(\|\Psi^{-1}(q)\|)}
{\|\Psi^{-1}(q)\|}\cdot\Psi^{-1}(q)}\\
\displaystyle{(q\in SO(n+1,\mathbb C)/SO(\mathbb C)),}
\end{array}$$
where $p$ is the base point of $\Psi^{-1}(q)$, $O$ is the origin of of the $(n+1)$-dimensional Euclidean space 
$\mathbb R^{n+1}$ including $S^n(1)(=SO(n+1)/SO(n))$, and $\overrightarrow{Op}$ and $\Psi^{-1}(q)$ are regarded as 
vectors of $\mathbb R^{n+1}$ (see Figure 2).  
Hence we have 
$$\iota(\Psi_p(v))=\cosh(\|v\|)\cdot\overrightarrow{Op}+\sqrt{-1}\cdot\frac{\sinh(\|v\|)}{\|v\|}\cdot v\quad\,\,
(v\in T_pS^n(1)).$$

\vspace{0.1truecm}

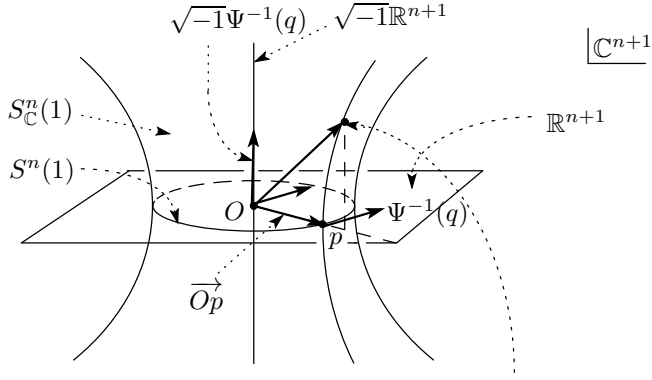
\centerline{
\unitlength 0.1in
\begin{picture}( 47.9400, 22.6100)( -6.8400,-30.2000)
%
\special{pn 8}%
\special{ar 2060 1852 528 132  6.2831853 6.2831853}%
\special{ar 2060 1852 528 132  0.0000000 3.1415927}%
%
\special{pn 8}%
\special{ar 2060 1852 528 132  3.1415927 3.3234108}%
\special{ar 2060 1852 528 132  3.4325017 3.6143199}%
\special{ar 2060 1852 528 132  3.7234108 3.9052290}%
\special{ar 2060 1852 528 132  4.0143199 4.1961381}%
\special{ar 2060 1852 528 132  4.3052290 4.4870472}%
\special{ar 2060 1852 528 132  4.5961381 4.7779563}%
\special{ar 2060 1852 528 132  4.8870472 5.0688654}%
\special{ar 2060 1852 528 132  5.1779563 5.3597745}%
\special{ar 2060 1852 528 132  5.4688654 5.6506836}%
\special{ar 2060 1852 528 132  5.7597745 5.9415927}%
\special{ar 2060 1852 528 132  6.0506836 6.2325017}%
%
\special{pn 8}%
\special{ar 120 1852 1410 1142  5.5347359 6.2831853}%
\special{ar 120 1852 1410 1142  0.0000000 0.7853982}%
%
\special{pn 8}%
\special{ar 3998 1852 1410 1142  2.3561945 3.8903513}%
%
\special{pn 8}%
\special{ar 3514 1976 1092 1336  2.5906354 3.8649700}%
%
\special{pn 20}%
\special{sh 1}%
\special{ar 2420 1948 10 10 0  6.28318530717959E+0000}%
\special{sh 1}%
\special{ar 2420 1948 10 10 0  6.28318530717959E+0000}%
\put(20.1500,-19.2200){\makebox(0,0)[rb]{$O$}}%
\put(24.5500,-19.9200){\makebox(0,0)[lt]{$p$}}%
\put(11.2500,-17.3700){\makebox(0,0)[rb]{$S^n(1)$}}%
\put(16.1800,-9.5600){\makebox(0,0)[lb]{$\sqrt{-1}\Psi^{-1}(q)$}}%
\put(27.5500,-18.1600){\makebox(0,0)[lt]{$\Psi^{-1}(q)$}}%
\put(11.1600,-14.3000){\makebox(0,0)[rb]{$S^n_{\mathbb C}(1)$}}%
\put(32.1000,-30.2000){\makebox(0,0)[lt]{$\,$}}%
%
\special{pn 13}%
\special{pa 2430 1958}%
\special{pa 2728 1870}%
\special{fp}%
\special{sh 1}%
\special{pa 2728 1870}%
\special{pa 2658 1870}%
\special{pa 2678 1884}%
\special{pa 2670 1908}%
\special{pa 2728 1870}%
\special{fp}%
%
\special{pn 20}%
\special{sh 1}%
\special{ar 2536 1412 10 10 0  6.28318530717959E+0000}%
\special{sh 1}%
\special{ar 2536 1412 10 10 0  6.28318530717959E+0000}%
\put(17.7000,-27.9000){\makebox(0,0)[lt]{{\small $\displaystyle{\iota(q)=\cosh(\|\Psi^{-1}(q)\|)\cdot\overrightarrow{Op}+\sqrt{-1}\cdot\frac{\sinh(\{\Psi^{-1}(q)\|}{\|\Psi^{-1}(q)\|}\cdot\Psi^{-1}(q)}$}}}%
%
\special{pn 20}%
\special{sh 1}%
\special{ar 2060 1852 10 10 0  6.28318530717959E+0000}%
\special{sh 1}%
\special{ar 2060 1852 10 10 0  6.28318530717959E+0000}%
%
\special{pn 13}%
\special{pa 2060 1852}%
\special{pa 2420 1948}%
\special{fp}%
\special{sh 1}%
\special{pa 2420 1948}%
\special{pa 2362 1912}%
\special{pa 2368 1934}%
\special{pa 2350 1950}%
\special{pa 2420 1948}%
\special{fp}%
%
\special{pn 13}%
\special{pa 2050 1844}%
\special{pa 2060 1466}%
\special{fp}%
\special{sh 1}%
\special{pa 2060 1466}%
\special{pa 2038 1532}%
\special{pa 2058 1518}%
\special{pa 2078 1532}%
\special{pa 2060 1466}%
\special{fp}%
%
\special{pn 13}%
\special{pa 2050 1844}%
\special{pa 2350 1756}%
\special{fp}%
\special{sh 1}%
\special{pa 2350 1756}%
\special{pa 2280 1756}%
\special{pa 2300 1770}%
\special{pa 2292 1794}%
\special{pa 2350 1756}%
\special{fp}%
%
\special{pn 8}%
\special{pa 2536 1976}%
\special{pa 2536 1404}%
\special{da 0.070}%
%
\special{pn 13}%
\special{pa 2060 1852}%
\special{pa 2518 1422}%
\special{fp}%
\special{sh 1}%
\special{pa 2518 1422}%
\special{pa 2456 1452}%
\special{pa 2478 1458}%
\special{pa 2482 1482}%
\special{pa 2518 1422}%
\special{fp}%
%
\special{pn 8}%
\special{pa 2060 1852}%
\special{pa 2800 2046}%
\special{da 0.070}%
%
\special{pn 8}%
\special{pa 1874 2204}%
\special{pa 2210 1904}%
\special{dt 0.045}%
\special{sh 1}%
\special{pa 2210 1904}%
\special{pa 2146 1934}%
\special{pa 2170 1940}%
\special{pa 2174 1964}%
\special{pa 2210 1904}%
\special{fp}%
\put(19.0000,-22.4700){\makebox(0,0)[rt]{$\overrightarrow{Op}$}}%
%
\special{pn 8}%
\special{pa 3810 1106}%
\special{pa 4110 1106}%
\special{fp}%
\put(38.3700,-10.7000){\makebox(0,0)[lb]{$\mathbb C^{n+1}$}}%
%
\special{pn 8}%
\special{pa 3810 1106}%
\special{pa 3810 922}%
\special{fp}%
%
\special{pn 8}%
\special{ar 2470 2730 954 1318  4.8733705 4.8839339}%
\special{ar 2470 2730 954 1318  4.9156240 4.9261874}%
\special{ar 2470 2730 954 1318  4.9578775 4.9684409}%
\special{ar 2470 2730 954 1318  5.0001311 5.0106944}%
\special{ar 2470 2730 954 1318  5.0423846 5.0529480}%
\special{ar 2470 2730 954 1318  5.0846381 5.0952015}%
\special{ar 2470 2730 954 1318  5.1268916 5.1374550}%
\special{ar 2470 2730 954 1318  5.1691452 5.1797085}%
\special{ar 2470 2730 954 1318  5.2113987 5.2219621}%
\special{ar 2470 2730 954 1318  5.2536522 5.2642156}%
\special{ar 2470 2730 954 1318  5.2959057 5.3064691}%
\special{ar 2470 2730 954 1318  5.3381592 5.3487226}%
\special{ar 2470 2730 954 1318  5.3804128 5.3909761}%
\special{ar 2470 2730 954 1318  5.4226663 5.4332297}%
\special{ar 2470 2730 954 1318  5.4649198 5.4754832}%
\special{ar 2470 2730 954 1318  5.5071733 5.5177367}%
\special{ar 2470 2730 954 1318  5.5494268 5.5599902}%
\special{ar 2470 2730 954 1318  5.5916804 5.6022437}%
\special{ar 2470 2730 954 1318  5.6339339 5.6444973}%
\special{ar 2470 2730 954 1318  5.6761874 5.6867508}%
\special{ar 2470 2730 954 1318  5.7184409 5.7290043}%
\special{ar 2470 2730 954 1318  5.7606944 5.7712578}%
\special{ar 2470 2730 954 1318  5.8029480 5.8135113}%
\special{ar 2470 2730 954 1318  5.8452015 5.8557649}%
\special{ar 2470 2730 954 1318  5.8874550 5.8980184}%
\special{ar 2470 2730 954 1318  5.9297085 5.9402719}%
\special{ar 2470 2730 954 1318  5.9719621 5.9825254}%
\special{ar 2470 2730 954 1318  6.0142156 6.0247790}%
\special{ar 2470 2730 954 1318  6.0564691 6.0670325}%
\special{ar 2470 2730 954 1318  6.0987226 6.1092860}%
\special{ar 2470 2730 954 1318  6.1409761 6.1515395}%
\special{ar 2470 2730 954 1318  6.1832297 6.1937930}%
\special{ar 2470 2730 954 1318  6.2254832 6.2360466}%
\special{ar 2470 2730 954 1318  6.2677367 6.2783001}%
%
\special{pn 8}%
\special{pa 2640 1440}%
\special{pa 2536 1422}%
\special{dt 0.045}%
\special{sh 1}%
\special{pa 2536 1422}%
\special{pa 2598 1452}%
\special{pa 2588 1430}%
\special{pa 2604 1414}%
\special{pa 2536 1422}%
\special{fp}%
%
\special{pn 8}%
\special{ar 2050 1238 230 370  1.7827244 1.8228582}%
\special{ar 2050 1238 230 370  1.9432595 1.9833933}%
\special{ar 2050 1238 230 370  2.1037946 2.1439284}%
\special{ar 2050 1238 230 370  2.2643297 2.3044635}%
\special{ar 2050 1238 230 370  2.4248649 2.4649986}%
\special{ar 2050 1238 230 370  2.5854000 2.6255337}%
\special{ar 2050 1238 230 370  2.7459351 2.7860689}%
\special{ar 2050 1238 230 370  2.9064702 2.9466040}%
\special{ar 2050 1238 230 370  3.0670053 3.1071391}%
%
\special{pn 8}%
\special{pa 1998 1598}%
\special{pa 2050 1614}%
\special{dt 0.045}%
\special{sh 1}%
\special{pa 2050 1614}%
\special{pa 1994 1576}%
\special{pa 2000 1598}%
\special{pa 1980 1614}%
\special{pa 2050 1614}%
\special{fp}%
%
\special{pn 8}%
\special{pa 1408 1668}%
\special{pa 844 2046}%
\special{fp}%
%
\special{pn 8}%
\special{pa 844 2046}%
\special{pa 2386 2046}%
\special{fp}%
%
\special{pn 8}%
\special{pa 2650 2046}%
\special{pa 2808 2046}%
\special{fp}%
%
\special{pn 8}%
\special{pa 2808 2046}%
\special{pa 3258 1668}%
\special{fp}%
%
\special{pn 8}%
\special{pa 3258 1668}%
\special{pa 2658 1668}%
\special{fp}%
%
\special{pn 8}%
\special{pa 2588 1668}%
\special{pa 2492 1668}%
\special{fp}%
%
\special{pn 8}%
\special{pa 2420 1668}%
\special{pa 2298 1668}%
\special{fp}%
%
\special{pn 8}%
\special{pa 2210 1668}%
\special{pa 2086 1668}%
\special{fp}%
%
\special{pn 8}%
\special{pa 2006 1668}%
\special{pa 1548 1668}%
\special{fp}%
%
\special{pn 8}%
\special{pa 1486 1668}%
\special{pa 1398 1668}%
\special{fp}%
%
\special{pn 8}%
\special{pa 2888 1668}%
\special{pa 2888 1756}%
\special{dt 0.045}%
\special{sh 1}%
\special{pa 2888 1756}%
\special{pa 2908 1688}%
\special{pa 2888 1702}%
\special{pa 2868 1688}%
\special{pa 2888 1756}%
\special{fp}%
\put(35.9200,-14.6500){\makebox(0,0)[lb]{$\mathbb R^{n+1}$}}%
%
\special{pn 8}%
\special{ar 3558 1702 670 318  3.2563675 3.2807083}%
\special{ar 3558 1702 670 318  3.3537306 3.3780714}%
\special{ar 3558 1702 670 318  3.4510937 3.4754344}%
\special{ar 3558 1702 670 318  3.5484567 3.5727975}%
\special{ar 3558 1702 670 318  3.6458198 3.6701606}%
\special{ar 3558 1702 670 318  3.7431829 3.7675237}%
\special{ar 3558 1702 670 318  3.8405460 3.8648868}%
\special{ar 3558 1702 670 318  3.9379091 3.9622499}%
\special{ar 3558 1702 670 318  4.0352722 4.0596129}%
\special{ar 3558 1702 670 318  4.1326352 4.1569760}%
\special{ar 3558 1702 670 318  4.2299983 4.2543391}%
\special{ar 3558 1702 670 318  4.3273614 4.3517022}%
\special{ar 3558 1702 670 318  4.4247245 4.4490653}%
\special{ar 3558 1702 670 318  4.5220876 4.5464284}%
\special{ar 3558 1702 670 318  4.6194507 4.6437914}%
%
\special{pn 8}%
\special{pa 1160 1378}%
\special{pa 1618 1474}%
\special{dt 0.045}%
\special{sh 1}%
\special{pa 1618 1474}%
\special{pa 1558 1442}%
\special{pa 1566 1464}%
\special{pa 1550 1480}%
\special{pa 1618 1474}%
\special{fp}%
%
\special{pn 8}%
\special{ar 1152 1904 494 246  4.7232450 4.7556774}%
\special{ar 1152 1904 494 246  4.8529747 4.8854072}%
\special{ar 1152 1904 494 246  4.9827045 5.0151369}%
\special{ar 1152 1904 494 246  5.1124342 5.1448666}%
\special{ar 1152 1904 494 246  5.2421639 5.2745964}%
\special{ar 1152 1904 494 246  5.3718937 5.4043261}%
\special{ar 1152 1904 494 246  5.5016234 5.5340558}%
\special{ar 1152 1904 494 246  5.6313531 5.6637855}%
\special{ar 1152 1904 494 246  5.7610828 5.7935153}%
\special{ar 1152 1904 494 246  5.8908126 5.9232450}%
\special{ar 1152 1904 494 246  6.0205423 6.0529747}%
%
\special{pn 8}%
\special{pa 1646 1878}%
\special{pa 1662 1940}%
\special{dt 0.045}%
\special{sh 1}%
\special{pa 1662 1940}%
\special{pa 1664 1870}%
\special{pa 1648 1888}%
\special{pa 1626 1880}%
\special{pa 1662 1940}%
\special{fp}%
%
\special{pn 8}%
\special{pa 2060 2678}%
\special{pa 2060 1000}%
\special{fp}%
%
\special{pn 8}%
\special{pa 2438 912}%
\special{pa 2068 1096}%
\special{dt 0.045}%
\special{sh 1}%
\special{pa 2068 1096}%
\special{pa 2138 1084}%
\special{pa 2116 1072}%
\special{pa 2120 1048}%
\special{pa 2068 1096}%
\special{fp}%
\put(24.7300,-9.2900){\makebox(0,0)[lb]{$\sqrt{-1}\mathbb R^{n+1}$}}%
%
\special{pn 8}%
\special{pa 2588 2046}%
\special{pa 2606 2046}%
\special{fp}%
%
\special{pn 8}%
\special{pa 2552 2046}%
\special{pa 2596 2046}%
\special{fp}%
%
\special{pn 8}%
\special{pa 1822 1246}%
\special{pa 1822 992}%
\special{dt 0.045}%
\end{picture}%
\hspace{4.5truecm}}

\vspace{0.1truecm}

\centerline{{\bf Figure 2.}}

\vspace{0.5truecm}

In 1993, M.B. Stenzel (\cite{St}) gave a construction of complete Ricci-flat metrics on the cotangent bundle 
$T^{\ast}(G/K)$ of $G/K$ in the case where the rank of $G/K$ is equal to one, where we note that $T^{\ast}(G/K)$ is 
identified with $T(G/K)(\approx G^{\Bbb C}/K^{\Bbb C})$ by the metric of $G/K$.  
In 2004, R. Bielawski (\cite{B2}) gave a construction of complete Ricci-flat metrics on $G^{\Bbb C}/K^{\Bbb C}$ 
in the case where the rank of $G/K$ is general.  
These complete Ricci metrics give Calabi-Yau structures on $G^{\Bbb C}/K^{\Bbb C}$ together with the 
natural complex structure $J_0$ and the natural non-vanishing closed holomorphic $(n,0)$-form $\Omega_0$ on 
$G^{\Bbb C}/K^{\Bbb C}$.  
H. Anciaux (\cite{An}) constructed special Lagrangian submanifolds of some phase in the complexification 
$SO(n+1,{\Bbb C})/SO(n,{\Bbb C})$ of the $n$-dimensional sphere $SO(n+1)/SO(n)$ which are invariant under 
the natural action $SO(n)\curvearrowright SO(n+1,{\Bbb C})/SO(n,{\Bbb C})$.  
M. Ionel and M. Min-Oo (\cite{IO}) constructed cohomogeneity one special Lagrangian submanifolds of some phase in 
$SO(4,{\Bbb C})/SO(3,{\Bbb C})$ which are invariant under the natural action 
$SO(2)\times SO(2)\curvearrowright SO(4,{\Bbb C})/SO(3,{\Bbb C})$.  
K. Hashimoto and T. Sakai (\cite{HS}) constructed cohomogeneity one special Lagrangian submanifolds of any phase 
in $SO(n+1,{\Bbb C})/SO(n,{\Bbb C})$ which are invariant under the natural action 
$SO(p)\times SO(n+1-p)\curvearrowright SO(n+1,{\Bbb C})/SO(n,{\Bbb C})$ ($1\leq p\leq[(n+1)/2]$).  
Later, K. Hashimoto and K. Mashimo (\cite{HM}) constructed cohomogeneity one special Lagrangian submanifolds 
of any phase in $SO(n+1,{\Bbb C})/SO(n,{\Bbb C})$ which are invariant under the natural action 
$K\curvearrowright SO(n+1,{\Bbb C})/SO(n,{\Bbb C})$ induced from the linear isotropy action 
$K\curvearrowright SO(n+1)/SO(n)(=S^n(1)\subset T_{eK}(G/K))$ of any irreducible rank two symmetric space $G/K$, 
where $n:={\rm dim}\,G/K-1$.  
Recently M. Arai and K. Baba (\cite{AB}) constructed cohomogeneity one special Lagarangian submanifolds of any 
phase and in the complexification $SL(n+1,\mathbb C)/(SL(1,\mathbb C)\times SL(n,\mathbb C))=T(\mathbb CP^n)$ 
of the complex projective space $\mathbb CP^n=SU(n+1)/S(U(1)\times U(n))$.  

In this paper, we first construct an almost Calabi-Yau structure $(J_0,\omega_{\psi_f},$\newline
$\Omega_0)$ 
on the complexification $G^{\Bbb C}/K^{\Bbb C}$, which is invariant under the natural action 
$G\curvearrowright G^{\Bbb C}/K^{\Bbb C}$, in terms of a $C^{\infty}$-function $f$ over 
$\mathbb R^l$ ($l:$ a natural number) and investigate in what case it is a Calabi-Yau structure, where 
$J_0$ and $\Omega_0$ are the natural complex structure and the natural non-vanishing closed holomorphic 
$(n,0)$-form on $G^{\mathbb C}/K^{\mathbb C}$ (Section 2).  
In Section 3, we investigate the $0$-level set of the moment map of a Hamiltonian action on the Calabi-Yau manifold 
$(G^{\mathbb C}/K^{\mathbb C},J_0,\omega_{\psi_f},\Omega_0)$.  
Let $H$ be a symmetric subgroup of $G$.  The natural action $H\curvearrowright G/K$ (which is called a 
{\it Hermann action}) is extended to the action on $G^{\mathbb C}/K^{\mathbb C}$ naturally.  
This extended action $H\curvearrowright G^{\mathbb C}/K^{\mathbb C}$ is a Hamiltonian action.  
In section 4, we investigate the orbit structure of this Hamiltonian action 
$H\curvearrowright G^{\mathbb C}/K^{\mathbb C}$.  
In Section 5, in the case where $\beta_{\psi_f}(\cdot,\cdot):=\omega_{\psi_f}(J_0(\cdot),\cdot)$ is the metric 
generalized the Stenzel metric, we first give a construction of an $H$-invariant special Lagrangian submanifold of 
cohomogeneity $r$ in $(G^{\mathbb C}/K^{\mathbb C},J_0,\omega_{\psi_f},\Omega_0)$, where $r$ denotes 
the cohomogeneity of $H\curvearrowright G/K$ (see Theorem 5.4 and Corollary 5.5).  

\section{Calabi-Yau structures on complexified symmetric spaces}
Let $G$ be a compact semi-simple Lie group and $\theta$ an involutive automorphism of $G$.  
Let $K$ be a closed subgroup of $G$ with $({\rm Fix}\,\theta)_0\subset K\subset{\rm Fix}\,\theta$, where 
${\rm Fix}\,\theta$ is the fixed point group of $\theta$ and $({\rm Fix}\,\theta)_0$ is the identity component 
of ${\rm Fix}\,\theta$.  Denote by $\mathfrak g$ (resp. $\mathfrak k$) the Lie algebra of $G$ (resp. $K$) and 
$B$ the Killing form of $\mathfrak g$.  
Denote by the same symbol $\theta$ the involution of $\mathfrak g$ induced from $\theta$.  
Set $\mathfrak p:={\rm Ker}(\theta+{\rm id}_{\mathfrak g})$, which is identified with 
the tangent space $T_o(G/K)$ of $G/K$ at $o:=eK$ ($e\,:\,$ the identity element of $G$), where 
${\rm id}_{\mathfrak g}$ is the identity transformation of $\mathfrak g$.  
Since $B\vert_{\mathfrak p\times\mathfrak p}$ is the ${\rm Ad}_G(K)\vert_{\mathfrak p}$-invariant, 
we obtain a $G$-invariant metric $\beta$ on $G/K$ with $\beta_{eK}=B$, where ${\rm Ad}_G$ is adjoint 
representation of $G$.  This Riemannian manifold $(G/K,\beta)$ is called a {\it (Riemannian) symmetric space of 
compact type}.  The dimension of maximal flat totally geodesic submanifold in $G/K$ is called the {\it rank} 
of $G/K$.  Denote by $\widetilde r$ the rank of $G/K$.  Also, assume that $G$ and $K$ admit faithful real 
representations.  Hence the complexifications $G^{\Bbb C}$ and $K^{\Bbb C}$ of $G$ and $K$ are defined.  
For the complexification $B^{\Bbb C}(:\mathfrak p^{\Bbb C}\times\mathfrak p^{\Bbb C}\to{\Bbb C}$) of $B$, 
its real part ${\rm Re}\,B^{\Bbb C}$ is 
a ${\rm Ad}_{G^{\Bbb C}}(K^{\Bbb C})\vert_{\mathfrak p^{\Bbb C}}$-invariant non-degenerate bilinear form 
(of half index) of $\mathfrak p^{\Bbb C}(=T_{eK^{\Bbb C}}(G^{\Bbb C}/K^{\Bbb C}))$ and hence we obtain 
a $G^{\Bbb C}$-invariant neutral metric $\beta_A$ on $G^{\Bbb C}/K^{\Bbb C}$ with 
$(\beta_A)_{eK}={\rm Re}\,B^{\Bbb C}$, where ${\rm Ad}_{G^{\Bbb C}}$ is adjoint representation of 
$G^{\Bbb C}$.  This pseudo-Riemannian manifold $(G^{\Bbb C}/K^{\Bbb C},\beta_A)$ is called an 
{\it anti-Kaehler symmetric space}, which is one of semi-simple pseudo-Riemannian symmetric spaces.  
Note that the terminology ``anti-Kaehler'' is used in \cite{BFV} and \cite{Koi3,Koi4} for example.  
Define $j:\mathfrak p^{\Bbb C}\to\mathfrak p^{\Bbb C}$ by $j(v):=\sqrt{-1}v$ 
($v\in\mathfrak p^{\Bbb C}$).  
Since $j$ is the ${\rm Ad}_{G^{\Bbb C}}(K^{\Bbb C})\vert_{\mathfrak p^{\Bbb C}}$-invariant, we obtain 
a $G^{\Bbb C}$-invariant almost complex structure $J_0$ of $G^{\Bbb C}/K^{\Bbb C}$ with $(J_0)_{eK^{\Bbb C}}=j$.  
Take an orthonormal base $(e_1,\cdots,e_n)$ of $\mathfrak p$ with respect to $B$ and let 
$(\theta^1,\cdots,\theta^n)$ be the dual base of $(e_1,\cdots,e_n)$.  
Set $(\theta^1)^{\Bbb C}\wedge\cdots\wedge(\theta^n)^{\Bbb C}$.  
Since $(\theta^1)^{\Bbb C}\wedge\cdots\wedge(\theta^n)^{\Bbb C}$ is 
${\rm Ad}_{G^{\Bbb C}}(K^{\Bbb C})\vert_{\mathfrak p^{\Bbb C}}$-invariant, 
we obtain a $G^{\Bbb C}$-invariant holomorphic $(n,0)$-form $\Omega_0$ on $G^{\Bbb C}/K^{\Bbb C}$ with 
$(\Omega_0)_{eK^{\Bbb C}}=(\theta^1)^{\Bbb C}\wedge\cdots\wedge(\theta^n)^{\Bbb C}$.  
Let $\psi$ be a strictly plurisubharmonic function over $G^{\Bbb C}/K^{\Bbb C}$, where we note that 
``strictly plurisubharmonicity'' means that 
the Hermitian matrix $\displaystyle{\left(\frac{\partial^2\psi}{\partial z_i\partial\bar z_j}\right)}$ 
is positive (or equivalently, $(\sqrt{-1}\partial\overline{\partial}\psi)(Z,\overline Z)>0$ holds for any 
nonzero $(1,0)$-vector $Z$).  Then $\omega_{\psi}:=\sqrt{-1}\partial\overline{\partial}\psi$ is a real 
non-degenerate closed $2$-form on $G^{\Bbb C}/K^{\Bbb C}$ and the symmetric $(0,2)$-tensor field $\beta_{\psi}$ 
associated with $J_0$ and $\omega_{\psi}$ is positive definite.  
Hence $(J_0,\omega_{\psi},\Omega_0)$ is an almost Calabi-Yau structure on $G^{\Bbb C}/K^{\Bbb C}$.  
Thus, from each strictly plurisubharmonic function over $G^{\Bbb C}/K^{\Bbb C}$, we obtain an almost Calabi-Yau 
structure on $G^{\Bbb C}/K^{\Bbb C}$.  
Hence we suffice to construct a strictly plurisubharmonic function on $G^{\Bbb C}/K^{\Bbb C}$ to construct 
an almost Calabi-Yau structure on $G^{\Bbb C}/K^{\Bbb C}$.  
Denote by ${\rm Exp}_p$ the exponential map of the anti-Kaehler manifold $(G^{\Bbb C}/K^{\Bbb C},\beta_A)$ 
at $p(\in G^{\Bbb C}/K^{\Bbb C})$ and $\exp$ the exponentional map of the Lie group $G^{\Bbb C}$.  
Set $\mathfrak g^d:=\mathfrak k\oplus\sqrt{-1}\mathfrak p(\subset\mathfrak g^{\Bbb C})$ and 
$G^d=\exp(\mathfrak g^d)$.  Denote by $\beta_{G/K}$ the $G$-invariant (Riemannian) metric on $G/K$ 
induced from $B\vert_{\mathfrak p\times\mathfrak p}$ and $\beta_{G^d/K}$ the $G^d$-invariant (Riemannian) metric 
on $G^d/K$ induced from $-({\rm Re}\,B^{\Bbb C})\vert_{\sqrt{-1}\mathfrak p\times\sqrt{-1}\mathfrak p}$.  
We may assume that the metric of $G/K$ is equal to $\beta_{G/K}$ by homothetically transforming the metric of 
$G/K$ if necessary.  On the other hand, the Riemannian manifold $(G^d/K,\beta_{G^d/K})$ is a (Riemannian) 
symmetric space of non-compact type.  
The orbit $G\cdot o$ is isometric to $(G/K,\beta_{G/K})$ and the normal umbrella 
${\rm Exp}_o(T^{\perp}_o(G\cdot o))(=G^d\cdot o)$ is isometric to $(G^d/K,\beta_{G^d/K})$.  
The complexification $\mathfrak p^{\Bbb C}$ of $\mathfrak p$ is identified with $T_o(G^{\Bbb C}/K^{\Bbb C})$ and 
$\sqrt{-1}\mathfrak p$ is identified with $T_o({\rm Exp}_o(T^{\perp}_o(G\cdot o)))$.  
Let $\mathfrak a$ be a maximal abelian subspace of $\mathfrak p$, where we note that 
${\rm dim}\,\mathfrak a=\widetilde r$.  Denote by $W$ the Weyl group of $G^d/K$ with 
respect to $\sqrt{-1}\mathfrak a$.  This group acts on $\sqrt{-1}\mathfrak a$.  
Let $C(\subset\sqrt{-1}\mathfrak a)$ be a Weyl domain (i.e., a fundamental domain of the action 
$W\curvearrowright\sqrt{-1}\mathfrak a$).  Then we have 
$G\cdot{\rm Exp}_o(\overline C)=G^{\Bbb C}/K^{\Bbb C}$, where $\overline C$ is the closure of $C$.  
For a connected open neighborhood $D$ of $0$ in $\sqrt{-1}\mathfrak a$, we define a neighborhood $U_1(D)$ of $o$ 
in ${\rm Exp}_o(\sqrt{-1}\mathfrak a)$ by $U_1(D):={\rm Exp}_o(D)$, a neighborhood $U_2(D)$ of $o$ in $G^d/K$ by 
$U_2(D):=K\cdot U_1(D)$ and a tubular neighborhood $U_3(D)$ of $G\cdot o$ in $G^{\Bbb C}/K^{\Bbb C}$ by 
$U_3(D):=G\cdot U_1(D)$ and (see Figure 3).  
Denote by ${\rm Conv}_W^+(D)$ the space of all $W$-invariant strictly convex ($C^{\infty}$-)functions over $D$, 
${\rm Conv}_K^+(U_2(D))$ the space of all $K$-invariant strictly convex ($C^{\infty}$-)functions 

\vspace{0.5truecm}

\centerline{
\unitlength 0.1in
\begin{picture}( 31.5000, 20.3000)(  6.3000,-26.0000)
%
\special{pn 8}%
\special{pa 1400 600}%
\special{pa 800 1200}%
\special{pa 800 2600}%
\special{pa 1400 2000}%
\special{pa 1400 2000}%
\special{pa 1400 600}%
\special{fp}%
%
\special{pn 8}%
\special{ar 1120 1620 130 340  0.0083331 6.2831853}%
%
\special{pn 20}%
\special{sh 1}%
\special{ar 1120 1640 10 10 0  6.28318530717959E+0000}%
\special{sh 1}%
\special{ar 1120 1640 10 10 0  6.28318530717959E+0000}%
%
\special{pn 8}%
\special{pa 850 1640}%
\special{pa 2830 1640}%
\special{fp}%
%
\special{pn 8}%
\special{pa 760 1640}%
\special{pa 630 1640}%
\special{fp}%
%
\special{pn 8}%
\special{pa 800 1940}%
\special{pa 1400 1370}%
\special{fp}%
%
\special{pn 13}%
\special{pa 1010 1760}%
\special{pa 1240 1530}%
\special{fp}%
%
\special{pn 8}%
\special{pa 840 1280}%
\special{pa 2820 1280}%
\special{fp}%
%
\special{pn 8}%
\special{pa 770 1280}%
\special{pa 640 1280}%
\special{fp}%
%
\special{pn 8}%
\special{pa 770 1960}%
\special{pa 640 1960}%
\special{fp}%
%
\special{pn 8}%
\special{pa 850 1960}%
\special{pa 2830 1960}%
\special{fp}%
%
\special{pn 8}%
\special{ar 2050 1620 130 340  0.0083331 6.2831853}%
%
\special{pn 4}%
\special{pa 700 1960}%
\special{pa 650 1910}%
\special{dt 0.027}%
\special{pa 760 1960}%
\special{pa 650 1850}%
\special{dt 0.027}%
\special{pa 820 1960}%
\special{pa 650 1790}%
\special{dt 0.027}%
\special{pa 880 1960}%
\special{pa 650 1730}%
\special{dt 0.027}%
\special{pa 940 1960}%
\special{pa 650 1670}%
\special{dt 0.027}%
\special{pa 1000 1960}%
\special{pa 650 1610}%
\special{dt 0.027}%
\special{pa 1060 1960}%
\special{pa 650 1550}%
\special{dt 0.027}%
\special{pa 1120 1960}%
\special{pa 650 1490}%
\special{dt 0.027}%
\special{pa 1180 1960}%
\special{pa 650 1430}%
\special{dt 0.027}%
\special{pa 1240 1960}%
\special{pa 650 1370}%
\special{dt 0.027}%
\special{pa 1300 1960}%
\special{pa 650 1310}%
\special{dt 0.027}%
\special{pa 1360 1960}%
\special{pa 680 1280}%
\special{dt 0.027}%
\special{pa 1420 1960}%
\special{pa 740 1280}%
\special{dt 0.027}%
\special{pa 1480 1960}%
\special{pa 800 1280}%
\special{dt 0.027}%
\special{pa 1540 1960}%
\special{pa 860 1280}%
\special{dt 0.027}%
\special{pa 1600 1960}%
\special{pa 920 1280}%
\special{dt 0.027}%
\special{pa 1660 1960}%
\special{pa 980 1280}%
\special{dt 0.027}%
\special{pa 1720 1960}%
\special{pa 1040 1280}%
\special{dt 0.027}%
\special{pa 1780 1960}%
\special{pa 1100 1280}%
\special{dt 0.027}%
\special{pa 1840 1960}%
\special{pa 1160 1280}%
\special{dt 0.027}%
\special{pa 1900 1960}%
\special{pa 1220 1280}%
\special{dt 0.027}%
\special{pa 1960 1960}%
\special{pa 1280 1280}%
\special{dt 0.027}%
\special{pa 2020 1960}%
\special{pa 1340 1280}%
\special{dt 0.027}%
\special{pa 2080 1960}%
\special{pa 1400 1280}%
\special{dt 0.027}%
\special{pa 2140 1960}%
\special{pa 1460 1280}%
\special{dt 0.027}%
\special{pa 2200 1960}%
\special{pa 1520 1280}%
\special{dt 0.027}%
\special{pa 2260 1960}%
\special{pa 1580 1280}%
\special{dt 0.027}%
\special{pa 2320 1960}%
\special{pa 1640 1280}%
\special{dt 0.027}%
\special{pa 2380 1960}%
\special{pa 1700 1280}%
\special{dt 0.027}%
\special{pa 2440 1960}%
\special{pa 1760 1280}%
\special{dt 0.027}%
%
\special{pn 4}%
\special{pa 2500 1960}%
\special{pa 1820 1280}%
\special{dt 0.027}%
\special{pa 2560 1960}%
\special{pa 1880 1280}%
\special{dt 0.027}%
\special{pa 2620 1960}%
\special{pa 1940 1280}%
\special{dt 0.027}%
\special{pa 2680 1960}%
\special{pa 2000 1280}%
\special{dt 0.027}%
\special{pa 2740 1960}%
\special{pa 2060 1280}%
\special{dt 0.027}%
\special{pa 2800 1960}%
\special{pa 2120 1280}%
\special{dt 0.027}%
\special{pa 2810 1910}%
\special{pa 2180 1280}%
\special{dt 0.027}%
\special{pa 2810 1850}%
\special{pa 2240 1280}%
\special{dt 0.027}%
\special{pa 2810 1790}%
\special{pa 2300 1280}%
\special{dt 0.027}%
\special{pa 2810 1730}%
\special{pa 2360 1280}%
\special{dt 0.027}%
\special{pa 2810 1670}%
\special{pa 2420 1280}%
\special{dt 0.027}%
\special{pa 2810 1610}%
\special{pa 2480 1280}%
\special{dt 0.027}%
\special{pa 2810 1550}%
\special{pa 2540 1280}%
\special{dt 0.027}%
\special{pa 2810 1490}%
\special{pa 2600 1280}%
\special{dt 0.027}%
\special{pa 2810 1430}%
\special{pa 2660 1280}%
\special{dt 0.027}%
\special{pa 2810 1370}%
\special{pa 2720 1280}%
\special{dt 0.027}%
\special{pa 2810 1310}%
\special{pa 2780 1280}%
\special{dt 0.027}%
%
\special{pn 8}%
\special{pa 1710 760}%
\special{pa 1280 940}%
\special{fp}%
\special{sh 1}%
\special{pa 1280 940}%
\special{pa 1350 934}%
\special{pa 1330 920}%
\special{pa 1334 896}%
\special{pa 1280 940}%
\special{fp}%
%
\special{pn 13}%
\special{pa 1500 2160}%
\special{pa 1070 1720}%
\special{fp}%
\special{sh 1}%
\special{pa 1070 1720}%
\special{pa 1102 1782}%
\special{pa 1108 1758}%
\special{pa 1132 1754}%
\special{pa 1070 1720}%
\special{fp}%
%
\special{pn 8}%
\special{pa 1230 2410}%
\special{pa 900 1850}%
\special{fp}%
\special{sh 1}%
\special{pa 900 1850}%
\special{pa 918 1918}%
\special{pa 928 1896}%
\special{pa 952 1898}%
\special{pa 900 1850}%
\special{fp}%
%
\special{pn 8}%
\special{pa 2660 2170}%
\special{pa 2380 1640}%
\special{fp}%
\special{sh 1}%
\special{pa 2380 1640}%
\special{pa 2394 1708}%
\special{pa 2406 1688}%
\special{pa 2430 1690}%
\special{pa 2380 1640}%
\special{fp}%
%
\special{pn 8}%
\special{pa 2610 1010}%
\special{pa 2390 1440}%
\special{fp}%
\special{sh 1}%
\special{pa 2390 1440}%
\special{pa 2438 1390}%
\special{pa 2414 1394}%
\special{pa 2404 1372}%
\special{pa 2390 1440}%
\special{fp}%
%
\special{pn 4}%
\special{pa 1114 1958}%
\special{pa 1000 1730}%
\special{fp}%
\special{pa 1150 1938}%
\special{pa 990 1620}%
\special{fp}%
\special{pa 1178 1904}%
\special{pa 1000 1550}%
\special{fp}%
\special{pa 1202 1866}%
\special{pa 1010 1478}%
\special{fp}%
\special{pa 1226 1820}%
\special{pa 1028 1426}%
\special{fp}%
\special{pa 1240 1760}%
\special{pa 1044 1368}%
\special{fp}%
\special{pa 1244 1680}%
\special{pa 1064 1318}%
\special{fp}%
\special{pa 1240 1580}%
\special{pa 1100 1298}%
\special{fp}%
\special{pa 1226 1462}%
\special{pa 1138 1284}%
\special{fp}%
%
\special{pn 13}%
\special{pa 1580 1090}%
\special{pa 1100 1460}%
\special{fp}%
\special{sh 1}%
\special{pa 1100 1460}%
\special{pa 1166 1436}%
\special{pa 1142 1428}%
\special{pa 1142 1404}%
\special{pa 1100 1460}%
\special{fp}%
\put(17.6000,-8.0000){\makebox(0,0)[lb]{$G^d\cdot o$}}%
\put(16.3000,-11.3000){\makebox(0,0)[lb]{$U_2(D)$}}%
\put(11.0000,-24.6000){\makebox(0,0)[lt]{${\rm Exp}_o(\sqrt{-1}\mathfrak a)$}}%
\put(15.3000,-21.3000){\makebox(0,0)[lt]{$U_1(D)$}}%
\put(25.5000,-9.9000){\makebox(0,0)[lb]{$U_3(D)$}}%
\put(25.9000,-21.9000){\makebox(0,0)[lt]{$G\cdot o$}}%
%
\special{pn 20}%
\special{sh 1}%
\special{ar 2050 1640 10 10 0  6.28318530717959E+0000}%
\special{sh 1}%
\special{ar 2050 1640 10 10 0  6.28318530717959E+0000}%
%
\special{pn 8}%
\special{pa 3340 800}%
\special{pa 3780 800}%
\special{fp}%
%
\special{pn 8}%
\special{pa 3350 800}%
\special{pa 3350 570}%
\special{fp}%
\put(34.0000,-7.7000){\makebox(0,0)[lb]{$G^{\mathbb C}/K^{\mathbb C}$}}%
\end{picture}%
\hspace{0truecm}}

\vspace{0.5truecm}

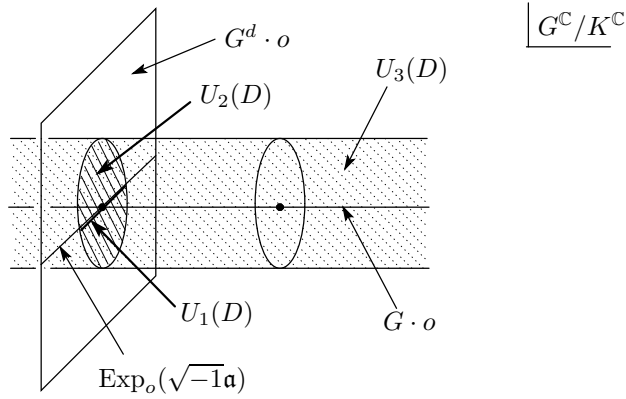
\centerline{{\bf Figure 3.}}

\vspace{0.5truecm}

\noindent
over $U_2(D)$ and $PH_G^+(U_3(D))$ the space of all $G$-invariant strictly plurisubharmonic ($C^{\infty}$-)functions 
over $U_3(D)$.  The restriction map from $U_3(D)$ to $U_2(D)$ gives an isomorphism of $PH_G^+(U_3(D))$ 
onto ${\rm Conv}_K^+(U_2(D))$ and the composition of the restriction map from $U_3(D)$ to $U_1(D)$ with 
${\rm Exp}_o$ gives an isomorphism of $PH_G^+(U_3(D))$ onto ${\rm Conv}_W^+(D)$ (see \cite{AL}).  
Hence we suffice to construct $W$-invariant strictly convex functions over $D$ or $K$-invariant strictly convex 
functions over $U_2(D)$ to construct strictly plurisubharmonic functions over $U_3(D)$.  Let $\psi$ be a 
$G$-invariant strictly plurisubharmonic ($C^{\infty}$-)functions over $U_3(D)$.  
Denote by $\bar{\psi}$ the restriction of $\psi$ to $U_2(D)$ and $\bar{\bar{\psi}}$ the composition of 
the restriction of $\psi$ to $U_1(D)$ with ${\rm Exp}_o$.  Denote by $Ric_{\psi}$ the Ricci form of $\beta_{\psi}$.  
By a result of R. Bielawski (Theorem 3.3 in \cite{B2}), we have 
$$Ric_{\psi}=-\sqrt{-1}\partial\overline{\partial}\log\,{\rm det}
\left(\frac{\partial^2\psi}{\partial z_i\partial\bar z_j}\right)
=-\sqrt{-1}\partial\overline{\partial}\log\left(\left(\frac{{\rm det}\,\nabla d\bar{\psi}}
{{\rm det}\,\beta_{G^d/K}}\right)^h\right),\leqno{(2.1)}$$
where $\nabla$ denotes the Riemannian connection of $\beta_{G^d/K}$, 
$(z_1,\cdots,z_n)$ is any complex coordinate of $G^{\Bbb C}/K^{\Bbb C}$ and 
$\displaystyle{\left(\frac{{\rm det}\,\nabla d\bar{\psi}}{{\rm det}\,\beta_{G^d/K}}\right)^h}$ is 
the $G$-invariant function over $G^{\Bbb C}/K^{\Bbb C}$ satisfying 
$$\left.\left(\frac{{\rm det}\,\nabla d\bar{\psi}}{{\rm det}\,\beta_{G^d/K}}\right)^h\right\vert_{G^d/K}
=\frac{{\rm det}\,\nabla d\bar{\psi}}{{\rm det}\,\beta_{G^d/K}}.$$
According to the result of \cite{B1}, for any given $K$-invariant positive $C^{\infty}$-function $\varphi$ on 
$G^d/K$, the Monge-Amp$\grave{\rm e}$re equation 
$$\frac{{\rm det}\,\nabla d\bar{\psi}}{{\rm det}\,\beta_{G^d/K}}=\varphi\leqno{(2.2)}$$
has a global $K$-invariant strictly convex $C^{\infty}$-solution.  

Furthermore, we can derive the following fact directly.  

\vspace{0.5truecm}

\noindent
{\bf Lemma 2.1.} {\sl {\rm (i)} For any $G$-invariant strictly plurisubharmonic ($C^{\infty}$-)function $\psi$ 
over $U_3(D)$, we have 
$$Ric_{\psi}=-\sqrt{-1}\partial\overline{\partial}\log\left(\left(\frac{{\rm det}\,\nabla^0d\bar{\bar{\psi}}}
{{\rm det}\,\beta_0}\right)^h\right),\leqno{(2.3)}$$
where $\beta_0$ is the Euclidean metric of $\sqrt{-1}\mathfrak a$ associated to 
$-{\rm Re}\,B^{\Bbb C}\vert_{\sqrt{-1}\mathfrak a\times\sqrt{-1}\mathfrak a}$  and $\nabla^0$ is 
the Euclidean connection of $\beta_0$ and 
$\displaystyle{\left(\frac{{\rm det}\,\nabla^0d\bar{\bar{\psi}}}{{\rm det}\,\beta_0}\right)^h}$ is 
the $G$-invariant function over $G^{\Bbb C}/K^{\Bbb C}$ satisfying 
$$\left.\left(\frac{{\rm det}\,\nabla^0d\bar{\bar{\psi}}}{{\rm det}\,\beta_0}\right)^h
\right\vert_{{\rm Exp}_o(\sqrt{-1}\mathfrak a)}\circ{\rm Exp}_o
=\frac{{\rm det}\,\nabla^0d\bar{\bar{\psi}}}{{\rm det}\,\beta_0}.$$

{\rm (ii)} For any given $W$-invariant positive $C^{\infty}$-function $\varphi$ on $\sqrt{-1}\mathfrak a$, 
the Monge\newline
-Amp$\grave{\rm e}$re equation 
$$\frac{{\rm det}\,\nabla^0d\rho}{{\rm det}\,\beta_0}=\varphi\leqno{(2.4)}$$
has a global $W$-invariant strictly convex $C^{\infty}$-solution.  
}

\vspace{0.5truecm}

\noindent
{\it Proof.} Since $\bar{\psi}$ is $K$-invariant, we have 
$$\left(\frac{{\rm det}\,\nabla d\bar{\psi}}{{\rm det}\,\beta_{G^d/K}}\right)^h
=\left(\frac{{\rm det}\,\nabla^0d\bar{\bar{\psi}}}{{\rm det}\,\beta_0}\right)^h.$$
Therefore the statement (i) is directly derived from the above result by R. Bielawski.  
The statement (ii) is trivial.  \qed

\vspace{0.5truecm}

From a global $W$-invariant strictly convex $C^{\infty}$-solution $\rho$ of the Monge-Amp$\grave{\rm e}$re 
equation 
$$\frac{{\rm det}\,\nabla^0d\rho}{{\rm det}\,\beta_0}=c\qquad(c:{\rm a}\,\,{\rm positive}\,\,{\rm constant})
\leqno{(2.5)}$$
we can construct a complete Ricci-flat metric $\beta_{\psi}$ on $G^{\Bbb C}/K^{\Bbb C}$, where 
$\psi$ is the $G$-invariant strictly plurisubharmonic $C^{\infty}$-function satisfying 
$\psi\vert_{{\rm Exp}_o(\sqrt{-1}\mathfrak a)}\circ{\rm Exp}_o=\rho$.  
Hence we obtain a Calabi-Yau structure $(J_0,\omega_{\psi},\Omega_0)$ on $G^{\Bbb C}/K^{\Bbb C}$ by replacing $\rho$ 
to a suitable positive constant-multiple of $\rho$ if necessary.  

We consider the case of $D=\sqrt{-1}\mathfrak a$.  
Then, according to the Schwarz's theorem ([Sc]), the ring $C^{\infty}_W(\sqrt{-1}\mathfrak a)$ of all $W$-invariant 
$C^{\infty}$-functions over $\sqrt{-1}\mathfrak a$ is given by 
$$C^{\infty}_W(\sqrt{-1}\mathfrak a)=\{f\circ(\rho_1,\cdots,\rho_l)\,|\,f\in C^{\infty}(\mathbb R^l)\},$$
where $\rho_1,\cdots,\rho_l$ are generators of $C^{\infty}_W(\sqrt{-1}\mathfrak a)$ of the ring 
${\rm Pol}_W(\sqrt{-1}\mathfrak a)$ of all $W$-invariant polynomials over $\sqrt{-1}\mathfrak a$.  
In the sequel, set $\overrightarrow{\rho}:=(\rho_1,\cdots,\rho_l)$ for simplicity. 
Let $\psi_i$ ($i=1,\cdots,l$) be the elements of $PH_G^+(G^{\mathbb C}/K^{\mathbb C})$ with 
$\bar{\bar{\psi}}_i=\rho_i$.  
In the sequel, set $\overrightarrow{\psi}:=(\psi_1,\cdots,\psi_l)$ for simplicity. 
Hence any element $\psi$ of $PH_G^+(G^{\mathbb C}/K^{\mathbb C})$ is described as 
$\psi=f\circ\overrightarrow{\psi}$ in terms of some $f\in C^{\infty}(\mathbb R^l)$.  
As the first generator $\rho_1$ of $C^{\infty}_W(\sqrt{-1}\mathfrak a)$, we take 
$$\rho_1(\sqrt{-1}v):=\|v\|^2+1\quad\,\,(v\in\mathfrak a).$$
In the following, set $\psi_f:=f\circ\overrightarrow{\psi}$.  
By using Lemma 2.1, we can derive the following fact.  

\vspace{0.4truecm}

\noindent
{\bf Theorem 2.2.} {\sl {\rm (i)} The triple $(J_0,\omega_{\psi_f},\Omega_0)$ is a Calabi-Yau structure 
of $G^{\mathbb C}/K^{\mathbb C}$ when 
{\small 
$${\rm det}\left(\sum_{k=1}^l\sum_{\widehat k=1}^l\left(
\left(\frac{\partial^2f}{\partial y_{\widehat k}\partial y_k}\circ\overrightarrow{\rho}\right)
\cdot\frac{\partial\rho_{\widehat k}}{\partial x_i}\cdot\frac{\partial\rho_k}{\partial x_j}
+\left(\frac{\partial f}{\partial y_k}\circ\overrightarrow{\rho}\right)\cdot
\frac{\partial^2\rho_k}{\partial x_i\partial x_j}\right)\right)=c,
\leqno{(2.6)}$$
}
where $c$ is a positive constant, and $(x_1,\cdots,x_r)$ and $(y_1,\cdots,y_l)$ are the natural coordinates of 
$\sqrt{-1}\mathfrak a$ and $\mathbb R^l$, respectively.  

{\rm (ii)} Assume that $\displaystyle{\frac{\partial f}{\partial y_2}=\cdots=\frac{\partial f}{\partial y_l}=0}$.  
Then $(J_0,\omega_{\psi_f},\Omega_0)$ is a Calabi-Yau structure of $G^{\mathbb C}/K^{\mathbb C}$ when 
$${\rm det}\left(2x_ix_j\cdot\left(\frac{\partial^2f}{\partial y_1^2}\circ\overrightarrow{\rho}\right)
+\left(\frac{\partial f}{\partial y_1}\circ\overrightarrow{\rho}\right)\cdot\delta_{ij}\right)=c,\leqno{(2.7)}$$
where $c$ is a positive constant, and $(x_1,\cdots,x_r)$ and $(y_1,\cdots,y_l)$ are as above.}

\vspace{0.5truecm}

\noindent
{\it Proof.} By a simple calculation, we have 
\begin{align*}
(\nabla^0d\bar{\bar{\psi}}_f)\left(\frac{\partial}{\partial x_i},\,\frac{\partial}{\partial x_j}\right)
=\sum_{k=1}^l\sum_{\widehat k=1}^l&\left(
\left(\frac{\partial^2f}{\partial y_{\widehat k}\partial y_k}\circ\overrightarrow{\rho}\right)
\cdot\frac{\partial\rho_{\widehat k}}{\partial x_i}\cdot\frac{\partial\rho_k}{\partial x_j}\right.\\
&\,\,\,\left.+\left(\frac{\partial f}{\partial y_k}\circ\overrightarrow{\rho}\right)\cdot
\frac{\partial^2\rho_k}{\partial x_i\partial x_j}\right).
\end{align*}
Hence, from $(2.6)$, we obtain 
$${\rm det}\left((\nabla^0d\bar{\bar{\psi}}_f)\left(\frac{\partial}{\partial x_i},\,\frac{\partial}{\partial x_j}
\right)\right)=c>0,$$
that is, $\bar{\bar{\psi}}_f$ is convex.  Also, we have 
$${\rm det}\left(\beta_0\left(\frac{\partial}{\partial x_i},\,\frac{\partial}{\partial x_j}\right)\right)=1.$$
Hence we have 
$$\frac{{\rm det}\,\nabla^0d\bar{\bar{\psi}}_f}{{\rm det}\,\beta_0}=c.$$
Therefore, from Lemma 2.1, we obtain ${\rm Ric}_{\psi_f}=0$.  Thus $(J_0,\omega_{\psi_f},\Omega_0)$ is a Calabi-Yau 
structure of $G^{\mathbb C}/K^{\mathbb C}$.  The statement (ii) follows from (i) direcctly.  \qed

\vspace{0.5truecm}

\noindent
{\it Remark 2.1.} (i) By using the result of \cite{B1}, we can show that the Monge-Amp$\grave{\rm e}$re type 
equation $(2.6)$ has global solution $f$.  

(ii) The Monge-Amp$\grave{\rm e}$re type equations $(2.6)$ and $(2.7)$ coincide in the case of ${\rm rank}\,G/K=1$.  

\vspace{0.5truecm}

From (ii) of Theorem 2.2, we can derive the following fact.  

\vspace{0.5truecm}

\noindent
{\bf Corollary 2.3.} {\sl Let $f$ be the $C^{\infty}$-function over $\mathbb R^l$ defined by 
$$f(y_1,\cdots,y_l):=\int_1^{y_1}(a\log\,s+b)^{\frac{1}{r}}\,ds+c,\leqno{(2.8)}$$
where $a,b$ and $c$ are positive constants.  Then $(J_0,\omega_{\psi_f},\Omega_0)$ is a Calabi-Yau structure of 
$G^{\mathbb C}/K^{\mathbb C}$.}

\vspace{0.5truecm}

\noindent
{\it Proof.} By a simple calculation, we have 
\begin{align*}
&{\rm det}\left(2x_ix_j\cdot\left(\frac{\partial^2f}{\partial y_1^2}\circ\overrightarrow{\rho}\right)
+\left(\frac{\partial f}{\partial y_1}\circ\overrightarrow{\rho}\right)\cdot\delta_{ij}\right)\\
=&2\left(\frac{\partial^2f}{\partial y_1^2}\circ\overrightarrow{\rho}\right)\cdot
\left(\frac{\partial f}{\partial y_1}\circ\overrightarrow{\rho}\right)^{r-1}\cdot\rho_1=\frac{2a}{r}.
\end{align*}
Hence, it follows from (ii) of Theorem 2.2 that $(J_0,\omega_{\psi_f},\Omega_0)$ is a Calabi-Yau structure of 
$G^{\mathbb C}/K^{\mathbb C}$.  \qed

\vspace{0.5truecm}

\noindent
{\it Remark 2.2.} For $f$ as in $(2.8)$, $\beta_{\psi_f}$ coincides with the Stenzel metric in the case where 
$G/K=SO(n+1)/SO(n)(=S^n)$.  

\section{Hamiltonian actions and the moment maps}
Let $(M,\omega)$ be a symplectic manifold and the action $H\curvearrowright M$ of a Lie group $H$ on $(M,\omega)$.  
This action $H\curvearrowright M$ is called a {\it Hamiltonian action} if it satisfies the following conditions 
(i)$\sim$(iii):

\vspace{0.25truecm}

(i) For any $h\in H$, $h^{\ast}\omega=\omega$ holds;

(ii) For any element $X$ of the Lie algebra $\mathfrak h$ of $H$, $i_{X^{\ast}}\omega$ is exact, 
where $X^{\ast}$ denote the fundamental vector field on $M$ associated to $X$, that is, 
$$X^{\ast}_p:=\left.\frac{d}{dt}\right|_{t=0}(\exp\,tX)\cdot p\,\,\,(p\in M)$$
and $i_{X^{\ast}}$ denotes the inner product operator by $X^{\ast}$; 

(iii) There exists a family $\{F_X\}_{X\in\mathfrak h}$ of $C^{\infty}$-functions over $M$ such that 
$dF_X=i_{X^{\ast}}\omega\,\,\,(X\in\mathfrak h)$ and that the correspndence 
$X\mapsto F_X\,\,\,(X\in\mathfrak h)$ is a Lie algebra homomorphism of $\mathfrak h$ into $C^{\infty}(M)$.  

\vspace{0.25truecm}

\noindent
Here we note that, by the condition (ii), it is assured that there exists a family $\{F_X\}_{X\in\mathfrak h}$ of 
$C^{\infty}$-functions over $M$ such that $dF_X=i_{X^{\ast}}\omega\,\,\,(X\in\mathfrak h)$ and that 
the correspndence $X\mapsto F_X\,\,\,(X\in\mathfrak h)$ is linear.  
For a function $F$ over $(M,\omega)$, the s-gradient vector field ${\rm sgrad}\,F$ is defined by 
$dF(Y)=\omega({\rm sgrad}\,F,Y)\,\,\,(Y\in TM)$.  
Clearly we have ${\rm sgrad}\,F_X=X^{\ast}$.  
The moment map $\mu:M\to\mathfrak h^{\ast}$ of this Hamiltonian action is defined by 
$$(\mu(p))(X):=F_X(p)\quad\,\,(p\in M,\,\,X\in\mathfrak h).$$
Hence the level set $\mu^{-1}(0)$ is given by 
$$\mu^{-1}(0)=\mathop{\cap}_{X\in\mathfrak h}F_X^{-1}(0).\leqno{(3.1)}$$

Let $(G^{\mathbb C}/K^{\mathbb C},J_0,\omega_{\psi_f},\Omega_0)$ be a Calabi-Yau manifold stated 
in the previous section and $H$ be a closed subgroup of $G$.  
Denote by $\mathfrak h$ the Lie algebra of $H$.  Let $n:={\rm dim}\,G/K$.  
For simplicity, set $M:=G\cdot o(\approx G/K),\,\,M^{\mathbb C}:=G^{\mathbb C}/K^{\mathbb C}$ and $M^d:=G^d\cdot o
(\approx G^d/K)$.  As stated in Introduction, set $\Psi_p={\rm Exp}_p\circ(J_0)_p$ ($p\in M$).  
Set $M^d_p:=\Psi_p(T_p(G\cdot o))$ ($p\in M$), which is the normal umbrella of $M$ in $(M^{\mathbb C},\beta_A)$.  
Note that $M^d_o=M^d$.  

\vspace{0.5truecm}

\noindent
{\bf Lemma 3.1.} {\sl {\rm (i)} The action $H\curvearrowright(M^{\mathbb C},J_0,\omega_{\psi_f},\Omega_0)$ is 
a Hamiltonian action and its moment map $\mu_{\psi_f}$ is given by 
$$((\mu_{\psi_f})(q))(X)=-({\rm Im}\,\overline{\partial}\psi_f)_q(X^{\ast}_q)\quad\,\,(q\in M^{\mathbb C},\,\,
X\in\mathfrak h),
\leqno{(3.2)}$$
where ${\rm Im}(\cdot)$ denotes the imaginary part of $(\cdot)$.  

{\rm (ii)}
The level set $\mu_{\psi_f}^{-1}(0)$ is given by 
$$\mu_{\psi_f}^{-1}(0)=\{q\in M^{\mathbb C}\,|\,({\rm Im}\,\overline{\partial}\psi)_q(X^{\ast}_q)=0\,\,\,
(\forall\,X\in\mathfrak h)\}.\leqno{(3.3)}$$
}

\vspace{0.5truecm}

\noindent
{\it Proof.} Since $\omega_{\psi_f}$ is $G$-invariant and $H$ is a closed subgroup, it is $H$-invariant.  
Set $\alpha:=-{\rm Im}\,\overline{\partial}\psi_f$.  For each $X\in\mathfrak h$, define a function 
$F_X$ over $M^{\mathbb C}$ by $F_X(q):=\alpha_q(X^{\ast}_q)$ ($q\in M^{\mathbb C}$).  
Then, for any tangent vector field $Y$ over $M^{\mathbb C}$, we have 
$$d\alpha(X^{\ast},Y)=X^{\ast}(\alpha(Y))-Y(\alpha(X^{\ast}))-\alpha({\mathcal L}_{X^{\ast}}Y)
=({\mathcal L}_{X^{\ast}}\alpha)(Y)-dF_X(Y),$$
where ${\mathcal L}_{X^{\ast}}$ denotes the Lie derivative with respect to $X^{\ast}$.  
Since $\alpha$ is $H$-invariant, we have ${\mathcal L}_{X^{\ast}}\alpha=0$.  
Also, we have $d\alpha=-\omega$.  Hence we obtain $dF_X=i_{X^{\ast}}\omega$.  
Also, it is clear that the correspndence $X\mapsto F_X\,\,\,(X\in\mathfrak h)$ is a Lie algebra homomorphism of 
$\mathfrak h$ into $C^{\infty}(M)$.  
Therefore the action $H\curvearrowright(M^{\mathbb C},J_0,\omega_{\psi_f},\Omega_0)$ is a Hamiltonian action and 
its moment map $\mu_{\psi_f}$ is given by 
$$(\mu_{\psi_f}(q))(X)=F_X(q)=-({\rm Im}\,\overline{\partial}\psi_f)_q(X^{\ast}_q)\quad
(q\in M^{\mathbb C},\,\,X\in\mathfrak h).$$
Thus the statement (i) has been shown.  The statement (ii) follows from (i) directly.  \qed

\vspace{0.5truecm}

By using this lemma, we obtain the following fact.  

\vspace{0.5truecm}

\noindent
{\bf Lemma 3.2.} {\sl Let $f$ be as in $(2.8)$.  Then the level set $\mu_{\psi_f}^{-1}(0)$ is given by 
$$\mu_{\psi_f}^{-1}(0)=\mathop{\amalg}_{p\in M}\Psi_p(T_p^{\perp}(H\cdot p)),\leqno{(3.4)}$$
where $T_p^{\perp}(H\cdot p)$ denotes the normal space of $H\cdot p$ in $M$ at $p$.  
Also, if ${\rm cohom}\,(H\curvearrowright G/K)=r$, then we have ${\rm dim}\,\mu_{\psi_f}^{-1}(0)=n+r$.}

\vspace{0.5truecm}

\noindent
{\it Proof.} Let $(U,(z_1=x_i+\sqrt{-1}y_i)_{i=1}^n)$ be a holomorphic coordinate of $M^{\mathbb C}$ 
such that ${\rm Span}\{(\frac{\partial}{\partial x_i})_p\,|\,i=1,\cdots,n\}=T_pM$ 
holds for any $p\in U\cap M$.  
Note that, for $q\in U\cap M_p^d$, the following relation holds:
$$q=\Psi_p\left(\sum_{i=1}^ny_i(q)\left(\frac{\partial}{\partial x_i}\right)_p\right).$$
Fix $p\in M$ and $q\in U\cap M_p^d$.  Take any $X\in\mathfrak h$.  
Then, by a simple calculation, we have 
$$\begin{array}{l}
\displaystyle{((\mu_{\psi_f})(q))(X)=-2\left(\log(\sum_{i=1}^ny_i(q)y_j(q)g_{ij}(p)+1)+b\right)^{\frac{1}{r}}}\\
\hspace{2.8truecm}\displaystyle{\times
\sum_{i=1}^n\sum_{j=1}^nX_i^{\ast}(q)y_j(q)g_{ij}(p),}
\end{array}\leqno{(3.5)}$$
where $g_{ij}:=g\left(\frac{\partial}{\partial x_i},\,\frac{\partial}{\partial x_j}\right)$, 
$X^{\ast}$ denotes the fundamental vector field on $M^{\mathbb C}$ associated to $X$ and 
$X^{\ast}_i$ is the function given by $X^{\ast}=\sum\limits_{i=1}^n\left(X^{\ast}_i\frac{\partial}{\partial x_i}
+{\widehat X}^{\ast}_i\frac{\partial}{\partial y_i}\right)$.  Hence $q\in\mu_{\psi_f}^{-1}(0)$ if and only if 
\begin{align*}
&g_p\left(\sum_{i=1}^nX_i^{\ast}(q)\left(\frac{\partial}{\partial x_i}\right)_p,\,
-(J_0)_p\left(\left({\rm Exp}_p|_{(J_0)_p(T_pM)}\right)^{-1}(q)\right)\right)\\
=&\sum_{i=1}^n\sum_{j=1}^nX^{\ast}_i(p)y_j(p)g_{ij}(p)=0
\end{align*}
holds for any $X\in\mathfrak h$.  
On the othe hand, $X$ moves over $\mathfrak h$, 
$\sum\limits_{i=1}^nX^{\ast}_i(p)\left(\frac{\partial}{\partial x_i}\right)_p$ moves over the whole of 
$T_p(H\cdot p)$.  Therefore $q\in\mu_{\psi_f}^{-1}(0)$ if and only if 
$$(J_0)_p\left(\left({\rm Exp}_p|_{(J_0)_p(T_pM)}\right)^{-1}(q)\right)\in T_p^{\perp}(H\cdot p)$$
holds.  From this fact, the relation $(3.4)$ follows.  

Let $U$ be the open subset of $G/K$ of all regular points of $H\curvearrowright G/K$.  
Then $\displaystyle{\mathop{\amalg}_{p\in U}\Psi_p\left(T_p^{\perp}(H\cdot p)\right)}$ 
is an open subset of $\mu_{\psi_f}^{-1}(0)$.  
It is clear that the dimension of this open subset is equal to $n+r$.  Hence we obtain 
${\rm dim}\,\mu_{\psi_f}^{-1}(0)=n+r$.  \qed

\section{The actions of symmetric subgroups on complexified symmetric spaces}
Let $(G^{\mathbb C}/K^{\mathbb C},J_0,\omega_{\psi_f},\Omega_0)$ be a Calabi-Yau manifold stated in Section 2.  
As in the previous section, set $M:=G\cdot o(=G/K),\,\,M^{\Bbb C}:=G^{\Bbb C}/K^{\Bbb C}$, 
$M^d:=G^d\cdot o(=G^d/K)$ and $M^d_p:=\Psi_p(T_p(G\cdot o))$.  
Let $H$ be a symmetric subgroup of $G$ and $\sigma$ the involutive automorphism of 
$G$ satisfying $({\rm Fix}\,\sigma)_0\subset H\subset{\rm Fix}\,\sigma$.  
The natural action $H$ of on $G/K(=M)$ is called a {\it Hermann action}.  
Assume that $\sigma\circ\theta=\theta\circ\sigma$.  Then the action is called a {\it commutative Hermann action}.  
Set $n:={\rm dim}\,M$ and denote by $r$ the cohomogeneity of the action $H\curvearrowright M$.  
The group $H$ acts on $M^{\Bbb C}$ as a subaction of the natural action 
$G\curvearrowright M^{\Bbb C}$, where we note that $G\curvearrowright M^{\Bbb C}$ is a Hermann type action 
(this terminology was used in \cite{Koi1}).  It is easy to show that the subaction $H\curvearrowright M^{\mathbb C}$ 
is a Hamiltonian action.  Set $\mathfrak q:={\rm Ker}(\sigma+{\rm id}_{\mathfrak g})$.  
From $\sigma\circ\theta=\theta\circ\sigma$, we have 
$$\mathfrak p=\mathfrak p\cap\mathfrak h\oplus\mathfrak p\cap\mathfrak q\quad\,\,{\rm and}\quad\,\,
\mathfrak k=\mathfrak k\cap\mathfrak h\oplus\mathfrak k\cap\mathfrak q.$$
Take a maximal abelian subspace $\mathfrak b$ of $\mathfrak p\cap\mathfrak q$ and 
a maximal abelian subsapce $\mathfrak a$ of $\mathfrak p$ including $\mathfrak b$.  
For $\beta\in\mathfrak b^{\ast}$, we define $\mathfrak p_{\beta}$ and $\mathfrak k_{\beta}$ by 
$$\mathfrak k_{\beta}:=\{v\in\mathfrak k\,\vert\,{\rm ad}(Z)^2(v)=-\beta(Z)^2v\,\,(\forall\,Z\in\mathfrak b)\}$$
and 
$$\mathfrak p_{\beta}:=\{v\in\mathfrak p\,\vert\,{\rm ad}(Z)^2(v)=-\beta(Z)^2v\,\,(\forall\,Z\in\mathfrak b)\}.$$
Also, we define $\triangle_{\mathfrak b}(\subset\mathfrak b^{\ast})$ by 
$$\triangle_{\mathfrak b}:=\{\beta\in\mathfrak b^{\ast}\,\vert\,\mathfrak p_{\beta}\not=\{0\}\,\,\},$$
which is the root system.  Let $(\triangle_{\mathfrak b})_+$ be the positive root subsystem of 
$\triangle_{\mathfrak b}$ with respect to some lexicographic ordering of $\mathfrak b^{\ast}$.  
Then we have 
$$\begin{array}{l}
\displaystyle{\mathfrak k=\mathfrak z_{\mathfrak k}(\mathfrak b)\oplus\left(
\mathop{\oplus}_{\beta\in(\triangle_{\mathfrak b})_+}\mathfrak k_{\beta}\right),}\\
\displaystyle{\mathfrak p=\mathfrak z_{\mathfrak p}(\mathfrak b)\oplus
\left(\mathop{\oplus}_{\beta\in(\triangle_{\mathfrak b})_+}\mathfrak p_{\beta}\right),}\\
\displaystyle{\mathfrak h=\mathfrak z_{\mathfrak h}(\mathfrak b)\oplus
\left(\mathop{\oplus}_{\beta\in(\triangle_{\mathfrak b})_+}\mathfrak h_{\beta}\right),}\\
\displaystyle{\mathfrak q=\mathfrak z_{\mathfrak q}(\mathfrak b)\oplus
\left(\mathop{\oplus}_{\beta\in(\triangle_{\mathfrak b})_+}\mathfrak q_{\beta}\right),}
\end{array}$$
where $\mathfrak z_{\bullet}(\mathfrak b)$ is the cetralizer of $\mathfrak b$ in $(\bullet)$.  
Set 
\begin{align*}
&\Sigma_{\mathfrak b}:={\rm Exp}_o(\mathfrak b),\,\,
\Sigma^d_{\mathfrak b}:={\rm Exp}_o(\sqrt{-1}\mathfrak b),\,\,
\Sigma^{\Bbb C}_{\mathfrak b}:={\rm Exp}_o(\mathfrak b^{\Bbb C}),\\
&\Sigma_{\mathfrak a}:={\rm Exp}_o(\mathfrak a),\,\,
\Sigma^d_{\mathfrak a}:={\rm Exp}_o(\sqrt{-1}\mathfrak a)\,\,\,{\rm and}\,\,\,
\Sigma^{\Bbb C}_{\mathfrak a}:={\rm Exp}_o(\mathfrak a^{\Bbb C}).
\end{align*}
Note that $\Sigma_{\mathfrak a}$ (resp. $\Sigma^d_{\mathfrak a}$) is included by $M$ (resp. $M^d$) because 
$\sqrt{-1}\mathfrak p$ is identified with $T_o(M^d)$.  
Set $H^d:=\exp((\mathfrak h\cap\mathfrak k)\oplus\sqrt{-1}(\mathfrak h\cap\mathfrak p))$, 
$\theta^d:=\theta^{\Bbb C}\vert_{{\mathfrak g}^d}$, $\sigma^d:=\sigma^{\Bbb C}\vert_{{\mathfrak g}^d}$, 
$L:={\rm Fix}(\sigma\circ\theta)$ and $L^d:={\rm Fix}(\sigma^d\circ\theta^d)$.  
The normal umbrella ${\rm Exp}_o(T^{\perp}_o(H^d\cdot o))$ of 
$H^d\cdot o$ in $M^d$ is isometric to the symmetric space $L^d/H\cap K$ and that 
the normal umbrella ${\rm Exp}_o(T^{\perp}_o(H\cdot o)\cap T_oM)$ of $H\cdot o$ in $M$ is 
isometric to the symmetric space $L/H\cap K$ (see \cite{Koi1,Koi3,Koi4}), where 
$T^{\perp}_o(H^d\cdot o)$ is the normal space of $H^d\cdot o$ in $M^d$ at $o$.  
It is shown that $T_o(L^d/H\cap K)=\sqrt{-1}(\mathfrak p\cap\mathfrak q)$, 
$T_o(H^d\cdot o)=\sqrt{-1}(\mathfrak p\cap\mathfrak h)$ and that all orbits of $G\curvearrowright M^{\Bbb C}$ 
meet $\Sigma^d_{\mathfrak a}$ orthogonally (see \cite{Koi1,Koi3,Koi4}).  
Denote by $H_p$ the isotropy group of $H\curvearrowright M$ at $p(\in M)$ and $\mathfrak h_p$ the Lie algebra of 
$H_p$.  Also, let $\mathfrak h_p^{\perp}$ be the orthogonal complement of $\mathfrak h_p$ in $\mathfrak h$.  
Set $H_p^{\perp}:=\{\exp\,X\,|\,X\in\mathfrak h_p^{\perp}\}$.  
The group $H_p$ acts on the normal umbrella $M^d_p$.  
First we prove the following fact.  

\vspace{0.5truecm}

\noindent
{\bf Lemma 4.1.} {\sl For $p\in M$ and $q\in M^d_p$, we have 
$$H\cdot q=\mathop{\cup}_{X\in\mathfrak h_p^{\perp}}H_{({\rm Exp}_o\,X)\cdot p}\cdot(({\rm Exp}_o\,X)\cdot q).$$
Hence $H\cdot q$ has the structure of the fiber bundle over $H\cdot p$ with the standard fibre 
$H_p\cdot q$ and the structure group $H_p$.}

\vspace{0.5truecm}

\noindent
{\it Proof.} Since $({\rm Exp}_oX)\cdot M_p^d=M_{({\rm Exp}_oX)\cdot p}^d$ holds for any 
$X\in\mathfrak h_p^{\perp}$, the first relation is derived.  
For any $X\in({\rm Exp}_o,X)$, $H_{({\rm Exp}_oX)\cdot p}$ is conjugate to $H_p$ and 
$H_{({\rm Exp}_oX)\cdot p}\cdot(({\rm Exp}_oX)\cdot q)$ is diffeomorphic to $H_p\cdot q$.  
Hence the second-half part of the statement is derived.  \qed

\vspace{0.5truecm}

\centerline{
\unitlength 0.1in
\begin{picture}( 81.8200, 19.9000)(-31.1000,-23.1000)
%
\special{pn 8}%
\special{pa 1288 600}%
\special{pa 776 1114}%
\special{pa 776 2310}%
\special{pa 1288 1798}%
\special{pa 1288 1798}%
\special{pa 1288 600}%
\special{fp}%
%
\special{pn 20}%
\special{sh 1}%
\special{ar 1048 1490 10 10 0  6.28318530717959E+0000}%
\special{sh 1}%
\special{ar 1048 1490 10 10 0  6.28318530717959E+0000}%
%
\special{pn 8}%
\special{pa 818 1490}%
\special{pa 2508 1490}%
\special{fp}%
%
\special{pn 8}%
\special{pa 742 1490}%
\special{pa 630 1490}%
\special{fp}%
%
\special{pn 8}%
\special{pa 810 1182}%
\special{pa 2498 1182}%
\special{fp}%
%
\special{pn 8}%
\special{pa 750 1182}%
\special{pa 640 1182}%
\special{fp}%
%
\special{pn 8}%
\special{pa 750 1764}%
\special{pa 640 1764}%
\special{fp}%
%
\special{pn 8}%
\special{pa 818 1764}%
\special{pa 2508 1764}%
\special{fp}%
\put(16.1000,-7.0000){\makebox(0,0)[lb]{$M^d_p$}}%
\put(23.9000,-19.5000){\makebox(0,0)[lt]{$({\rm Exp}_oX)\cdot p$}}%
%
\special{pn 20}%
\special{sh 1}%
\special{ar 1842 1490 10 10 0  6.28318530717959E+0000}%
\special{sh 1}%
\special{ar 1842 1490 10 10 0  6.28318530717959E+0000}%
\put(10.1000,-15.4000){\makebox(0,0)[lt]{$p$}}%
\put(27.7000,-12.0000){\makebox(0,0)[lb]{$M$}}%
%
\special{pn 20}%
\special{sh 1}%
\special{ar 1840 1180 10 10 0  6.28318530717959E+0000}%
\special{sh 1}%
\special{ar 1840 1180 10 10 0  6.28318530717959E+0000}%
%
\special{pn 20}%
\special{sh 1}%
\special{ar 1050 1180 10 10 0  6.28318530717959E+0000}%
\special{sh 1}%
\special{ar 1050 1180 10 10 0  6.28318530717959E+0000}%
\put(28.1000,-8.3000){\makebox(0,0)[lb]{$({\rm Exp}_oX)\cdot q$}}%
\put(10.1000,-11.3000){\makebox(0,0)[lb]{$q$}}%
%
\special{pn 8}%
\special{pa 810 910}%
\special{pa 960 1260}%
\special{dt 0.045}%
\special{sh 1}%
\special{pa 960 1260}%
\special{pa 952 1192}%
\special{pa 940 1212}%
\special{pa 916 1208}%
\special{pa 960 1260}%
\special{fp}%
\put(28.3000,-10.0000){\makebox(0,0)[lb]{$H\cdot q$}}%
\put(18.4000,-4.9000){\makebox(0,0)[rb]{$H_{({\rm Exp}_oX)\cdot p}\cdot(({\rm Exp}_oX)\cdot q)$}}%
%
\special{pn 8}%
\special{ar 3338 1584 942 284  3.4058527 5.1840909}%
%
\special{pn 8}%
\special{pa 2750 1150}%
\special{pa 2330 1490}%
\special{dt 0.045}%
\special{sh 1}%
\special{pa 2330 1490}%
\special{pa 2394 1464}%
\special{pa 2372 1456}%
\special{pa 2370 1434}%
\special{pa 2330 1490}%
\special{fp}%
\put(29.4000,-13.7000){\makebox(0,0)[lt]{in fact}}%
%
\special{pn 8}%
\special{pa 3734 1338}%
\special{pa 3834 1348}%
\special{fp}%
\special{sh 1}%
\special{pa 3834 1348}%
\special{pa 3768 1322}%
\special{pa 3780 1342}%
\special{pa 3766 1362}%
\special{pa 3834 1348}%
\special{fp}%
%
\special{pn 8}%
\special{pa 4166 1244}%
\special{pa 4198 1240}%
\special{pa 4230 1236}%
\special{pa 4262 1234}%
\special{pa 4294 1234}%
\special{pa 4324 1232}%
\special{pa 4356 1232}%
\special{pa 4388 1230}%
\special{pa 4420 1232}%
\special{pa 4452 1234}%
\special{pa 4484 1236}%
\special{pa 4516 1238}%
\special{pa 4548 1240}%
\special{pa 4580 1244}%
\special{pa 4612 1248}%
\special{pa 4644 1254}%
\special{pa 4676 1258}%
\special{pa 4706 1264}%
\special{pa 4738 1272}%
\special{pa 4770 1278}%
\special{pa 4800 1286}%
\special{pa 4832 1294}%
\special{pa 4862 1304}%
\special{pa 4892 1312}%
\special{pa 4922 1324}%
\special{pa 4952 1334}%
\special{pa 4982 1346}%
\special{pa 5012 1360}%
\special{pa 5040 1374}%
\special{pa 5068 1388}%
\special{pa 5072 1390}%
\special{sp}%
%
\special{pn 8}%
\special{pa 3830 1620}%
\special{pa 3842 1590}%
\special{pa 3856 1562}%
\special{pa 3870 1532}%
\special{pa 3884 1504}%
\special{pa 3902 1478}%
\special{pa 3918 1450}%
\special{pa 3938 1424}%
\special{pa 3958 1400}%
\special{pa 3980 1376}%
\special{pa 4002 1354}%
\special{pa 4026 1332}%
\special{pa 4050 1312}%
\special{pa 4076 1292}%
\special{pa 4102 1276}%
\special{pa 4130 1260}%
\special{pa 4160 1248}%
\special{pa 4176 1242}%
\special{sp}%
%
\special{pn 8}%
\special{pa 3820 1620}%
\special{pa 3852 1618}%
\special{pa 3884 1614}%
\special{pa 3916 1612}%
\special{pa 3948 1610}%
\special{pa 3980 1608}%
\special{pa 4012 1608}%
\special{pa 4044 1608}%
\special{pa 4076 1608}%
\special{pa 4108 1610}%
\special{pa 4140 1612}%
\special{pa 4172 1616}%
\special{pa 4204 1618}%
\special{pa 4236 1622}%
\special{pa 4268 1626}%
\special{pa 4298 1632}%
\special{pa 4330 1636}%
\special{pa 4362 1642}%
\special{pa 4392 1650}%
\special{pa 4424 1656}%
\special{pa 4456 1664}%
\special{pa 4486 1672}%
\special{pa 4516 1682}%
\special{pa 4548 1690}%
\special{pa 4578 1702}%
\special{pa 4608 1712}%
\special{pa 4638 1724}%
\special{pa 4666 1736}%
\special{pa 4696 1750}%
\special{pa 4724 1764}%
\special{pa 4728 1766}%
\special{sp}%
%
\special{pn 8}%
\special{pa 4720 1760}%
\special{pa 4732 1730}%
\special{pa 4746 1702}%
\special{pa 4760 1672}%
\special{pa 4774 1644}%
\special{pa 4792 1618}%
\special{pa 4808 1590}%
\special{pa 4828 1564}%
\special{pa 4848 1540}%
\special{pa 4870 1516}%
\special{pa 4892 1494}%
\special{pa 4916 1472}%
\special{pa 4940 1452}%
\special{pa 4966 1432}%
\special{pa 4992 1416}%
\special{pa 5020 1400}%
\special{pa 5050 1388}%
\special{pa 5066 1382}%
\special{sp}%
%
\special{pn 20}%
\special{sh 1}%
\special{ar 4290 1530 10 10 0  6.28318530717959E+0000}%
\special{sh 1}%
\special{ar 4290 1530 10 10 0  6.28318530717959E+0000}%
%
\special{pn 8}%
\special{ar 4430 1440 280 110  0.0000000 6.2831853}%
%
\special{pn 20}%
\special{sh 1}%
\special{ar 4530 1340 10 10 0  6.28318530717959E+0000}%
\special{sh 1}%
\special{ar 4530 1340 10 10 0  6.28318530717959E+0000}%
\put(43.5000,-17.5000){\makebox(0,0)[rt]{$M$}}%
%
\special{pn 8}%
\special{pa 4430 1140}%
\special{pa 4290 1520}%
\special{dt 0.045}%
\special{sh 1}%
\special{pa 4290 1520}%
\special{pa 4332 1464}%
\special{pa 4308 1470}%
\special{pa 4294 1452}%
\special{pa 4290 1520}%
\special{fp}%
\put(44.1000,-11.0000){\makebox(0,0)[lb]{$p$}}%
\put(46.7000,-8.7000){\makebox(0,0)[lb]{$({\rm Exp}_oX)\cdot p$}}%
\put(47.8000,-11.7000){\makebox(0,0)[lb]{$H\cdot p$}}%
\put(44.7000,-14.0000){\makebox(0,0)[lt]{$o$}}%
%
\special{pn 20}%
\special{sh 1}%
\special{ar 4430 1440 10 10 0  6.28318530717959E+0000}%
\special{sh 1}%
\special{ar 4430 1440 10 10 0  6.28318530717959E+0000}%
\put(9.1000,-8.9000){\makebox(0,0)[rb]{$H_p\cdot q$}}%
%
\special{pn 8}%
\special{pa 2080 600}%
\special{pa 1568 1114}%
\special{pa 1568 2310}%
\special{pa 2080 1798}%
\special{pa 2080 1798}%
\special{pa 2080 600}%
\special{fp}%
\put(22.1000,-5.1000){\makebox(0,0)[lb]{$M^d_{({\rm Exp}_oX)\cdot p}$}}%
%
\special{pn 8}%
\special{ar 1840 1470 100 286  0.5819051 6.2831853}%
%
\special{pn 8}%
\special{pa 1450 500}%
\special{pa 1770 1250}%
\special{dt 0.045}%
\special{sh 1}%
\special{pa 1770 1250}%
\special{pa 1762 1182}%
\special{pa 1750 1202}%
\special{pa 1726 1198}%
\special{pa 1770 1250}%
\special{fp}%
%
\special{pn 8}%
\special{pa 1590 640}%
\special{pa 1200 840}%
\special{dt 0.045}%
\special{sh 1}%
\special{pa 1200 840}%
\special{pa 1268 828}%
\special{pa 1248 816}%
\special{pa 1250 792}%
\special{pa 1200 840}%
\special{fp}%
%
\special{pn 8}%
\special{pa 4770 910}%
\special{pa 4540 1330}%
\special{dt 0.045}%
\special{sh 1}%
\special{pa 4540 1330}%
\special{pa 4590 1282}%
\special{pa 4566 1284}%
\special{pa 4554 1262}%
\special{pa 4540 1330}%
\special{fp}%
%
\special{pn 8}%
\special{pa 4900 1180}%
\special{pa 4670 1380}%
\special{dt 0.045}%
\special{sh 1}%
\special{pa 4670 1380}%
\special{pa 4734 1352}%
\special{pa 4710 1346}%
\special{pa 4708 1322}%
\special{pa 4670 1380}%
\special{fp}%
%
\special{pn 8}%
\special{ar 2790 1240 430 330  3.1798093 3.2113883}%
\special{ar 2790 1240 430 330  3.3061251 3.3377041}%
\special{ar 2790 1240 430 330  3.4324409 3.4640199}%
\special{ar 2790 1240 430 330  3.5587567 3.5903356}%
\special{ar 2790 1240 430 330  3.6850725 3.7166514}%
\special{ar 2790 1240 430 330  3.8113883 3.8429672}%
\special{ar 2790 1240 430 330  3.9377041 3.9692830}%
\special{ar 2790 1240 430 330  4.0640199 4.0955988}%
\special{ar 2790 1240 430 330  4.1903356 4.2219146}%
\special{ar 2790 1240 430 330  4.3166514 4.3482304}%
\special{ar 2790 1240 430 330  4.4429672 4.4745462}%
\special{ar 2790 1240 430 330  4.5692830 4.6008620}%
\special{ar 2790 1240 430 330  4.6955988 4.7123890}%
%
\special{pn 8}%
\special{pa 2360 1220}%
\special{pa 2350 1280}%
\special{fp}%
\special{sh 1}%
\special{pa 2350 1280}%
\special{pa 2382 1218}%
\special{pa 2360 1228}%
\special{pa 2342 1212}%
\special{pa 2350 1280}%
\special{fp}%
%
\special{pn 8}%
\special{ar 2760 1160 920 390  3.2484567 3.2667773}%
\special{ar 2760 1160 920 390  3.3217392 3.3400598}%
\special{ar 2760 1160 920 390  3.3950216 3.4133422}%
\special{ar 2760 1160 920 390  3.4683040 3.4866247}%
\special{ar 2760 1160 920 390  3.5415865 3.5599071}%
\special{ar 2760 1160 920 390  3.6148689 3.6331895}%
\special{ar 2760 1160 920 390  3.6881514 3.7064720}%
\special{ar 2760 1160 920 390  3.7614338 3.7797544}%
\special{ar 2760 1160 920 390  3.8347163 3.8530369}%
\special{ar 2760 1160 920 390  3.9079987 3.9263193}%
\special{ar 2760 1160 920 390  3.9812811 3.9996018}%
\special{ar 2760 1160 920 390  4.0545636 4.0728842}%
\special{ar 2760 1160 920 390  4.1278460 4.1461666}%
\special{ar 2760 1160 920 390  4.2011285 4.2194491}%
\special{ar 2760 1160 920 390  4.2744109 4.2927315}%
\special{ar 2760 1160 920 390  4.3476934 4.3660140}%
\special{ar 2760 1160 920 390  4.4209758 4.4392964}%
\special{ar 2760 1160 920 390  4.4942582 4.5125789}%
\special{ar 2760 1160 920 390  4.5675407 4.5858613}%
\special{ar 2760 1160 920 390  4.6408231 4.6591437}%
\special{ar 2760 1160 920 390  4.7141056 4.7200989}%
%
\special{pn 8}%
\special{pa 1850 1110}%
\special{pa 1830 1180}%
\special{fp}%
\special{sh 1}%
\special{pa 1830 1180}%
\special{pa 1868 1122}%
\special{pa 1846 1130}%
\special{pa 1830 1110}%
\special{pa 1830 1180}%
\special{fp}%
%
\special{pn 8}%
\special{ar 2070 530 510 270  6.2831853 6.3139545}%
\special{ar 2070 530 510 270  6.4062622 6.4370315}%
\special{ar 2070 530 510 270  6.5293392 6.5601084}%
\special{ar 2070 530 510 270  6.6524161 6.6831853}%
\special{ar 2070 530 510 270  6.7754930 6.8062622}%
\special{ar 2070 530 510 270  6.8985699 6.9293392}%
\special{ar 2070 530 510 270  7.0216468 7.0524161}%
\special{ar 2070 530 510 270  7.1447238 7.1754930}%
\special{ar 2070 530 510 270  7.2678007 7.2985699}%
\special{ar 2070 530 510 270  7.3908776 7.4216468}%
\special{ar 2070 530 510 270  7.5139545 7.5447238}%
\special{ar 2070 530 510 270  7.6370315 7.6678007}%
\special{ar 2070 530 510 270  7.7601084 7.7908776}%
%
\special{pn 8}%
\special{pa 2060 800}%
\special{pa 2000 800}%
\special{fp}%
\special{sh 1}%
\special{pa 2000 800}%
\special{pa 2068 820}%
\special{pa 2054 800}%
\special{pa 2068 780}%
\special{pa 2000 800}%
\special{fp}%
%
\special{pn 8}%
\special{pa 2350 1920}%
\special{pa 1840 1500}%
\special{dt 0.045}%
\special{sh 1}%
\special{pa 1840 1500}%
\special{pa 1880 1558}%
\special{pa 1882 1534}%
\special{pa 1904 1528}%
\special{pa 1840 1500}%
\special{fp}%
%
\special{pn 8}%
\special{ar 1050 1470 110 290  0.2265960 6.2297479}%
\end{picture}%
\hspace{8.5truecm}}

\vspace{0.5truecm}

\centerline{{\bf Figure 4.}}

\vspace{0.5truecm}

\noindent
{\bf Lemma 4.2} {\sl Let $q\in\Psi_p\left(T_p^{\perp}(H\cdot p)\right)$ 
and denote by 
${\rm Hol}^{\perp}_{-\Psi_p^{-1}(q)}(H\cdot p)$ the normal holonomy bundle of the submanifold $H\cdot p$ in $M$ 
through $-\Psi_p^{-1}(q)$.  Then we have 
$$H\cdot q=\Psi\left({\rm Hol}^{\perp}_{-\Psi_p^{-1}(q)}(H\cdot p)\right).$$
}

\vspace{0.5truecm}

\noindent
{\it Proof.} It is clear that 
$\{({\rm Exp}_oX)\cdot p\,|\,X\in\mathfrak h_p^{\perp}\}=H\cdot p$.  
Since $H\curvearrowright G/K$ is hyperpolar, 
${\rm Exp}_{({\rm Exp}_oX)\cdot p}$\newline
$\displaystyle{\left(T^{\perp}_{({\rm Exp}_oX)\cdot p}(H\cdot p)\right)}$ 
is totally geodesic in $M$.  From this fact, 
we can show that the orbit $H_{({\rm Exp}_oX)\cdot p}$\newline
$(({\rm Exp}_oX)\cdot q)$ is equal to 
the image of the fibre of the normal holonomy bundle 
${\rm Hol}^{\perp}_{-\Psi_p^{-1}(q)}(H\cdot p)$ over $({\rm Exp}_oX)\cdot p$ by 
$\Psi_{({\rm Exp}_oX)\cdot p}$.  
Hence it follows from Lemma 4.1 that $H\cdot q$ is described as in the statement.  \qed

\section{Special Lagrangian submanifolds in complexified symmetric spaces}
Let $(G^{\mathbb C}/K^{\mathbb C},J_0,\omega_{\psi_f},\Omega_0)$ be the Calabi-Yau manifold stated in Section 2, 
where $f$ is as in (i) of Theorem 2.2.  
As in the previous section, set $M:=G\cdot o(=G/K),\,\,M^{\Bbb C}:=G^{\Bbb C}/K^{\Bbb C}$, 
$M^d:=G^d/K(=G^d\cdot o)$ and $M^d_p:=\Psi_p(T_p(G\cdot o))$.  
Let $H$ be a symmetric subgroup of $G$ and $r$ be the cohomogeneity of the Hermann action $H\curvearrowright G/K$.  
The naturally extended action of $H$ on $(M^{\mathbb C},J_0,\omega_{\psi_f},\Omega_0)$ is a Hamiltonian action.  
Denote by $\mu_{\psi_f}$ the moment map of this Hamiltonian action.  Let $Z(\mathfrak h^{\ast})$ be the center of 
$\mathfrak g^{\ast}$, that is, 
$$Z(\mathfrak h^{\ast}):=\{X\in\mathfrak g^{\ast}\,|\,{\rm Ad}^{\ast}(h)(X)=X\,\,(\forall\,h\in H)\},$$
where ${\rm Ad}^{\ast}$ denotes the coadjoint representation of $H$.  
It is clear that $\mu_{\psi_f}^{-1}(c)$ is $H$-invariant if and only if $c$ belongs to $Z(\mathfrak h^{\ast})$.  
According to Proposition 2.5 of \cite{HS}, the following fact holds.  

\vspace{0.5truecm}

\noindent
{\bf Proposition 5.1(\cite{HS}).} 
{\sl Assume that $L$ is a $H$-invariant connected isotropic submanifold in \newline
$(M^{\mathbb C},J_0,\omega_{\psi_f},\Omega_0)$, where ``isotropic'' means that $\omega_{\psi_f}(TL,TL)=0$ holds.  
Then $L\subset\mu_{\psi_f}^{-1}(c)$ holds for some $c\in Z(\mathfrak h^{\ast})$.}

\vspace{0.5truecm}

In the method of the proof of Proposition 2.6 of \cite{HS}, we can show the following fact.  

\vspace{0.5truecm}

\noindent
{\bf Proposition 5.2.} {\sl Let $L$ be a $H$-invariant connected submanifold in $M^{\mathbb C}$ and 
$r_0$ be the cohomogeneity of the action $H\curvearrowright L$.  Assume that $L\subset\mu_{\psi_f}^{-1}(c)$ for some 
$c\in Z(\mathfrak h^{\ast})$ and that there exists a $r_0$-dimensional isotropic submanifold $L_0$ in 
$(M^{\mathbb C},J_0,\omega_{\psi_f},\Omega_0)$ satisfying the following conditions:

(i) $L_0\subset L$,

(ii) $L_0$ is transversal to the principal orbits of the action $H\curvearrowright L$,

(iii) $H\cdot L_0=L$,

\noindent
Then $L$ also is an isotropic submanifold in $(M^{\mathbb C},J_0,\omega_{\psi_f},\Omega_0)$.}

\vspace{0.5truecm}

\noindent
{\it Proof.} Take any $X\in\mathfrak h$ and any $Y\in T_pL$.  
From $L\subset\mu_{\psi_f}^{-1}(c)$, we have $d(\mu_{\psi_f})_p(Y)=0$.  
On the other hand, we have $(d(\mu_{\psi_f})_p(Y))(X)=$\newline
$(\omega_{\psi_f})_p(Y,X^{\ast}_p)$, 
where $X^{\ast}$ is the vector field on $M^{\mathbb C}$ associated to 
the one-parameter transformation group $\{\exp\,tX\}_{t\in\mathbb R}$ of $M^{\mathbb C}$ 
($\exp\,:\,$ the exponential map of $H$).  Hence we have $(\omega_{\psi_f})_p(Y,X^{\ast}_p)=0$.  Therefore, 
it follows from the arbitrariness of $X$ and $Y$ that $(\omega_{\psi_f})_p(T_pL,T_p(H\cdot p))=0$.  
Also, since $L_0$ is isotropic, we have 
$(\omega_{\psi_f})_p(T_pL_0,T_pL_0)=0$.  Hence we obtain $(\omega_{\psi_f})_p(T_pL,T_pL)=0$.  
Therefore, it follows from the arbitariness of $p$ that $L$ is isotropic.  \qed

\vspace{0.5truecm}

By Proposition 2.4 of \cite{HS}, we can show the following fact.  

\vspace{0.5truecm}

\noindent
{\bf Proposition 5.3.} {\sl Let $L$ be a $n$-dimensional connected submanifold in \newline
$(M^{\mathbb C},J_0,\omega_{\psi_f},\Omega_0)$.  
Then $L$ is a special Lagrangian submanifold of phase $\theta$ if and only if \newline
$\omega_{\psi_f}|_{TL\times TL}=0$ and ${\rm Im}\left(e^{\sqrt{-1}\theta}\Omega_0|_{(TL)^n}\right)=0$.}

\vspace{0.5truecm}

Let $f$ be as in $(2.8)$.  We give constructions of special Lagrangian submanifolds in the Calabi-Yau manifold 
$(M^{\mathbb C},J_0,\omega_{\psi_f},\Omega_0)$.  
Let $U$ be the open subset of $M$ of all regular points of $H\curvearrowright M$.  
Then, as stated in the proof of Lemma 3.2, 
$\displaystyle{\widetilde{\Sigma}:=\mathop{\amalg}_{p\in U}\Psi_p\left(T_p^{\perp}(H\cdot p)\right)}$ is 
an open subset of $\mu_{\psi_f}^{-1}(0)$.  
Since $H\curvearrowright M$ is a Hermann action, it is hyperpolar (see Subsections 3.1 in [HPTT]).  
Hence the principal orbit $H\cdot p_0$ ($p_0\in U$) is an equifocal submanifold in $M$ and 
its section $\displaystyle{\Sigma:={\rm Exp}_{p_0}\left(T^{\perp}_{p_0}(H\cdot p_0)\right)}$ is an $r$-dimensional 
flat torus $T^r=S^1\times\cdots\times S^1$ ($r$-times) embedded totally geodesically into $M$.  
Without loss of generality, we may assume that $\Sigma$ passes through $o$.  
Let $C$ be the component of $U\cap\Sigma$ containing $p_0$.  Then we have $H\cdot C=U$.  
Set $\displaystyle{\widehat{\Sigma}:=\mathop{\amalg}_{p\in U\cap\Sigma}\Psi_p\left(T^{\perp}_p(H\cdot p)\right)}$.  
It is clear that $\widehat{\Sigma}$ is a dense open subset of 
$\displaystyle{\Sigma^{\mathbb C}:=\mathop{\amalg}_{p\in\Sigma}\Psi_p(T_p\Sigma)\,
(\approx(T^r)^{\mathbb C}=S^1_{\mathbb C}\times\cdots\times S^1_{\mathbb C})}$, where 
$S^1_{\mathbb C}\times\cdots\times S^1_{\mathbb C})$ denotes the $r$-times of $S^1_{\mathbb C}$'s.  
We identify $T_o\Sigma\,(\subset\mathfrak p)$ and $T_o(\Sigma^{\mathbb C})\,(\subset\mathfrak p^{\mathbb C})$ with 
$\mathbb R^r$ and $\mathbb C^r$, respectively.  
Let $\tau_i:I_i\to \mathbb C$ ($i=1,\cdots,r$) be regular curves, where $I_i$ is an open interval.  Define 
an immersion $\tau:I_1\times\cdots\times I_r\hookrightarrow\mathbb C^r$ by $\tau:=\tau_1\times\cdots\times\tau_r$.  
Set $\overline{\tau}:={\rm Exp}_o\circ\tau 
(:I_1\times\cdots\times I_r\to S^1_{\mathbb C}\times\cdots\times S^1_{\mathbb C}(=\Sigma^{\mathbb C}))$.  
Assume that $(L_{\tau})_0:=\overline{\tau}(I_1\times\cdots\times I_r)$ is included by $\widehat{\Sigma}$.  
It is clear that $(L_{\tau})_0$ is an isotropic submanifold in $\Sigma^{\mathbb C}$ (hence in 
$(M^{\mathbb C},J_0,\omega_{\psi_f},\Omega_0)$).  
Set $L_{\tau}:=H\cdot(L_{\tau})_0$.  
For any $p\in U$ and any $\displaystyle{q\in\Psi_p\left(T_p^{\perp}(H\cdot p)\right)\,
(\subset\widetilde{\Sigma})}$, since $H\cdot p$ is an equifocal submanifold in $M$, the normal connection of 
the submanifold $H\cdot p$ in $M$ is flat and hence 
the (restricted) normal holonomy representation 
$$(H_{({\rm Exp}_oX)\cdot p})_0\curvearrowright T_{({\rm Exp}_oX)\cdot p}^{\perp}(H\cdot p)$$
is trivial, where 
$(H_{({\rm Exp}_oX)\cdot p})_0$ denotes the identity component of $H_{({\rm Exp}_oX)\cdot p}$.  
Hence the action 
$$(H_{({\rm Exp}_oX)\cdot p})_0\curvearrowright\Psi_{({\rm Exp}_oX)\cdot p}\left(T_{({\rm Exp}_oX)\cdot p}^{\perp}
(H\cdot p)\right)$$
also is trivial.  Therefore it follows from Lemma 4.1 that each component of $H\cdot q$ is diffeomorphic to 
$\displaystyle{\mathop{\cup}_{X\in\mathfrak h_p^{\perp}}({\rm Exp}_oX)\cdot q}$.  
From this fact, ${\rm dim}\,H\cdot q=n-r$ follows.  
Since $(L_{\tau})_0$ is included by $\widetilde{\Sigma}$, $L_{\tau}$ is an $n$-dimensional 
submanifold of cohomogeneity $r$ in $M^{\mathbb C}$.  By Proposition 5.2, $L_{\tau}$ is 
a Lagrangian submanifold.  

Here we shall explain that the cohomogeneity of the Hamiltonian action $H\curvearrowright M^{\mathbb C}$ is 
possible to be smaller than $(n+r)$.  
For $\displaystyle{\widehat q\in\Psi_p\left(T_p(H\cdot p)\right)\,(\subset\widetilde{\Sigma})}$, 
the (restricted) holonomy representation 
$$(H_{({\rm Exp}_oX)\cdot p})_0\curvearrowright T_{({\rm Exp}_oX)\cdot p}(H\cdot p)$$
of the Riemannian manifold $H\cdot p$ at $({\rm Exp}_oX)\cdot p$ is not 
necessarily trivial.  Hence the action 
$$(H_{({\rm Exp}_oX)\cdot p})_0\curvearrowright\Psi_{({\rm Exp}_oX)\cdot p}\left(T_{({\rm Exp}_oX)\cdot p}(H\cdot p)
\right)$$
also is not necessarily trivial.  
On the other hand, we can show that $H\cdot\widehat q$ is equal to the image of 
the holonomy bundle ${\rm Hol}_{\Psi_p^{-1}(\widehat q)}(H\cdot p)$ of $H\cdot p$ throught 
$\Psi_p^{-1}(\widehat q)$ by $\Psi$.  
From these facts, it follows that ${\rm dim}\,(H\cdot\widehat q)$ is possible to be larger than $(n-r)$.  
That is, the cohomogeneity of the action $H\curvearrowright M^{\mathbb C}$ is possible to be smaller than $(n+r)$.  
Set $\mathfrak b:=T_o\Sigma$, which is a maximal abelian subspace of $\mathfrak p\cap\mathfrak q$.  
Note that $\tau$ is rearded as a regular curve in $\mathfrak b^{\mathbb C}$ under the identification of 
$\mathbb C^r$ with $\mathfrak b^{\mathbb C}$.  
Let $\triangle_{\mathfrak b},\,(\triangle_{\mathfrak b})_+,\,\mathfrak k_{\beta},\,\mathfrak p_{\beta},\,
\mathfrak h_{\beta}$ and $\mathfrak q_{\beta}$ be as in Section 4.  Define $(\triangle_{\mathfrak b}^V)_+$ and 
$(\triangle_{\mathfrak b}^H)_+$ by 
$$(\triangle_{\mathfrak b}^V)_+:=\{\beta\in(\triangle_{\mathfrak b})_+\,|\,
\mathfrak p_{\beta}\cap\mathfrak q\not=\{0\}\,\}$$
and 
$$(\triangle_{\mathfrak b}^H)_+:=\{\beta\in(\triangle_{\mathfrak b})_+\,|\,
\mathfrak p_{\beta}\cap\mathfrak h\not=\{0\}\,\},$$
respectively.  Note that 
${\rm dim}(\mathfrak p_{\beta}\cap\mathfrak q)={\rm dim}(\mathfrak k_{\beta}\cap\mathfrak h)$ and 
${\rm dim}(\mathfrak p_{\beta}\cap\mathfrak h)={\rm dim}(\mathfrak k_{\beta}\cap\mathfrak q)$.  
Set $m_{\beta}^V:={\rm dim}(\mathfrak p_{\beta}\cap\mathfrak q)$ ($\beta\in(\triangle_{\mathfrak b}^V)_+$) and 
$m_{\beta}^H:={\rm dim}(\mathfrak p_{\beta}\cap\mathfrak h)$ ($\beta\in(\triangle_{\mathfrak b}^H)_+$).  
Let $\{X_{\beta,i}^V\,|\,i=1,\cdots,m_{\beta}^V\}$ be a basis of $\mathfrak k_{\beta}\cap\mathfrak h$ 
($\beta\in(\triangle_{\mathfrak b}^V)_+$) and 
$\{X_{\beta,i}^H\,|\,i=1,\cdots,m_{\beta}^H\}$ be a basis of $\mathfrak p_{\beta}\cap\mathfrak h$ 
($\beta\in(\triangle_{\mathfrak b}^V)_+$).  Also, let $Y_{\beta,i}^V$ be the element of 
$\mathfrak k_{\beta}\cap\mathfrak h$ such that ${\rm ad}(Z)(X_{\beta,i}^V)=\beta(Z)Y_{\beta,i}^V$ holds for any 
$Z\in\mathfrak b$.  Define a Killing vector field $(Y^V_{\beta,i})^{\ast}$ over $M^{\mathbb C}$ by 
$$(Y^V_{\beta,i})^{\ast}_p:=\left.\frac{d}{dt}\right|_{t=0}\exp(tY^V_{\beta,i})(p)
\quad(p\in M^{\mathbb C})$$
and 
a Killing vector field $(X^H_{\beta,i})^{\ast}$ over $M^{\mathbb C}$ by 
$$(X^H_{\beta,i})^{\ast}_p:=\left.\frac{d}{dt}\right|_{t=0}\exp(tX^H_{\beta,i})(p)
\quad(p\in M^{\mathbb C}).$$

\vspace{0.5truecm}

\centerline{
\unitlength 0.1in
\begin{picture}( 20.5700, 23.6000)( 14.1700,-25.9500)
%
\special{pn 8}%
\special{ar 2146 1748 512 154  6.2831853 6.2831853}%
\special{ar 2146 1748 512 154  0.0000000 3.1415927}%
%
\special{pn 8}%
\special{ar 2146 1748 512 154  3.1415927 3.1776829}%
\special{ar 2146 1748 512 154  3.2859536 3.3220438}%
\special{ar 2146 1748 512 154  3.4303145 3.4664047}%
\special{ar 2146 1748 512 154  3.5746754 3.6107656}%
\special{ar 2146 1748 512 154  3.7190363 3.7551265}%
\special{ar 2146 1748 512 154  3.8633972 3.8994874}%
\special{ar 2146 1748 512 154  4.0077581 4.0438483}%
\special{ar 2146 1748 512 154  4.1521190 4.1882092}%
\special{ar 2146 1748 512 154  4.2964799 4.3325701}%
\special{ar 2146 1748 512 154  4.4408408 4.4769310}%
\special{ar 2146 1748 512 154  4.5852017 4.6212919}%
\special{ar 2146 1748 512 154  4.7295626 4.7656528}%
\special{ar 2146 1748 512 154  4.8739235 4.9100137}%
\special{ar 2146 1748 512 154  5.0182844 5.0543746}%
\special{ar 2146 1748 512 154  5.1626453 5.1987355}%
\special{ar 2146 1748 512 154  5.3070062 5.3430964}%
\special{ar 2146 1748 512 154  5.4513671 5.4874573}%
\special{ar 2146 1748 512 154  5.5957280 5.6318182}%
\special{ar 2146 1748 512 154  5.7400889 5.7761791}%
\special{ar 2146 1748 512 154  5.8844498 5.9205400}%
\special{ar 2146 1748 512 154  6.0288107 6.0649009}%
\special{ar 2146 1748 512 154  6.1731716 6.2092618}%
%
\special{pn 8}%
\special{pa 2618 1166}%
\special{pa 2610 1198}%
\special{pa 2588 1220}%
\special{pa 2562 1240}%
\special{pa 2534 1256}%
\special{pa 2506 1268}%
\special{pa 2476 1278}%
\special{pa 2446 1290}%
\special{pa 2414 1298}%
\special{pa 2382 1304}%
\special{pa 2352 1312}%
\special{pa 2320 1318}%
\special{pa 2288 1322}%
\special{pa 2258 1328}%
\special{pa 2226 1332}%
\special{pa 2194 1334}%
\special{pa 2162 1336}%
\special{pa 2130 1338}%
\special{pa 2098 1340}%
\special{pa 2066 1340}%
\special{pa 2034 1340}%
\special{pa 2002 1338}%
\special{pa 1970 1338}%
\special{pa 1938 1336}%
\special{pa 1906 1332}%
\special{pa 1874 1330}%
\special{pa 1842 1324}%
\special{pa 1810 1320}%
\special{pa 1780 1314}%
\special{pa 1748 1306}%
\special{pa 1718 1298}%
\special{pa 1688 1286}%
\special{pa 1658 1274}%
\special{pa 1632 1256}%
\special{pa 1608 1234}%
\special{pa 1596 1206}%
\special{pa 1596 1204}%
\special{sp}%
%
\special{pn 8}%
\special{pa 1596 1204}%
\special{pa 1606 1174}%
\special{pa 1626 1150}%
\special{pa 1652 1130}%
\special{pa 1680 1114}%
\special{pa 1708 1102}%
\special{pa 1740 1090}%
\special{pa 1770 1080}%
\special{pa 1800 1072}%
\special{pa 1832 1066}%
\special{pa 1862 1058}%
\special{pa 1894 1052}%
\special{pa 1926 1048}%
\special{pa 1958 1042}%
\special{pa 1990 1038}%
\special{pa 2022 1036}%
\special{pa 2054 1034}%
\special{pa 2084 1032}%
\special{pa 2116 1030}%
\special{pa 2148 1030}%
\special{pa 2180 1030}%
\special{pa 2212 1032}%
\special{pa 2244 1032}%
\special{pa 2276 1034}%
\special{pa 2308 1038}%
\special{pa 2340 1040}%
\special{pa 2372 1044}%
\special{pa 2404 1048}%
\special{pa 2436 1056}%
\special{pa 2466 1064}%
\special{pa 2498 1072}%
\special{pa 2528 1082}%
\special{pa 2556 1096}%
\special{pa 2584 1112}%
\special{pa 2606 1136}%
\special{pa 2618 1164}%
\special{pa 2618 1166}%
\special{sp -0.045}%
%
\special{pn 8}%
\special{ar 1114 1740 522 1386  5.3365688 6.2831853}%
\special{ar 1114 1740 522 1386  0.0000000 0.5077668}%
%
\special{pn 8}%
\special{ar 1988 1748 670 1804  5.5116623 6.2831853}%
\special{ar 1988 1748 670 1804  0.0000000 0.3882855}%
%
\special{pn 8}%
\special{ar 1476 1808 670 1804  5.5116623 6.2831853}%
\special{ar 1476 1808 670 1804  0.0000000 0.3882855}%
%
\special{pn 13}%
\special{pa 1940 1894}%
\special{pa 2442 1928}%
\special{fp}%
\special{sh 1}%
\special{pa 2442 1928}%
\special{pa 2378 1904}%
\special{pa 2390 1924}%
\special{pa 2374 1944}%
\special{pa 2442 1928}%
\special{fp}%
%
\special{pn 13}%
\special{pa 1950 1338}%
\special{pa 2452 1356}%
\special{fp}%
\special{sh 1}%
\special{pa 2452 1356}%
\special{pa 2386 1334}%
\special{pa 2400 1354}%
\special{pa 2386 1374}%
\special{pa 2452 1356}%
\special{fp}%
%
\special{pn 13}%
\special{pa 2526 1270}%
\special{pa 2834 1072}%
\special{fp}%
\special{sh 1}%
\special{pa 2834 1072}%
\special{pa 2766 1092}%
\special{pa 2788 1102}%
\special{pa 2788 1126}%
\special{pa 2834 1072}%
\special{fp}%
%
\special{pn 13}%
\special{pa 2518 1860}%
\special{pa 2842 1722}%
\special{fp}%
\special{sh 1}%
\special{pa 2842 1722}%
\special{pa 2774 1730}%
\special{pa 2794 1744}%
\special{pa 2788 1766}%
\special{pa 2842 1722}%
\special{fp}%
%
\special{pn 13}%
\special{pa 1680 1816}%
\special{pa 2016 1962}%
\special{fp}%
\special{sh 1}%
\special{pa 2016 1962}%
\special{pa 1962 1918}%
\special{pa 1966 1942}%
\special{pa 1946 1954}%
\special{pa 2016 1962}%
\special{fp}%
%
\special{pn 13}%
\special{pa 1644 1270}%
\special{pa 1988 1406}%
\special{fp}%
\special{sh 1}%
\special{pa 1988 1406}%
\special{pa 1932 1364}%
\special{pa 1938 1386}%
\special{pa 1918 1400}%
\special{pa 1988 1406}%
\special{fp}%
%
\special{pn 20}%
\special{sh 1}%
\special{ar 2098 1192 10 10 0  6.28318530717959E+0000}%
\special{sh 1}%
\special{ar 2098 1192 10 10 0  6.28318530717959E+0000}%
%
\special{pn 20}%
\special{sh 1}%
\special{ar 2146 1766 10 10 0  6.28318530717959E+0000}%
\special{sh 1}%
\special{ar 2146 1766 10 10 0  6.28318530717959E+0000}%
\put(24.5200,-25.9500){\makebox(0,0)[lt]{$\,$}}%
%
\special{pn 8}%
\special{pa 2926 824}%
\special{pa 2442 890}%
\special{dt 0.045}%
\special{sh 1}%
\special{pa 2442 890}%
\special{pa 2512 900}%
\special{pa 2496 882}%
\special{pa 2506 860}%
\special{pa 2442 890}%
\special{fp}%
%
\special{pn 13}%
\special{pa 2944 2356}%
\special{pa 3150 2244}%
\special{fp}%
\special{sh 1}%
\special{pa 3150 2244}%
\special{pa 3082 2258}%
\special{pa 3102 2270}%
\special{pa 3100 2294}%
\special{pa 3150 2244}%
\special{fp}%
\put(31.8600,-22.2700){\makebox(0,0)[lt]{'s $:$ $(Y^V_{\beta,i})^{\ast}|_M$}}%
\put(28.3300,-7.3000){\makebox(0,0)[lb]{$H\cdot o$}}%
\put(29.6300,-8.9300){\makebox(0,0)[lb]{$H\cdot p$}}%
\put(21.8200,-17.9100){\makebox(0,0)[lb]{$o$}}%
%
\special{pn 20}%
\special{sh 1}%
\special{ar 2518 1860 10 10 0  6.28318530717959E+0000}%
\special{sh 1}%
\special{ar 2518 1860 10 10 0  6.28318530717959E+0000}%
\put(25.1700,-19.0200){\makebox(0,0)[lt]{$p$}}%
%
\special{pn 13}%
\special{pa 2228 1022}%
\special{pa 1838 1004}%
\special{fp}%
\special{sh 1}%
\special{pa 1838 1004}%
\special{pa 1904 1028}%
\special{pa 1892 1006}%
\special{pa 1906 988}%
\special{pa 1838 1004}%
\special{fp}%
%
\special{pn 13}%
\special{pa 2256 1604}%
\special{pa 1838 1578}%
\special{fp}%
\special{sh 1}%
\special{pa 1838 1578}%
\special{pa 1904 1602}%
\special{pa 1892 1580}%
\special{pa 1906 1562}%
\special{pa 1838 1578}%
\special{fp}%
%
\special{pn 8}%
\special{pa 3308 440}%
\special{pa 3474 440}%
\special{fp}%
\put(33.3500,-4.0500){\makebox(0,0)[lb]{$M$}}%
%
\special{pn 8}%
\special{pa 3298 440}%
\special{pa 3298 278}%
\special{fp}%
%
\special{pn 8}%
\special{pa 2814 662}%
\special{pa 2034 790}%
\special{dt 0.045}%
\special{sh 1}%
\special{pa 2034 790}%
\special{pa 2102 800}%
\special{pa 2086 782}%
\special{pa 2096 760}%
\special{pa 2034 790}%
\special{fp}%
\end{picture}%
\hspace{0truecm}}

\vspace{0.5truecm}

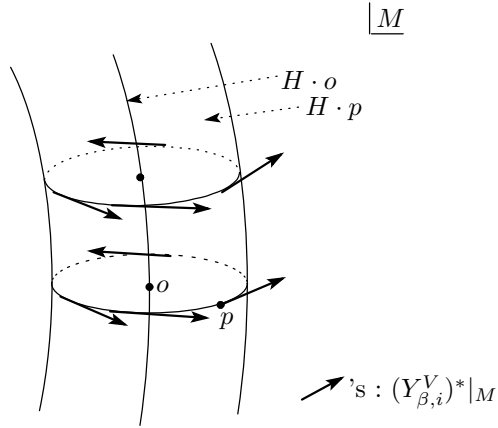
\centerline{{\bf Figure 5.}}

\vspace{0.5truecm}

\centerline{
\unitlength 0.1in
\begin{picture}( 19.7300, 22.2000)( 14.1700,-24.6000)
%
\special{pn 8}%
\special{ar 2138 1694 506 150  6.2831853 6.2831853}%
\special{ar 2138 1694 506 150  0.0000000 3.1415927}%
%
\special{pn 8}%
\special{ar 2138 1694 506 148  3.1415927 3.1782899}%
\special{ar 2138 1694 506 148  3.2883816 3.3250789}%
\special{ar 2138 1694 506 148  3.4351706 3.4718679}%
\special{ar 2138 1694 506 148  3.5819596 3.6186569}%
\special{ar 2138 1694 506 148  3.7287486 3.7654459}%
\special{ar 2138 1694 506 148  3.8755376 3.9122349}%
\special{ar 2138 1694 506 148  4.0223266 4.0590238}%
\special{ar 2138 1694 506 148  4.1691156 4.2058128}%
\special{ar 2138 1694 506 148  4.3159046 4.3526018}%
\special{ar 2138 1694 506 148  4.4626936 4.4993908}%
\special{ar 2138 1694 506 148  4.6094826 4.6461798}%
\special{ar 2138 1694 506 148  4.7562716 4.7929688}%
\special{ar 2138 1694 506 148  4.9030605 4.9397578}%
\special{ar 2138 1694 506 148  5.0498495 5.0865468}%
\special{ar 2138 1694 506 148  5.1966385 5.2333358}%
\special{ar 2138 1694 506 148  5.3434275 5.3801248}%
\special{ar 2138 1694 506 148  5.4902165 5.5269138}%
\special{ar 2138 1694 506 148  5.6370055 5.6737027}%
\special{ar 2138 1694 506 148  5.7837945 5.8204917}%
\special{ar 2138 1694 506 148  5.9305835 5.9672807}%
\special{ar 2138 1694 506 148  6.0773725 6.1140697}%
\special{ar 2138 1694 506 148  6.2241615 6.2608587}%
%
\special{pn 8}%
\special{pa 2606 1134}%
\special{pa 2598 1164}%
\special{pa 2576 1188}%
\special{pa 2550 1206}%
\special{pa 2522 1222}%
\special{pa 2492 1234}%
\special{pa 2462 1244}%
\special{pa 2430 1254}%
\special{pa 2400 1262}%
\special{pa 2368 1268}%
\special{pa 2338 1276}%
\special{pa 2306 1280}%
\special{pa 2274 1286}%
\special{pa 2242 1290}%
\special{pa 2210 1294}%
\special{pa 2178 1296}%
\special{pa 2146 1298}%
\special{pa 2114 1300}%
\special{pa 2082 1302}%
\special{pa 2050 1302}%
\special{pa 2018 1302}%
\special{pa 1988 1300}%
\special{pa 1956 1298}%
\special{pa 1924 1296}%
\special{pa 1892 1294}%
\special{pa 1860 1290}%
\special{pa 1828 1286}%
\special{pa 1796 1280}%
\special{pa 1764 1274}%
\special{pa 1734 1266}%
\special{pa 1704 1258}%
\special{pa 1674 1246}%
\special{pa 1644 1232}%
\special{pa 1618 1214}%
\special{pa 1600 1188}%
\special{pa 1596 1170}%
\special{sp}%
%
\special{pn 8}%
\special{pa 1596 1170}%
\special{pa 1604 1140}%
\special{pa 1626 1116}%
\special{pa 1652 1098}%
\special{pa 1680 1084}%
\special{pa 1710 1070}%
\special{pa 1740 1060}%
\special{pa 1770 1050}%
\special{pa 1802 1042}%
\special{pa 1832 1036}%
\special{pa 1864 1028}%
\special{pa 1896 1024}%
\special{pa 1928 1020}%
\special{pa 1958 1014}%
\special{pa 1990 1012}%
\special{pa 2022 1008}%
\special{pa 2054 1006}%
\special{pa 2086 1004}%
\special{pa 2118 1004}%
\special{pa 2150 1004}%
\special{pa 2182 1004}%
\special{pa 2214 1004}%
\special{pa 2246 1006}%
\special{pa 2278 1008}%
\special{pa 2310 1012}%
\special{pa 2342 1014}%
\special{pa 2374 1018}%
\special{pa 2406 1024}%
\special{pa 2436 1030}%
\special{pa 2468 1038}%
\special{pa 2498 1048}%
\special{pa 2528 1060}%
\special{pa 2556 1074}%
\special{pa 2584 1090}%
\special{pa 2602 1116}%
\special{pa 2606 1134}%
\special{sp -0.045}%
%
\special{pn 8}%
\special{ar 1116 1686 516 1334  5.3360424 6.2831853}%
\special{ar 1116 1686 516 1334  0.0000000 0.5076606}%
%
\special{pn 8}%
\special{ar 1982 1694 664 1738  5.5117918 6.2831853}%
\special{ar 1982 1694 664 1738  0.0000000 0.3886543}%
%
\special{pn 8}%
\special{ar 1476 1752 664 1738  5.5117918 6.2831853}%
\special{ar 1476 1752 664 1738  0.0000000 0.3879158}%
%
\special{pn 13}%
\special{pa 2148 1728}%
\special{pa 2156 1432}%
\special{fp}%
\special{sh 1}%
\special{pa 2156 1432}%
\special{pa 2134 1498}%
\special{pa 2154 1484}%
\special{pa 2174 1498}%
\special{pa 2156 1432}%
\special{fp}%
%
\special{pn 13}%
\special{pa 2654 1712}%
\special{pa 2664 1414}%
\special{fp}%
\special{sh 1}%
\special{pa 2664 1414}%
\special{pa 2642 1480}%
\special{pa 2662 1468}%
\special{pa 2682 1482}%
\special{pa 2664 1414}%
\special{fp}%
%
\special{pn 13}%
\special{pa 1642 1712}%
\special{pa 1650 1414}%
\special{fp}%
\special{sh 1}%
\special{pa 1650 1414}%
\special{pa 1628 1480}%
\special{pa 1648 1468}%
\special{pa 1668 1482}%
\special{pa 1650 1414}%
\special{fp}%
%
\special{pn 13}%
\special{pa 2102 1134}%
\special{pa 2074 822}%
\special{fp}%
\special{sh 1}%
\special{pa 2074 822}%
\special{pa 2060 890}%
\special{pa 2078 876}%
\special{pa 2100 888}%
\special{pa 2074 822}%
\special{fp}%
%
\special{pn 13}%
\special{pa 2616 1134}%
\special{pa 2590 822}%
\special{fp}%
\special{sh 1}%
\special{pa 2590 822}%
\special{pa 2576 890}%
\special{pa 2594 876}%
\special{pa 2616 888}%
\special{pa 2590 822}%
\special{fp}%
%
\special{pn 13}%
\special{pa 1604 1176}%
\special{pa 1576 864}%
\special{fp}%
\special{sh 1}%
\special{pa 1576 864}%
\special{pa 1562 932}%
\special{pa 1582 916}%
\special{pa 1602 928}%
\special{pa 1576 864}%
\special{fp}%
\put(24.6900,-24.6000){\makebox(0,0)[lt]{$\,$}}%
%
\special{pn 8}%
\special{pa 2682 566}%
\special{pa 2010 732}%
\special{dt 0.045}%
\special{sh 1}%
\special{pa 2010 732}%
\special{pa 2080 736}%
\special{pa 2062 718}%
\special{pa 2070 696}%
\special{pa 2010 732}%
\special{fp}%
%
\special{pn 8}%
\special{pa 2728 740}%
\special{pa 2350 780}%
\special{dt 0.045}%
\special{sh 1}%
\special{pa 2350 780}%
\special{pa 2418 794}%
\special{pa 2404 774}%
\special{pa 2414 754}%
\special{pa 2350 780}%
\special{fp}%
\put(27.1800,-6.1600){\makebox(0,0)[lb]{$H\cdot o$}}%
\put(27.6400,-7.9700){\makebox(0,0)[lb]{$H\cdot p$}}%
%
\special{pn 20}%
\special{sh 1}%
\special{ar 2138 1728 10 10 0  6.28318530717959E+0000}%
\special{sh 1}%
\special{ar 2138 1728 10 10 0  6.28318530717959E+0000}%
%
\special{pn 20}%
\special{sh 1}%
\special{ar 2488 1810 10 10 0  6.28318530717959E+0000}%
\special{sh 1}%
\special{ar 2488 1810 10 10 0  6.28318530717959E+0000}%
\put(24.6900,-18.5100){\makebox(0,0)[lt]{$p$}}%
\put(21.9300,-17.0300){\makebox(0,0)[lt]{$o$}}%
%
\special{pn 13}%
\special{pa 2488 1794}%
\special{pa 2488 1498}%
\special{fp}%
\special{sh 1}%
\special{pa 2488 1498}%
\special{pa 2468 1564}%
\special{pa 2488 1550}%
\special{pa 2508 1564}%
\special{pa 2488 1498}%
\special{fp}%
%
\special{pn 13}%
\special{pa 1972 1826}%
\special{pa 1972 1522}%
\special{fp}%
\special{sh 1}%
\special{pa 1972 1522}%
\special{pa 1952 1588}%
\special{pa 1972 1574}%
\special{pa 1992 1588}%
\special{pa 1972 1522}%
\special{fp}%
%
\special{pn 13}%
\special{pa 1964 1300}%
\special{pa 1926 938}%
\special{fp}%
\special{sh 1}%
\special{pa 1926 938}%
\special{pa 1914 1006}%
\special{pa 1932 990}%
\special{pa 1954 1002}%
\special{pa 1926 938}%
\special{fp}%
%
\special{pn 13}%
\special{pa 2432 1258}%
\special{pa 2396 920}%
\special{fp}%
\special{sh 1}%
\special{pa 2396 920}%
\special{pa 2384 988}%
\special{pa 2402 974}%
\special{pa 2424 984}%
\special{pa 2396 920}%
\special{fp}%
%
\special{pn 13}%
\special{pa 2940 2000}%
\special{pa 2948 1852}%
\special{fp}%
\special{sh 1}%
\special{pa 2948 1852}%
\special{pa 2924 1916}%
\special{pa 2946 1904}%
\special{pa 2964 1920}%
\special{pa 2948 1852}%
\special{fp}%
\put(30.0300,-18.6700){\makebox(0,0)[lt]{'s $:$ $(X^H_{\beta,i})^{\ast}|_M$}}%
%
\special{pn 20}%
\special{sh 1}%
\special{ar 2092 1134 10 10 0  6.28318530717959E+0000}%
\special{sh 1}%
\special{ar 2092 1134 10 10 0  6.28318530717959E+0000}%
%
\special{pn 8}%
\special{pa 3216 436}%
\special{pa 3390 436}%
\special{fp}%
\put(32.4200,-4.1000){\makebox(0,0)[lb]{$M$}}%
%
\special{pn 8}%
\special{pa 3206 436}%
\special{pa 3206 270}%
\special{fp}%
\end{picture}%
\hspace{0truecm}}

\vspace{0.5truecm}

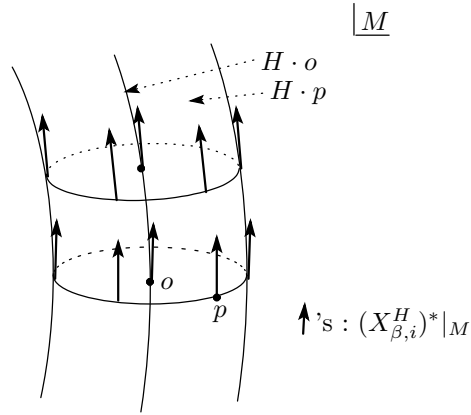
\centerline{{\bf Figure 6.}}

\vspace{0.5truecm}

\noindent
For $Z\in\mathfrak b^{\mathbb C}$, 
$(Y^V_{\beta,i})^{\ast}_{{\rm Exp}_o(Z)}$ and $(X^H_{\beta,i})^{\ast}_{{\rm Exp}_o(Z)}$ are described as 
$$(Y^V_{\beta,i})^{\ast}_{{\rm Exp}_o(Z)}=-\sin(\beta^{\mathbb C}(Z))(\exp\,Z)_{\ast}(X^V_{\beta,i})\leqno{(5.1)}$$
and 
$$(X^H_{\beta,i})^{\ast}_{{\rm Exp}_o(Z)}=\cos(\beta^{\mathbb C}(Z))(\exp\,Z)_{\ast}(X^H_{\beta,i}),\leqno{(5.2)}$$
respectively.  
A basis of $T_{\overline{\tau}({\bf s})}(H\cdot\overline{\tau}({\bf s}))
(=T_{\overline{\tau}({\bf s})}L_{\tau}\cap 
T^{\perp}_{\overline{\tau}({\bf s})}(L_{\tau})_0)$ is given by 
\begin{align*}
&\left(\mathop{\cup}_{\beta\in(\triangle_{\mathfrak b}^V)_+}
\{(Y^V_{\beta,i})^{\ast}_{\overline{\tau}({\bf s})}\,|\,i=1,\cdots,m^V_{\beta}\}\right)\\
&\cup
\left(\mathop{\cup}_{\beta\in(\triangle_{\mathfrak b}^H)_+}\{(X^H_{\beta,i})^{\ast}_{\overline{\tau}({\bf s})}\,|\,
i=1,\cdots,m^H_{\beta}\}\right).
\end{align*}
On the other hand, a basis of $T_{\overline{\tau}({\bf s})}(L_{\tau})_0(\subset 
T_{\overline{\tau}({\bf s})}(\Sigma^{\mathbb C})\,(\approx\mathbb C^r))$ is given by 
$$\left\{\frac{d\tau_1}{ds_1}\,{\bf e}_1,\cdots,\frac{d\tau_r}{ds_r}\,{\bf e}_r\right\},$$
where $\displaystyle{{\bf e}_i:=(0,\cdots,0,\mathop{1}^i,0,\cdots 0)}$ 
($\displaystyle{\mathop{1}^i}$ means that $i$-component is equal to $1$).  
Let $(\triangle_{\mathfrak b}^V)_+=$\newline
$\{\beta^V_i\,|\,i=1,\cdots,k_V\}$ and 
$(\triangle_{\mathfrak b}^H)_+=\{\beta^H_i\,|\,i=1,\cdots,k_H\}$.  
From $(5.1)$ and $(5.2)$, we have 
\begin{align*}
&(\Omega_0)_{\overline{\tau}({\bf s})}\left((Y^V_{\beta^V_1,1})^{\ast}_{\overline{\tau}({\bf s})},\cdots,
(Y^V_{\beta^V_1,m_{\beta^V_1}^V})^{\ast}_{\overline{\tau}({\bf s})},\cdots,\right.\\
&\hspace{1.7truecm}(Y^V_{\beta^V_{k_V},1})^{\ast}_{\overline{\tau}({\bf s})},
\cdots,(Y^V_{\beta^V_{k_V},m_{\beta^V_{k_V}}^V})^{\ast}_{\overline{\tau}({\bf s})},\\
&\hspace{1.7truecm}(X^H_{\beta^H_1,1})^{\ast}_{\overline{\tau}({\bf s})},\cdots,
(X^H_{\beta^H_1,m_{\beta^H_1}^H})^{\ast}_{\overline{\tau}({\bf s})},\cdots,\\
&\hspace{1.7truecm}(X^H_{\beta^H_{k_H},1})^{\ast}_{\overline{\tau}({\bf s})},
\cdots,(X^H_{\beta^H_{k_H},m_{\beta^H_{k_H}}^H})^{\ast}_{\overline{\tau}({\bf s})},\\
&\hspace{1.7truecm}
\left.\frac{d\tau_1}{ds_1}\,{\bf e}_1,\cdots,\frac{d\tau_r}{ds_r}\,{\bf e}_r\right)\\
=&\mathop{\Pi}_{\beta\in(\triangle_{\mathfrak b}^V)_+}\sin^{m_{\beta}^V}(-\beta^{\mathbb C}(\tau({\bf s})))
\cdot\mathop{\Pi}_{\beta\in(\triangle_{\mathfrak b}^H)_+}\cos^{m_{\beta}^H}(\beta^{\mathbb C}(\tau({\bf s})))
\cdot\mathop{\Pi}_{i=1}^r\frac{d\tau_i}{ds_i}\\
&\times
(\Omega_0)_o\left(X^V_{\beta^V_1,1},\cdots,X^V_{\beta^V_1,m_{\beta^V_1}^V},\cdots,
X^V_{\beta^V_{k_V},1},\cdots,X^V_{\beta^V_{k_V},m_{\beta^V_{k_V}}^V},\right.\\
&\hspace{1.7truecm}X^H_{\beta^H_1,1},\cdots,X^H_{\beta^H_1,m_{\beta^H_1}^H},
\cdots,X^H_{\beta^H_{k_H},1},\cdots,X^H_{\beta^H_{k_H},m_{\beta^H_{k_H}}^H},\\
&\hspace{1.7truecm}\left.(\exp(\tau({\bf s})))_{\ast}^{-1}({\bf e}_1),\cdots,
(\exp(\tau({\bf s})))_{\ast}^{-1}({\bf e}_r)\right)
\end{align*}
\newpage
\begin{align*}
=&-\mathop{\Pi}_{\beta\in(\triangle_{\mathfrak b}^V)_+}\sin^{m_{\beta}^V}
\left(\sum_{i=1}^r\tau_i(s_i)\beta({\bf e}_i)\right)
\cdot\mathop{\Pi}_{\beta\in(\triangle_{\mathfrak b}^H)_+}\cos^{m_{\beta}^H}
\left(\sum_{i=1}^r\tau_i(s_i)\beta({\bf e}_i)\right)\\
&\hspace{0.5truecm}\times\mathop{\Pi}_{i=1}^r\frac{d\tau_i}{ds_i}\cdot
(\Omega_0)_o\left(X^V_{\beta^V_1,1},\cdots,X^V_{\beta^V_1,m_{\beta^V_1}^V},\cdots,
X^V_{\beta^V_{k_V},1},\cdots,X^V_{\beta^V_{k_V},m_{\beta^V_{k_V}}^V},\right.\\
&\hspace{3.5truecm}X^H_{\beta^H_1,1},\cdots,X^H_{\beta^H_1,m_{\beta^H_1}^H},
\cdots,X^H_{\beta^H_{k_H},1},\cdots,\\
&\hspace{3.5truecm}\left.X^H_{\beta^H_{k_H},m_{\beta^H_{k_H}}^H},{\bf e}_1,\cdots,{\bf e}_r\right).
\end{align*}
In the last equality, we used the fact that 
$\exp(\tau({\bf s})))_{\ast}^{-1}({\bf e}_i)={\bf e}_i$ ($i=1,\cdots,r$) hold under the identification 
$T_{\overline{\tau}({\bf s})}\Sigma=T_o\Sigma=\mathbb C^r$ because $\exp(\tau({\bf s})))_{\ast}^{-1}$ is 
the parallel translation along the geodesic $t\mapsto{\rm Exp}_o(t\tau({\bf s}))$ in $\Sigma^{\mathbb C}$.  
It is clear that 
\begin{align*}
&(\Omega_0)_o\left(X^V_{\beta^V_1,1},\cdots,X^V_{\beta^V_1,m_{\beta^V_1}^V},\cdots,
X^V_{\beta^V_{k_V},1},\cdots,X^V_{\beta^V_{k_V},m_{\beta^V_{k_V}}^V},\right.\\
&\left.
X^H_{\beta^H_1,1},\cdots,X^H_{\beta^H_1,m_{\beta^H_1}^H},\cdots,X^H_{\beta^H_{k_H},1},\cdots,
X^H_{\beta^H_{k_H},m_{\beta^H_{k_H}}^H},{\bf e}_1,\cdots,{\bf e}_r\right)
\end{align*}
is a nonzero real constant independent of ${\bf s}=(s_1,\cdots,s_r)$.  
From these facts and Proposition 5.3, we obtain the following fact for $L_{\tau}$.  

\vspace{0.5truecm}

\noindent
{\bf Theorem 5.4.} {\sl The submanifold $L_{\tau}$ is a special Lagrangian submanifold of phase $\theta$ 
if and only if $\tau_1,\cdots,\tau_r$ satisfy the following ordinary differential equation:
{\small
$$\begin{array}{l}
\displaystyle{{\rm Im}\left(e^{\sqrt{-1}\theta}\cdot
\mathop{\Pi}_{\beta\in(\triangle_{\mathfrak b}^V)_+}\sin^{m_{\beta}^V}
\left(\sum_{i=1}^r\tau_i(s_i)\beta({\bf e}_i)\right)\right.}\\
\hspace{1truecm}\displaystyle{\left.\times\mathop{\Pi}_{\beta\in(\triangle_{\mathfrak b}^H)_+}\cos^{m_{\beta}^H}
\left(\sum_{i=1}^r\tau_i(s_i)\beta({\bf e}_i)\right)
\cdot\mathop{\Pi}_{i=1}^r\frac{d\tau_i}{ds_i}\right)=0.}
\end{array}\leqno{(5.3)}$$}
}

\vspace{0.5truecm}

Next we shall give solutions of the ordinary differential equation $(5.3)$.  
Let $\tau_i(s_i)=\varphi_i(s_i)+\sqrt{-1}\rho_i(s_i)$ ($i=1,\cdots,r$), 
where $\varphi_i$ and $\rho_i$ are real-valued functions.  
Set 
\begin{align*}
&F(\tau_1(s_1),\cdots,\tau_r(s_r))(=F(\varphi_1(s_1),\rho_1(s_1),\cdots,\varphi_r(s_r),\rho_r(s_r)))\\
:=&e^{\sqrt{-1}\theta}\cdot
\mathop{\Pi}_{\beta\in(\triangle_{\mathfrak b}^V)_+}\sin^{m_{\beta}^V}
\left(\sum_{i=1}^r\tau_i(s_i)\beta({\bf e}_i)\right)\\
&\times\mathop{\Pi}_{\beta\in(\triangle_{\mathfrak b}^H)_+}\cos^{m_{\beta}^H}
\left(\sum_{i=1}^r\tau_i(s_i)\beta({\bf e}_i)\right),
\end{align*}
Also, let 
$$F(\varphi_1,\rho_1,\cdots,\varphi_r,\rho_r)=u_0(\varphi_1,\rho_1,\cdots,\varphi_r,\rho_r)
+\sqrt{-1}v_0(\varphi_1,\rho_1,\cdots,\varphi_r,\rho_r)$$
and 
\begin{align*}
&\int F(\varphi_1,\rho_1,\cdots,\varphi_r,\rho_r)\,d\varphi_1\\
=&U_0(\varphi_1,\rho_1,\cdots,\varphi_r,\rho_r)+\sqrt{-1}V_0(\varphi_1,\rho_1,\cdots,\varphi_r,\rho_r),
\end{align*}
where $u_0,\,v_0,\,U_0$ and $V_0$ are real-valued functions.  
Define $u_1$ and $v_1$ by 
\begin{align*}
&u_1(\varphi_1(s_1),\rho_1(s_1),\varphi'_1(s_1),\rho'_1(s_1),\varphi_2(s_2),\rho_2(s_2),\cdots,\varphi_r(s_r),
\rho_r(s_r))\\
:=&\frac{\partial}{\partial s_1}
\left(U_0(\varphi_1(s_1),\rho_1(s_1),\cdots,\varphi_r(s_r),\rho_r(s_r))\right)
\end{align*}
and 
\begin{align*}
&v_1(\varphi_1(s_1),\rho_1(s_1),\varphi'_1(s_1),\rho'_1(s_1),\varphi_2(s_2),\rho_2(s_2),\cdots,\varphi_r(s_r),
\rho_r(s_r))\\
:=&\frac{\partial}{\partial s_1}
\left(V_0(\varphi_1(s_1),\rho_1(s_1),\cdots,\varphi_r(s_r),\rho_r(s_r))\right).
\end{align*}
It is clear that 
\begin{align*}
&u_1(\varphi_1(s_1),\rho_1(s_1),\varphi'_1(s_1),\rho'_1(s_1),\varphi_2(s_2),\rho_2(s_2),\cdots,\varphi_r(s_r),
\rho_r(s_r))\\
&+\sqrt{-1}v_1(\varphi_1(s_1),\rho_1(s_1),\varphi'_1(s_1),\rho'_1(s_1),\varphi_2(s_2),\rho_2(s_2),\cdots,
\varphi_r(s_r),\rho_r(s_r))\\
=&\frac{\partial}{\partial s_1}\left(\left(\int\,F(\varphi_1,\rho_1,\cdots,\varphi_r,\rho_r)d\varphi_1\right)
(\varphi_1(s_1),\rho_1(s_1),\cdots,\varphi_r(s_r),\rho_r(s_r)\right).
\end{align*}
Let 
\begin{align*}
&\int(u_1(\varphi_1,\varphi'_1,\rho_1,\rho'_1,\varphi_2,\rho_2,\cdots,\varphi_r,\rho_r)\\
&\hspace{0.6truecm}+\sqrt{-1}v_1(\varphi_1,\varphi'_1,\rho_1,\rho'_1,\varphi_2,\rho_2,\cdots,\varphi_r,\rho_r))\,
d\varphi_2\\
=&U_1(\varphi_1,\varphi'_1,\rho_1,\rho'_1,\varphi_2,\rho_2,\cdots,\varphi_r,\rho_r)\\
&+\sqrt{-1}V_1(\varphi_1,\varphi'_1,\rho_1,\rho'_1,\varphi_2,\rho_2,\cdots,\varphi_r,\rho_r),\\
\end{align*}
where $U_1$ and $V_1$ are real-valued functions.  
In the sequel, we define $u_i,v_i,U_i$ and $V_i$ ($i=2,\cdots,r$) by repeating the same process.  
Set 
$$\widehat F(\varphi_1,\rho_1,\cdots,\varphi_r,\rho_r):=\int\cdots\int F(\varphi_1,\rho_1,\cdots,\varphi_r,\rho_r)\,
d\varphi_1\cdots d\varphi_r.$$
It is easy to show that 
\begin{align*}
&(u_r+\sqrt{-1}v_r)(\varphi_1(s_1),\rho_1(s_1),\varphi'_1(s_1),\rho'_1(s_1),\cdots,\\
&\hspace{2.6truecm}\varphi_r(s_r),\rho_r(s_r),\varphi'_r(s_r),\rho'_r(s_r))\\
=&\frac{\partial^r}{\partial s_1\cdots\partial s_r}
\left(\widehat F(\varphi_1(s_1),\rho_1(s_1),\cdots,\varphi_r(s_r),\rho_r(s_r))\right)
\end{align*}

\vspace{0.5truecm}

\noindent
{\bf Corollary 5.5.} {\sl Let $F$ be the complex-valued function over $\mathbb R^{2r}$ defined by 
\begin{align*}
F(\varphi_1,\rho_1,\cdots,\varphi_r,\rho_r)
:=&e^{\sqrt{-1}\theta}\cdot
\mathop{\Pi}_{\beta\in(\triangle_{\mathfrak b}^V)_+}\sin^{m_{\beta}^V}
\left(\sum_{i=1}^r(\varphi_i+\sqrt{-1}\rho_i)\cdot\beta({\bf e}_i)\right)\\
&\cdot\mathop{\Pi}_{\beta\in(\triangle_{\mathfrak b}^H)_+}\cos^{m_{\beta}^H}
\left(\sum_{i=1}^r(\varphi_i+\sqrt{-1}\rho_i)\cdot\beta({\bf e}_i)\right).
\end{align*}
If $\tau_i(s_i)=\varphi_i(s_i)+\sqrt{-1}\rho_i(s_i)$ ($i=1,\cdots,r$) satisfy 
$${\rm Im}\left(\widehat F(\varphi_1(s_1),\rho_1(s_1),\cdots,\varphi_r(s_r),\rho_r(s_r))\right)=0,$$
then they are a solution of $(5.3)$ and hence $L_{\tau}$ ($\tau:=\tau_1\times\cdots\times\tau_r$) is 
a special Lagrangian submanifold of phase $\theta$.  
}

\vspace{0.5truecm}

\noindent
{\it Proof.} Since $F$ is a holomorphic function, we have 
$\displaystyle{\frac{\partial u_0}{\partial\varphi_1}=\frac{\partial v_0}{\partial\rho_1}}$ and 
$\displaystyle{\frac{\partial u_0}{\partial\rho_1}=}$\newline
$\displaystyle{-\frac{\partial v_0}{\partial\varphi_1}}$.  
From these relations and the definitions of $U_0$ and $V_0$, we have 
$$\frac{\partial U_0}{\partial\varphi_1}=\frac{\partial V_0}{\partial\rho_1}=u_0,\,\,
\frac{\partial U_0}{\partial\rho_1}=-v_0\,\,\,{\rm and}\,\,\,\frac{\partial V_0}{\partial\varphi_1}=v_0.$$
Hence we obtain 
\begin{align*}
&F(\tau_1(s_1),\cdots,\tau_r(s_r))\cdot\tau'_1(s_1)\\
=&\frac{\partial}{\partial s_1}\left(\left(\int\,F(\varphi_1,\rho_1,\cdots,\varphi_r,\rho_r)\,d\varphi_1\right)
(\varphi_1(s_1),\rho_1(s_1),\cdots,\varphi_r(s_r),\rho_r(s_r))\right).
\end{align*}
Since $u_1+\sqrt{-1}v_1$ is holomorphic with respect to $\tau_2(=\varphi_2+\sqrt{-1}\rho_2)$, 
we have 
$\displaystyle{\frac{\partial u_1}{\partial\varphi_2}=\frac{\partial v_1}{\partial\rho_2}}$ and 
$\displaystyle{\frac{\partial u_1}{\partial\rho_2}=-\frac{\partial v_1}{\partial\varphi_2}}$.  
From these relations and the definitions of $U_1$ and $V_1$, we have 
$$\frac{\partial U_1}{\partial\varphi_2}=\frac{\partial V_1}{\partial\rho_2}=u_1,\,\,
\frac{\partial U_1}{\partial\rho_2}=-v_1\,\,\,{\rm and}\,\,\,\frac{\partial V_1}{\partial\varphi_2}=v_1.$$
Hence we obtain 
{\small\begin{align*}
&F(\tau_1(s_1),\cdots,\tau_r(s_r))\cdot\tau'_1(s_1)\cdot\tau'_2(s_2)=\\
&\frac{\partial^2}{\partial s_1\partial s_2}\hspace{-0.15truecm}\left(\left(\hspace{-0.1truecm}\int
\hspace{-0.2truecm}\int F(\varphi_1,\rho_1,\cdots,\varphi_r,\rho_r)d\varphi_1d\varphi_2\right)
(\varphi_1(s_1),\rho_1(s_1),\cdots,\varphi_r(s_r),\rho_r(s_r))\hspace{-0.1truecm}\right).
\end{align*}
}
In the sequel, by repeating the same discussion, we obtain 
\begin{align*}
&F(\tau_1(s_1),\cdots,\tau_r(s_r))\cdot\tau'_1(s_1)\cdots\tau'_r(s_r)\\
=&\frac{\partial^r}{\partial s_1\cdots\partial s_r}
\left(\widehat F(\varphi_1(s_1),\rho_1(s_1),\cdots,\varphi_r(s_r),\rho_r(s_r))\right).
\end{align*}
From this relation, we can derive the statement of this corollary directly.  \qed

\vspace{0.25truecm}

We consider the case where $N=G/K$ is an $dm$-dimensional simply connected rank one symmetirc space of compact type 
and constant maximal sectional curvature $4c$, that is, 
$N=\mathbb FP^m(4c)$ ($\mathbb F=\mathbb C,\,\mathbb Q$ or $\mathbb O$), where $\mathbb Q$ (resp. $\mathbb O$) 
denotes the quaternionic algebra (resp. the Octonian) and $d$ is given by 
$d=2$ (when $\mathbb F=\mathbb C$), $d=4$ (when $\mathbb F=\mathbb Q$) or $d=8$ (when $\mathbb F=\mathbb O$).  
Note that 
$$\mathbb FP^m(4c)=\left\{
\begin{array}{ll}
SU(m+1)/S(U(1)\times U(m)) & (\mathbb F=\mathbb C)\\
Sp(m+1)/(Sp(1)\times Sp(m)) & (\mathbb F=\mathbb Q)\\
F_4/Spin(9) & (\mathbb F=\mathbb O,\,\,m=2).
\end{array}\right.$$
In these cases, we have $\triangle_+=\{\sqrt ce^{\ast},\,2\sqrt ce^{\ast}\}$, where 
$e^{\ast}$ denotes the dual $1$-form of the unit vector $e$ of $\mathfrak b$ (${\rm dim}\,\mathfrak b=1$ 
in these cases).  Also, we have 
${\rm dim}\,\mathfrak p_{\sqrt c{\bf e}^{\ast}}=d(m-1)$ and ${\rm dim}\,\mathfrak p_{2\sqrt c{\bf e}^{\ast}}=d-1$.  
Hence, as a corollary of Theorem 5.4, we obtain the following fact.  

\vspace{0.25truecm}

\noindent
{\bf Corollary 5.6.} {\sl Let $H\curvearrowright\mathbb FP^m(4c)$ be a Hermann action.  
Then the submanifold $L_{\tau}$ is a special Lagrangian submanifold of phase $\theta$ 
if and only if $\tau$ satisfies the following ordinary differential equation:
{\small
$$\begin{array}{l}
\displaystyle{{\rm Im}\left(e^{\sqrt{-1}\theta}\cdot
\sin^{m^V}(\sqrt c\tau(s))\cdot\sin^{d-1}(2\sqrt c\tau(s))
\cdot\cos^{m^H}(\sqrt c\tau(s))\cdot\frac{d\tau}{ds}\right)=0,}
\end{array}\leqno{(5.4)}$$}
where $m^V$ (resp. $m^H$) denotes $m_{\sqrt c {\bf e}^{\ast}}^V$ (resp. $m_{\sqrt c{\bf e}^{\ast}}^H$).}

\vspace{0.25truecm}

\noindent
{\it Proof.} Let $\{J_1,\cdots,J_{d-1}\}$ be the complex structure, the canonical local basis of 
the quaternionic structure or the Cayley structure of $\mathbb FP^m(4c)$.  
Then, since $\mathbb FP^{dm}(4c)$ is of rank one, this action is of cohomogeneity one.  
It is shown that this action is commutative (i.e., $\theta\circ\sigma=\sigma\circ\theta$).  
In fact, Hermann actions on $\mathbb FP^m(4c)$ are classified as in Table 1 and all of Hermann actions in Table 1 
are commutative.  
Since $H\curvearrowright\mathbb FP^{dm}(4c)$ is commutative, it is shown that $H\cdot o$ and the normal umbrella 
${\rm Exp}_o(T^{\perp}_o(H\cdot o))$ are reflective submanifolds, that is, they are $J_i$-invariant 
($i=1,\cdots,d-1$).  This implies that $T^{\perp}_o(H\cdot o)=\mathfrak p\cap\mathfrak q$ includes 
$\mathfrak p_{2\sqrt ce^{\ast}}$.  Hence we have $m_{2\sqrt ce^{\ast}}^V=d-1$ and $m_{2\sqrt ce^{\ast}}^H=0$.  
Therefore, the statement of this corollary follows from Theorem 5.4 directly.  \qed

\vspace{0.25truecm}

K. Hashimoto and K. Mashimo (\cite{HM}) gave the ordinary differential equation corresponding to $(5.3)$ 
for the Hamiltonian action $K\curvearrowright TS^n(1)(=SO(n+1,{\Bbb C})/SO(n,{\Bbb C}))$ induced from 
the restricted action $K\curvearrowright S^n(1)(=SO(n+1)/SO(n))$ of the linear isotropy action of any irreducible 
rank two symmetric space $G/K$ (see Theorem 5.2 of \cite{HM}), where $n:={\rm dim}\,G/K-1$ and $S^n(1)$ is the unit 
sphere of $T_o(G/K)$.  
Here we note that the action $K\curvearrowright S^n(1)$ is a non-Hermann action of cohomogeneity one 
(stated in (ii) of Theorem A in \cite{Kol}).  

Let $K\curvearrowright\mathbb CP^n(4)(=SU(n+1)/S(U(1)\times U(n)))$ be the action induced from the restricted action 
$K\curvearrowright S^{2n+1}$ of the linear isotropy action of any irreducible rank two Hermitian symmetric space 
$G/K$ through the Hopf fibration $\pi:S^{2n+1}(1)\to\mathbb CP^n(4)$, where $n=\frac{1}{2}\cdot{\rm dim}\,G/K-1$ and 
$S^{2n+1}(1)$ is the unit sphere of $T_o(G/K)$.  Here we note that the action $K\curvearrowright\mathbb CP^n(4)$ is 
a non-Hermann action of cohomogeneity one (stated in (iii) of Theorem A in \cite{Kol}).  
Recently, M. Arai and K. Baba (\cite{AB}) gave the ordinary differential equation corresponding to $(5.3)$ 
for the Hamiltonian action $K\curvearrowright T\mathbb CP^n(4)
(=SL(n+1,{\Bbb C})/(SL(1,{\Bbb C})\times SL(n,{\Bbb C})))$ induced from the action 
$K\curvearrowright\mathbb CP^n(4)$ (see Theorems 2.1-2.4 of \cite{AB}).  

According to the classification of cohomogeneity one actions on irreducible symmetric spaces $G/K$ of compact type 
such that $G$ is simple (i.e., $G/K$ is of type I in \cite{H}) by A. Kollross (see Theorem B of \cite{Kol}), 
any Hermann action on $\mathbb FP^m(4c)$ is orbit equivalent to one of Hermann actions in Table 1.  

\vspace{0truecm}

{\small
$$\begin{tabular}{|c|c|}
\hline
$H\curvearrowright G/K(=\mathbb FP^m(4c))$ & \\
\hline
$S(U(1)\times U(m))\curvearrowright SU(m+1)/S(U(1)\times U(m))$ & isotropy action\\
\hline
$S(U(m)\times U(1))\curvearrowright SU(m+1)/S(U(1)\times U(m))$ & $-$\\
\hline
$SO(m+1)\curvearrowright SU(m+1)/S(U(1)\times U(m))$ & $-$\\
\hline
$Sp(1)\times Sp(m)\curvearrowright Sp(m+1)/(Sp(1)\times Sp(m))$ & isotropy action\\
\hline
$Sp(m)\times Sp(1)\curvearrowright Sp(m+1)/(Sp(1)\times Sp(m))$ & $-$\\
\hline
$U(m+1)\curvearrowright Sp(m+1)/(Sp(1)\times Sp(m))$ & $-$\\
\hline
$Spin(9)\curvearrowright F_4/Spin(9)$ & isotropy action\\
\hline
$Sp(3)\cdot Sp(1)\curvearrowright F_4/Spin(9)$ & $-$\\
\hline
\end{tabular}$$
}

\vspace{0.1truecm}

\centerline{\bf Table 1 : Hermann actions on $\mathbb FP^m(4c)$}

\vspace{0.2truecm}


For all Hermann actions of cohomogeneity two on irreducible rank two symmetric spaces 
of compact type, we shall give the following datas:
$$(\triangle_{\mathfrak b})_+,\,\,(\triangle_{\mathfrak b})_+^V,\,\,(\triangle_{\mathfrak b})^H_+,\,\,
m_{\beta}^V\,\,(\beta\in(\triangle_{\mathfrak b})_+^V),\,\,\,\,m_{\beta}^H\,\,
(\beta\in\triangle_{\mathfrak b})_+^H).$$
All of such Hermann actions and the above datas for the actions are given as in Table 2.  
By using Table 2, we can explicitly describe the ordinary differential equation $(5.3)$ 
for the Hermann actions of cohomogeneity two on irreducible rank two symmetric spaces of compact type.  
In Table 2, in the case where 
$(\triangle_{\mathfrak b})_+=(\triangle_{\mathfrak b})_+^V\cup(\triangle_{\mathfrak b})_+^H$ is 
$\{\beta_1,\beta_2,\beta_1+\beta_2\}$, it implies a positive 
root system of the root system of (${\mathfrak a}_2$)-type ($(\beta_1({\bf e}_1),\beta_1({\bf e}_2))=(2,0),
(\beta_2({\bf e}_1),\beta_2({\bf e}_2))=(-1,\sqrt 3)$), 
in the case where $(\triangle_{\mathfrak b})_+$ is $\{\beta_1,\beta_2,\beta_1+\beta_2, 2\beta_1+\beta_2\}$, it 
implies a positive root system of the root system of (${\mathfrak b}_2$)(=(${\mathfrak c}_2$))-type 
($(\beta_1({\bf e}_1),\beta_1({\bf e}_2))=(1,0),(\beta_2({\bf e}_1),\beta_2({\bf e}_2))=(-1,1)$) and 
in the case where $(\triangle_{\mathfrak b})_+$ is 
$\{\beta_1,\beta_2,\beta_1+\beta_2,\beta_1+2\beta_2,\beta_1+3\beta_2,2\beta_1+3\beta_2\}$, it implies 
a positive root system of the root system of (${\mathfrak g}_2$)-type 
($(\beta_1({\bf e}_1),\beta_1({\bf e}_2))=(2\sqrt3,0),(\beta_2({\bf e}_1),\beta_2({\bf e}_2))=(-\sqrt3,1)$).  
Also, $\rho_i$ ($i=1,\cdots,16$) imply automorphisms of $G$ whose dual actions are given as in Table 3 and 
$(\bullet)^2$ implies the product Lie group $(\bullet)\times(\bullet)$ of a Lie group $(\bullet)$, 
$\displaystyle{\mathop{\bullet}_{(m)}}$ in the column of $(\triangle_{\mathfrak b})_+$, 
$(\triangle_{\mathfrak b})^V_+$ and $(\triangle_{\mathfrak b})^H_+$ imply 
$m_{\bullet}=m,\,m_{\bullet}^V=m$ and $m_{\bullet}^H=m$, respectively.  
Note that Tables 2 and 3 are based on Tables 1 and 2 in \cite{Koi2}.  


\vspace{0.15truecm}

{\scriptsize
$$\begin{tabular}{|c|c|c|}
\hline
{\scriptsize$H\curvearrowright G/K$} & {\scriptsize$(\triangle_{\mathfrak b})^V_+,\,\,\,\,m_{\bullet}^V$} & 
{\scriptsize$(\triangle_{\mathfrak b})^H_+,\,\,\,\,m_{\bullet}^H$}\\
\hline
{\scriptsize$\rho_1(SO(3))\curvearrowright SU(3)/SO(3)$} & 
{\scriptsize$\displaystyle{\{\mathop{\beta_1}_{(1)}\}}$} & 
{\scriptsize$\displaystyle{\{\mathop{\beta_2}_{(1)},
\mathop{\beta_1+\beta_2}_{(1)}\}}$}\\
\hline
{\scriptsize$SO(6)\curvearrowright SU(6)/Sp(3)$} & 
{\scriptsize$\displaystyle{\{\mathop{\beta_1}_{(2)},\mathop{\beta_2}_{(2)},
\mathop{\beta_1+\beta_2}_{(2)}\}}$} & 
{\scriptsize$\displaystyle{\{\mathop{\beta_1}_{(2)},\mathop{\beta_2}_{(2)},
\mathop{\beta_1+\beta_2}_{(2)}\}}$}\\
\hline
{\scriptsize$\rho_2(Sp(3))\curvearrowright SU(6)/Sp(3)$} & 
{\scriptsize$\displaystyle{\{\mathop{\beta_1}_{(4)}\}}$} & 
{\scriptsize$\displaystyle{\{\mathop{\beta_2}_{(4)},
\mathop{\beta_1+\beta_2}_{(4)}\}}$}\\
\hline
{\scriptsize$SO(q+2)\curvearrowright$} & 
{\scriptsize$\displaystyle{\{\mathop{\beta_1}_{(q-2)},\mathop{\beta_2}_{(1)},
\mathop{\beta_1+\beta_2}_{(q-2)},}$} & 
{\scriptsize$\displaystyle{\{\mathop{\beta_1}_{(q-2)},\mathop{\beta_2}_{(1)},
\mathop{\beta_1+\beta_2}_{(q-2)},}$}\\
{\scriptsize$SU(q+2)/S(U(2)\times U(q))$} & 
{\scriptsize$\displaystyle{\mathop{2\beta_1+\beta_2}_{(1)}\}}$} & 
{\scriptsize$\displaystyle{\mathop{2\beta_1+\beta_2}_{(1)},\mathop{2\beta_1}_{(1)},
\mathop{2\beta_1+2\beta_2}_{(1)}\}}$}\\
{\scriptsize$(q\,>\,2)$}&&\\
\hline
{\scriptsize$S(U(j+1)\times U(q-j+1))\curvearrowright$} & 
{\scriptsize$\displaystyle{\{\mathop{\beta_1}_{(2j-2)},
\mathop{\beta_1+\beta_2}_{(2q-2j-2)},}$} & 
{\scriptsize$\displaystyle{\{\mathop{\beta_1}_{(2q-2j-6)},
\mathop{\beta_2}_{(2)},}$}\\
{\scriptsize$SU(q+2)/S(U(2)\times U(q))$} & 
{\scriptsize$\displaystyle{\mathop{2\beta_1}_{(1)},
\mathop{2\beta_1+2\beta_2}_{(1)}\}}$} & 
{\scriptsize$\displaystyle{\mathop{\beta_1+\beta_2}_{(2j-2)},
\mathop{2\beta_1+\beta_2}_{(2)}\}}$}\\
{\scriptsize$(q>2)$} & & \\
\hline
{\scriptsize$S(U(2)\times U(2))\curvearrowright$} & 
{\scriptsize$\displaystyle{\{\mathop{\beta_1}_{(1)},
\mathop{\beta_1+\beta_2}_{(1)}\}}$} & 
{\scriptsize$\displaystyle{\{\mathop{\beta_1}_{(1)},
\mathop{\beta_2}_{(1)},}$}\\
{\scriptsize$SU(4)/S(U(2)\times U(2))$} & 
{\scriptsize} & 
{\scriptsize$\displaystyle{\mathop{\beta_1+\beta_2}_{(1)},
\mathop{2\beta_1+\beta_2}_{(1)}\}}$}\\
{\scriptsize (non-isotropy gr. act.)} & &\\
\hline
{\scriptsize$SO(j+1)\times SO(q-j+1)\curvearrowright$} & 
{\scriptsize$\displaystyle{\{\mathop{\beta_1}_{(j-1)},
\mathop{\beta_1+\beta_2}_{(q-j-1)}\}}$} & 
{\scriptsize$\displaystyle{\{\mathop{\beta_1}_{(q-j-1)},
\mathop{\beta_2}_{(1)},}$}\\
{\scriptsize$SO(q+2)/SO(2)\times SO(q)$} & 
{\scriptsize} & 
{\scriptsize$\displaystyle{\mathop{\beta_1+\beta_2}_{(j-1)},
\mathop{2\beta_1+\beta_2}_{(1)}\}}$}\\
\hline
{\scriptsize$SO(4)\times SO(4)\curvearrowright$} & 
{\scriptsize$\displaystyle{\{\mathop{\beta_1}_{(2)},\mathop{\beta_2}_{(1)},
\mathop{\beta_1+\beta_2}_{(2)},}$} & 
{\scriptsize$\displaystyle{\{\mathop{\beta_1}_{(2)},
\mathop{\beta_1+\beta_2}_{(2)}\}}$}\\
{\scriptsize$SO(8)/U(4)$} & 
{\scriptsize$\displaystyle{\mathop{2\beta_1+\beta_2}_{(1)}\}}$} & 
{\scriptsize}\\
\hline
{\scriptsize$\rho_3(SO(4)\times SO(4))\curvearrowright$} & 
{\scriptsize$\displaystyle{\{\mathop{\beta_1}_{(2)},\mathop{\beta_1+\beta_2}_{(2)}
\}}$} & 
{\scriptsize$\displaystyle{\{\mathop{\beta_1}_{(2)},\mathop{\beta_2}_{(1)},
\mathop{\beta_1+\beta_2}_{(2)},}$}\\
{\scriptsize$SO(8)/U(4)$} & 
{\scriptsize} & 
{\scriptsize$\displaystyle{\mathop{2\beta_1+\beta_2}_{(1)}}$}\\
\hline
{\scriptsize$\rho_4(U(4))\curvearrowright SO(8)/U(4)$} & 
{\scriptsize$\displaystyle{\{\mathop{\beta_1}_{(1)},\mathop{\beta_1+\beta_2}_{(1)}
\}}$} & 
{\scriptsize$\displaystyle{\{\mathop{\beta_1}_{(3)},\mathop{\beta_2}_{(1)},
\mathop{\beta_1+\beta_2}_{(3)},}$}\\
{\scriptsize} & 
{\scriptsize} & 
{\scriptsize$\displaystyle{\mathop{2\beta_1+\beta_2}_{(1)}}$}\\
\hline
{\scriptsize$SO(4)\times SO(6)\curvearrowright$} & 
{\scriptsize$\displaystyle{\{\mathop{\beta_1}_{(2)},\mathop{\beta_2}_{(2)},
\mathop{\beta_1+\beta_2}_{(2)},}$} & 
{\scriptsize$\displaystyle{\{\mathop{\beta_1}_{(2)},\mathop{\beta_2}_{(2)},
\mathop{\beta_1+\beta_2}_{(2)},}$}\\
{\scriptsize$SO(10)/U(5)$} & 
{\scriptsize$\displaystyle{\mathop{2\beta_1+\beta_2}_{(2)},
\mathop{2\beta_1}_{(1)},\mathop{2\beta_1+2\beta_2}_{(1)}\}}$} & 
{\scriptsize$\displaystyle{\mathop{2\beta_1+\beta_2}_{(2)}}$}\\
\hline
{\scriptsize$SO(5)\times SO(5)\curvearrowright$} & 
{\scriptsize$\displaystyle{\{\mathop{\beta_1}_{(2)},\mathop{\beta_2}_{(2)},
\mathop{\beta_1+\beta_2}_{(2)},}$} & 
{\scriptsize$\displaystyle{\{\mathop{\beta_1}_{(2)},\mathop{\beta_2}_{(2)},
\mathop{\beta_1+\beta_2}_{(2)},}$}\\
{\scriptsize$SO(10)/U(5)$} & 
{\scriptsize$\displaystyle{\mathop{2\beta_1+\beta_2}_{(2)}\}}$} & 
{\scriptsize$\displaystyle{\mathop{2\beta_1+\beta_2}_{(2)},
\mathop{2\beta_1}_{(1)},\mathop{2\beta_1+2\beta_2}_{(1)}}$}\\
\hline
\end{tabular}$$
}

\vspace{0.15truecm}

\centerline{{\bf Table 2$\,$:$\,$Hermann actions on rank two symmetric spaces}}

\newpage


{\small 
$$\begin{tabular}{|c|c|c|}
\hline
{\scriptsize$H\curvearrowright G/K$}
& {\scriptsize$(\triangle_{\mathfrak b})^V_+,\,\,\,\,m_{\bullet}^V$} & 
{\scriptsize$(\triangle_{\mathfrak b})^H_+,\,\,\,\,m_{\bullet}^H$}\\
\hline
{\scriptsize$\rho_5(U(5))\curvearrowright SO(10)/U(5)$} & 
{\scriptsize$\displaystyle{\{\mathop{\beta_1}_{(4)},\mathop{2\beta_1}_{(1)},
\mathop{2\beta_1+2\beta_2}_{(1)}\}}$} & 
{\scriptsize$\displaystyle{\{\mathop{\beta_2}_{(4)},\mathop{\beta_1+\beta_2}_{(4)},
\mathop{2\beta_1+\beta_2}_{(4)}\}}$}\\
{\scriptsize} & 
{\scriptsize} & {\scriptsize}\\
\hline
{\scriptsize$SO(2)^2\times SO(3)^2\curvearrowright$} & 
{\scriptsize$\displaystyle{\{\mathop{\beta_1}_{(1)},
\mathop{\beta_2}_{(1)},\mathop{\beta_1+\beta_2}_{(1)},}$} & 
{\scriptsize$\displaystyle{\{\mathop{\beta_1}_{(1)},\mathop{\beta_2}_{(1)},
\mathop{\beta_1+\beta_2}_{(1)},}$}\\
{\scriptsize$(SO(5)\times SO(5))/SO(5)$} & 
{\scriptsize$\displaystyle{\mathop{2\beta_1+\beta_2}_{(1)}\}}$} & 
{\scriptsize$\displaystyle{\mathop{2\beta_1+\beta_2}_{(1)}\}}$}\\
\hline
{\scriptsize$\rho_6(SO(5))\curvearrowright$} & 
{\scriptsize$\displaystyle{\{\mathop{\beta_1}_{(2)}\}}$} & 
{\scriptsize$\displaystyle{\{\mathop{\beta_2}_{(2)},
\mathop{\beta_1+\beta_2}_{(2)},\mathop{2\beta_1+\beta_2}_{(2)}\}}$}\\
{\scriptsize$(SO(5)\times SO(5))/SO(5)$} & 
{\scriptsize} & {\scriptsize}\\
\hline
{\scriptsize$\rho_7(U(2))\curvearrowright Sp(2)/U(2)$} & 
{\scriptsize$\displaystyle{\{\mathop{\beta_1}_{(1)},
\mathop{\beta_1+\beta_2}_{(1)}\}}$} & 
{\scriptsize$\displaystyle{\{\mathop{\beta_2}_{(1)},
\mathop{2\beta_1+\beta_2}_{(1)}\}}$}\\
{\scriptsize} & 
{\scriptsize} & {\scriptsize}\\
\hline
{\scriptsize$SU(q+2)\curvearrowright$} & 
{\scriptsize$\displaystyle{\{\mathop{\beta_1}_{(2q-4)},\mathop{\beta_2}_{(2)},
\mathop{\beta_1+\beta_2}_{(2q-4)}\}}$} & 
{\scriptsize$\displaystyle{\{\mathop{\beta_1}_{(2q-4)},\mathop{\beta_2}_{(2)},
\mathop{\beta_1+\beta_2}_{(2q-4)}\}}$}\\
{\scriptsize$Sp(q+2)/Sp(2)\times Sp(q)$} & 
{\scriptsize$\displaystyle{\mathop{2\beta_1+\beta_2}_{(2)},
\mathop{2\beta_1}_{(1)},\mathop{2\beta_1+2\beta_2}_{(1)}\}}$} & 
{\scriptsize$\displaystyle{\mathop{2\beta_1+\beta_2}_{(2)},
\mathop{2\beta_1}_{(2)},\mathop{2\beta_1+2\beta_2}_{(2)}\}}$}\\
{\scriptsize$(q>2)$}&&\\
\hline
{\scriptsize$SU(4)\curvearrowright$} & 
{\scriptsize$\displaystyle{\{\mathop{\beta_1}_{(2)},\mathop{\beta_2}_{(1)},
\mathop{\beta_1+\beta_2}_{(2)}\}}$} & 
{\scriptsize$\displaystyle{\{\mathop{\beta_1}_{(2)},\mathop{\beta_2}_{(2)},
\mathop{\beta_1+\beta_2}_{(1)}\}}$}\\
{\scriptsize$Sp(4)/Sp(2)\times Sp(2)$} & 
{\scriptsize$\displaystyle{\mathop{2\beta_1+\beta_2}_{(1)}\}}$} & 
{\scriptsize$\displaystyle{\mathop{2\beta_1+\beta_2}_{(3)}\}}$}\\
\hline
{\scriptsize$U(4)\curvearrowright$} & 
{\scriptsize$\displaystyle{\{\mathop{\beta_1}_{(2)},\mathop{\beta_2}_{(2)},
\mathop{\beta_1+\beta_2}_{(2)}\}}$} & 
{\scriptsize$\displaystyle{\{\mathop{\beta_1}_{(2)},\mathop{\beta_2}_{(1)},
\mathop{\beta_1+\beta_2}_{(1)}\}}$}\\
{\scriptsize$Sp(4)/Sp(2)\times Sp(2)$} & 
{\scriptsize$\displaystyle{\mathop{2\beta_1+\beta_2}_{(2)}\}}$} & 
{\scriptsize$\displaystyle{\mathop{2\beta_1+\beta_2}_{(2)}\}}$}\\
\hline
{\scriptsize$Sp(j+1)\times Sp(q-j+1)\curvearrowright$} & 
{\scriptsize$\displaystyle{\{\mathop{\beta_1}_{(2j-4)},\mathop{2\beta_1}_{(3)},
\mathop{\beta_1+\beta_2}_{(4q-4j-4)}\}}$} & 
{\scriptsize$\displaystyle{\{\mathop{\beta_1}_{(4q-4j-4)},\mathop{\beta_2}_{(4)},
\mathop{\beta_1+\beta_2}_{(4j-4)}\}}$}\\
{\scriptsize$Sp(q+2)/Sp(2)\times Sp(q)$} & 
{\scriptsize$\displaystyle{\mathop{2\beta_1+2\beta_2}_{(3)}\}}$} & 
{\scriptsize$\displaystyle{\mathop{2\beta_1+\beta_2}_{(4)}\}}$}\\
{\scriptsize$(q>2)$}&&\\
\hline
{\scriptsize$Sp(2)\times Sp(2)\curvearrowright$} & 
{\scriptsize$\displaystyle{\{\mathop{\beta_1}_{(3)},
\mathop{\beta_1+\beta_2}_{(3)}\}}$} & 
{\scriptsize$\displaystyle{\{\mathop{\beta_1}_{(1)},\mathop{\beta_2}_{(3)},
\mathop{2\beta_1+\beta_2}_{(4)}\}}$}\\
{\scriptsize$Sp(4)/Sp(2)\times Sp(2)$} & 
{\scriptsize} & {\scriptsize}\\
\hline
{\scriptsize$SU(2)^2\cdot SO(2)^2\curvearrowright$} & 
{\scriptsize$\displaystyle{\{\mathop{\beta_1}_{(1)},\mathop{\beta_2}_{(1)},
\mathop{\beta_1+\beta_2}_{(1)},}$} & 
{\scriptsize$\displaystyle{\{\mathop{\beta_1}_{(1)},\mathop{\beta_2}_{(1)},
\mathop{\beta_1+\beta_2}_{(1)},}$}\\
{\scriptsize$(Sp(2)\times Sp(2))/Sp(2)$} & 
{\scriptsize$\displaystyle{\mathop{2\beta_1+\beta_2}_{(1)}\}}$} & 
{\scriptsize$\displaystyle{\mathop{2\beta_1+\beta_2}_{(1)}\}}$}\\
\hline
{\scriptsize$\rho_8(Sp(2))\curvearrowright$} & 
{\scriptsize$\displaystyle{\{\mathop{\beta_1}_{(2)},
\mathop{\beta_1+\beta_2}_{(2)}\}}$} & 
{\scriptsize$\displaystyle{\{\mathop{\beta_2}_{(2)},
\mathop{2\beta_1+\beta_2}_{(2)}\}}$}\\
{\scriptsize$(Sp(2)\times Sp(2))/Sp(2)$} & 
{\scriptsize} & {\scriptsize}\\
\hline
\end{tabular}$$
}

\vspace{0.3truecm}

\centerline{{\bf Table 2$\,$:$\,$Hermann actions on rank two symmetric spaces (continued)}}

\newpage

{\small
$$\begin{tabular}{|c|c|c|}
\hline
{\scriptsize$H\curvearrowright G/K$} 
& {\scriptsize$(\triangle_{\mathfrak b})^V_+,\,\,\,\,m_{\bullet}^V$} & 
{\scriptsize$(\triangle_{\mathfrak b})^H_+,\,\,\,\,m_{\bullet}^H$}\\
\hline
{\scriptsize$\rho_9(Sp(2))\curvearrowright$} & 
{\scriptsize$\displaystyle{\{\mathop{\beta_1}_{(2)},
\mathop{\beta_1+\beta_2}_{(2)}\}}$} & 
{\scriptsize$\displaystyle{\{\mathop{\beta_2}_{(2)},
\mathop{2\beta_1+\beta_2}_{(2)}\}}$}\\
{\scriptsize$(Sp(2)\times Sp(2))/Sp(2)$} & 
{\scriptsize} & {\scriptsize}\\
\hline
{\scriptsize$Sp(4)\curvearrowright$} & 
{\scriptsize$\displaystyle{\{\mathop{\beta_1}_{(4)},\mathop{\beta_2}_{(3)},
\mathop{\beta_1+\beta_2}_{(3)},}$} & 
{\scriptsize$\displaystyle{\{\mathop{\beta_1}_{(4)},\mathop{\beta_2}_{(3)},
\mathop{\beta_1+\beta_2}_{(6)},}$}\\
{\scriptsize$\displaystyle{E_6/Spin(10)\cdot U(1)}$} & 
{\scriptsize$\displaystyle{\mathop{2\beta_1+\beta_2}_{(4)}\}}$} & 
{\scriptsize$\displaystyle{\mathop{2\beta_1+\beta_2}_{(1)}\}}$}\\
\hline
{\scriptsize$SU(6)\cdot SU(2)\curvearrowright$} & 
{\scriptsize$\displaystyle{\{\mathop{\beta_1}_{(4)},\mathop{\beta_2}_{(2)},
\mathop{\beta_1+\beta_2}_{(4)},}$} & 
{\scriptsize$\displaystyle{\{\mathop{\beta_1}_{(4)},\mathop{\beta_2}_{(4)},
\mathop{\beta_1+\beta_2}_{(5)},}$}\\
{\scriptsize$E_6/Spin(10)\cdot U(1)$} & 
{\scriptsize$\displaystyle{\mathop{2\beta_1+\beta_2}_{(2)},
\mathop{2\beta_1}_{(1)},\mathop{2\beta_1+2\beta_2}_{(1)}\}}$} & 
{\scriptsize$\displaystyle{\mathop{2\beta_1+\beta_2}_{(3)}\}}$}\\
\hline
{\scriptsize$\rho_{10}(SU(6)\cdot SU(2))\curvearrowright$} & 
{\scriptsize$\displaystyle{\{\mathop{\beta_1}_{(4)},\mathop{\beta_2}_{(4)},
\mathop{\beta_1+\beta_2}_{(4)},}$} & 
{\scriptsize$\displaystyle{\{\mathop{\beta_1}_{(4)},\mathop{\beta_2}_{(2)},
\mathop{\beta_1+\beta_2}_{(5)},}$}\\
{\scriptsize$E_6/Spin(10)\cdot U(1)$} & 
{\scriptsize$\displaystyle{\mathop{2\beta_1+\beta_2}_{(4)},
\mathop{2\beta_1}_{(1)},\mathop{2\beta_1+2\beta_2}_{(1)}\}}$} & 
{\scriptsize$\displaystyle{\mathop{2\beta_1+\beta_2}_{(1)}\}}$}\\
\hline
{\scriptsize$\rho_{11}(Spin(10)\cdot U(1))\curvearrowright$} & 
{\scriptsize$\displaystyle{\{\mathop{\beta_1}_{(8)},\mathop{2\beta_1}_{(1)},
\mathop{2\beta_1+2\beta_2}_{(1)}\}}$} & 
{\scriptsize$\displaystyle{\{\mathop{\beta_2}_{(6)},\mathop{\beta_1+\beta_2}_{(9)},
\mathop{2\beta_1+\beta_2}_{(5)}\}}$}\\
{\scriptsize$E_6/Spin(10)\cdot U(1)$} & 
{\scriptsize} & {\scriptsize}\\
\hline
{\scriptsize$\rho_{12}(Spin(10)\cdot U(1))\curvearrowright$} & 
{\scriptsize$\displaystyle{\{\mathop{\beta_1}_{(6)},\mathop{\beta_2}_{(1)},
\mathop{\beta_1+\beta_2}_{(6)},}$} & 
{\scriptsize$\displaystyle{\{\mathop{\beta_1}_{(2)},\mathop{\beta_2}_{(5)},
\mathop{\beta_1+\beta_2}_{(3)},}$}\\
{\scriptsize$E_6/Spin(10)\cdot U(1)$} & 
{\scriptsize$\displaystyle{\mathop{2\beta_1+\beta_2}_{(1)}\}}$} & 
{\scriptsize$\displaystyle{\mathop{2\beta_1+\beta_2}_{(4)},
\mathop{2\beta_1}_{(1)},\mathop{2\beta_1+2\beta_2}_{(1)}\}}$}\\
\hline
{\scriptsize$Sp(4)\curvearrowright E_6/F_4$} & 
{\scriptsize$\displaystyle{\{\mathop{\beta_1}_{(4)},\mathop{\beta_2}_{(4)},
\mathop{\beta_1+\beta_2}_{(4)}\}}$} 
& {\scriptsize$\displaystyle{\{\mathop{\beta_1}_{(4)},\mathop{\beta_2}_{(4)},
\mathop{\beta_1+\beta_2}_{(4)}\}}$}\\
\hline
{\scriptsize$\rho_{13}(F_4)\curvearrowright E_6/F_4$} & 
{\scriptsize$\displaystyle{\{\mathop{\beta_1}_{(8)}\}}$} 
& {\scriptsize$\displaystyle{\{\mathop{\beta_2}_{(8)},
\mathop{\beta_1+\beta_2}_{(8)}\}}$}\\
\hline
{\scriptsize$\rho_{14}(SO(4))\curvearrowright$} & 
{\scriptsize$\displaystyle{\{\mathop{\beta_1}_{(1)},
\mathop{3\beta_1+2\beta_2}_{(1)}\}}$} & 
{\scriptsize$\displaystyle{\{\mathop{\beta_2}_{(1)},\mathop{\beta_1+\beta_2}_{(1)},
\mathop{2\beta_1+\beta_2}_{(1)},}$}\\
{\scriptsize$\displaystyle{G_2/SO(4)}$} & 
{\scriptsize} & {\scriptsize$\displaystyle{\mathop{3\beta_1+\beta_2}_{(1)}\}}$}\\
\hline
{\scriptsize$\rho_{15}(SO(4))\curvearrowright$} & 
{\scriptsize$\displaystyle{\{\mathop{\beta_1}_{(1)},
\mathop{3\beta_1+2\beta_2}_{(1)}\}}$} & 
{\scriptsize$\displaystyle{\{\mathop{\beta_2}_{(1)},\mathop{\beta_1+\beta_2}_{(1)},
\mathop{2\beta_1+\beta_2}_{(1)},}$}\\
{\scriptsize$\displaystyle{G_2/SO(4)}$} & 
{\scriptsize} & {\scriptsize$\displaystyle{\mathop{3\beta_1+\beta_2}_{(1)}\}}$}\\
\hline
{\scriptsize$\rho_{16}(G_2)\curvearrowright$} & 
{\scriptsize$\displaystyle{\{\mathop{\beta_1}_{(2)},
\mathop{3\beta_1+2\beta_2}_{(2)}\}}$} & 
{\scriptsize$\displaystyle{\{\mathop{\beta_2}_{(2)},\mathop{\beta_1+\beta_2}_{(2)},
\mathop{2\beta_1+\beta_2}_{(2)},}$}\\
{\scriptsize$\displaystyle{(G_2\times G_2)/G_2}$} & 
{\scriptsize} & {\scriptsize$\displaystyle{\mathop{3\beta_1+\beta_2}_{(2)}\}}$}\\
\hline
{\scriptsize$SU(2)^4\curvearrowright$} & 
{\scriptsize$\displaystyle{\{\mathop{\beta_1}_{(1)},\mathop{\beta_2}_{(1)},
\mathop{\beta_1+\beta_2}_{(1)},}$} & 
{\scriptsize$\displaystyle{\{\mathop{\beta_1}_{(1)},\mathop{\beta_2}_{(1)},
\mathop{\beta_1+\beta_2}_{(1)},}$}\\
{\scriptsize$\displaystyle{(G_2\times G_2)/G_2}$} & 
{\scriptsize$\displaystyle{\mathop{2\beta_1+\beta_2}_{(1)},
\mathop{3\beta_1+\beta_2}_{(1)},\mathop{3\beta_1+2\beta_2}_{(1)}\}}$} & 
{\scriptsize$\displaystyle{\mathop{2\beta_1+\beta_2}_{(1)},
\mathop{3\beta_1+\beta_2}_{(1)},\mathop{3\beta_1+2\beta_2}_{(1)}\}}$}\\
\hline
\end{tabular}$$
}

\vspace{0.15truecm}

\centerline{{\bf Table 2$\,$:$\,$Hermann actions on rank two symmetric spaces (continued${}^2$)}}

\newpage


$$\begin{tabular}{|c|c|}
\hline
{\scriptsize$H\curvearrowright G/K$} & {\scriptsize$H^{\ast}\curvearrowright G^{\ast}/K$}\\
\hline
{\scriptsize$\rho_1(SO(3))\curvearrowright SU(3)/SO(3)$} & 
{\scriptsize$SO_0(1,2)\curvearrowright SL(3,{\Bbb R})/SO(3)$}\\
\hline
{\scriptsize$\rho_2(Sp(3))\curvearrowright SU(6)/Sp(3)$} & 
{\scriptsize$Sp(1,2)\curvearrowright SU^{\ast}(6)/Sp(3)$}\\
\hline
{\scriptsize$\rho_3(SO(4)\times SO(4))\curvearrowright SO(8)/U(4)$} 
& {\scriptsize$SO(4,{\Bbb C})\curvearrowright SO^{\ast}(8)/U(4)$}\\
\hline
{\scriptsize$\rho_4(U(4))\curvearrowright SO(8)/U(4)$} 
& {\scriptsize$U(2,2)\curvearrowright SO^{\ast}(8)/U(4)$}\\
\hline
{\scriptsize$\rho_5(U(5))\curvearrowright SO(10)/U(5)$} 
& {\scriptsize$U(2,3)\curvearrowright SO^{\ast}(10)/U(5)$}\\
\hline
{\scriptsize$\rho_6(SO(5))\curvearrowright(SO(5)\times SO(5))/SO(5)$} 
& {\scriptsize$SO_0(2,3)\curvearrowright SO(5,{\Bbb C})/SO(5)$}\\
\hline
{\scriptsize$\rho_7(U(2))\curvearrowright Sp(2)/U(2)$} 
& {\scriptsize$U(1,1)\curvearrowright Sp(2,{\Bbb R})/U(2)$}\\
\hline
{\scriptsize$\rho_8(Sp(2))\curvearrowright(Sp(2)\times Sp(2))/Sp(2)$} 
& {\scriptsize$Sp(2,{\Bbb R})\curvearrowright Sp(2,{\Bbb C})/Sp(2)$}\\
\hline
{\scriptsize$\rho_9(Sp(2))\curvearrowright(Sp(2)\times Sp(2))/Sp(2)$} 
& {\scriptsize$Sp(1,1)\curvearrowright Sp(2,{\Bbb C})/Sp(2)$}\\
\hline
{\scriptsize$\rho_{10}(SU(6)\cdot SU(2))\curvearrowright E_6/Spin(10)\cdot U(1)$} 
& {\scriptsize$SU(1,5)\cdot SL(2,{\Bbb R})\curvearrowright E_6^{-14}/Spin(10)\cdot U(1)$}\\
\hline
{\scriptsize$\rho_{11}(Spin(10)\cdot U(1))\curvearrowright E_6/Spin(10)\cdot U(1)$} 
& {\scriptsize$SO^{\ast}(10)\cdot U(1)\curvearrowright E_6^{-14}/Spin(10)\cdot U(1)$}\\
\hline
{\scriptsize$\rho_{12}(Spin(10)\cdot U(1))\curvearrowright E_6/Spin(10)\cdot U(1)$} 
& {\scriptsize$SO_0(2,8)\cdot U(1)\curvearrowright E_6^{-14}/Spin(10)\cdot U(1)$}\\
\hline
{\scriptsize$\rho_{13}(F_4)\curvearrowright E_6/F_4$} 
& {\scriptsize$F_4^{-20}\curvearrowright E_6^{-26}/F_4$}\\
\hline
{\scriptsize$\rho_{14}(SO(4))\curvearrowright G_2/SO(4)$} 
& {\scriptsize$SL(2,{\Bbb R})\times SL(2,{\Bbb R})\curvearrowright G_2^2/SO(4)$}\\
\hline
{\scriptsize$\rho_{15}(SO(4))\curvearrowright G_2/SO(4)$} 
& {\scriptsize$\rho^{\ast}_{15}(SO(4))\curvearrowright G_2^2/SO(4)$}\\
\hline
{\scriptsize$\rho_{16}(G_2)\curvearrowright(G_2\times G_2)/G_2$} 
& {\scriptsize$G_2^2\curvearrowright G_2^{\bf C}/G_2$}\\
\hline
\end{tabular}$$

\vspace{0.3truecm}

\centerline{{\bf Table 3$\,$:$\,$The dual actions of $\rho_i$}}

\end{document}